\documentclass{article}

\usepackage{latexsym}
\usepackage{amsmath}
\usepackage{amsfonts}
\usepackage{amssymb}
\usepackage{maple2eACA} 
\usepackage{mapletab}
\def\today{January 30, 2006}

\def\Grobner{Gr\"{o}bner }
\def\Bezier{B\'{e}zier }

\newcommand{\mytaba}{\mbox{\hspace{1.5cm}}}
\newcommand{\mytabb}{\mbox{\hspace{1.2cm}}}

\newenvironment{smallmaplegroup}%
{\small \begin{maplegroup}}{\end{maplegroup} \vskip-0.5ex}

\newcommand{\ed}{\end{document}}

\newcommand{\BP}{\mathbb{P}}

\newcommand{\BR}{\mathbb{R}}

\newcommand{\grad}{\nabla}

\newcommand{\bi}{{\bf i}}
\newcommand{\bj}{{\bf j}}
\newcommand{\bV}{{\bf V}}

\newcommand{\beq}{\begin{equation}}
\newcommand{\eeq}{\end{equation}}

\newtheorem{example}{Example}

\newtheorem{remark}{Remark}

\begin{document}
\title{On the Parallel Lines for Nondegenerate Conics}
\author{Rafa\l \ Ab\l amowicz$^{1)}$
\and
Jane Liu$^{2)}$}
{{\renewcommand{\thefootnote}{\fnsymbol{footnote}}
\footnotetext{\kern-15.3pt AMS Subject Classification: 
\textbf{13P10, 14Q05, 68W30}\par
\vskip 8pt\noindent\vbox{
\noindent\small\sc
1) Department of Mathematics, Box 5054, Tennessee Technological University, Cookeville, TN 38505, USA\par
\noindent
E-mail: {\tt rablamowicz@tntech.edu}, 
URL: {\tt http://math.tntech.edu/rafal/}
\par
\noindent\small\sc
2) Department of Civil and Environmental Engineering,  Box 5015, Tennessee Technological University, Cookeville, TN 38505, USA\par
\noindent
{\tt jliu@tntech.edu},
{\tt http://http://www.tntech.edu/cee/faculty/jane\_liu.html}
}}}}
\def\date{\today} 
\maketitle

\begin{abstract}
Computation of parallel lines (envelopes) to parabolas, ellipses, and hyperbolas is of importance in structure engineering and theory of mechanisms. Homogeneous polynomials that implicitly define parallel lines for the given offset to a conic are found by computing \Grobner bases for an elimination ideal of a suitably defined affine variety. Singularity of the lines is discussed and their singular points are explicitly found as functions of the offset and the parameters of the conic. Critical values of the offset are linked to the maximum curvature of each conic. Application to a finite element analysis is shown.\\

\end{abstract}

\tableofcontents

\section{Introduction}
\label{intro}
The present paper introduces the Reader to a computation of parallel lines (envelopes) to three nondegenerate conics: a parabola, an ellipse, and a hyperbola. In each case, we show how to find a single polynomial $g \in \BR[x,y]$ that implicitly generates a suitable elimination ideal of an affine variety containing the parallel lines for the given offset $r$ and for the given parameters of the conic. This amounts to a computation of a reduced \Grobner basis  \cite{Cox},  \cite{Cox2}, \cite{Singular},  \cite{Wright} for the ideal for an appropriately chosen total order. We use {\tt Maple} \cite{Maple} and {\tt Singular} \cite{Singular} to accomplish computations. We show an application of the results presented here to a finite element analysis. We also make a connection with a design of surface cams in mechanical engineering \cite{Norton}, \cite{Uicker}.  

In each of the next three sections we discuss three cases of parallel lines for the specific values of the offset $0 < r < r_{crit},$ $r = r_{crit},$ and $r > r_{crit}$ for a parabola (Section~\ref{parabola}), an ellipse  (Section~\ref{ellipse}), and a hyperbola (Section~\ref{hyperbola}). Then, in each section, we consider a general case of each conic while imposing essentially no condition on the offset except that $r>0.$ In 
Section~\ref{applications} we show one possible application of the results to a finite element analysis of a composite material of elliptical shape. In Appendices~A and B we collect some of the more cumbersome expressions while in Appendix~C we show in detail how we found the singular points for the parallel lines of the ellipse. 

\section{Parallel Lines to a Parabola}
\label{parabola}

In this section we compute parallel lines to a parabola defined by a polynomial~$f_1:$ 
\beq
f_1 = 4 p y_0 - x_0^2 = 0
\label{eq:f1}
\eeq
where $|p|$ denotes the distance between the focus $F=(0,p)$ and the vertex $V=(0,0).$\footnote{Without any loss of generality, we assume that the parabola has its vertex at the origin and its focus on the $y$-axis.} Polynomial $f_2$ defines a circle of radius (offset) $r$ centered at a point 
$(x_0,y_0)$ on the parabola $f_1,$ 
\beq
f_2=(y-y_0)^2+(x-x_0)^2-r^2 = 0
\label{eq:f2}
\eeq
while a polynomial $f_3$ 
\beq
f_3 = 2 x p - 2 x_0 p+ x_0 y-x_0 y_0 = 0
\label{eq:f3}
\eeq
gives a condition that a point $P(x,y)$ on the circle $f_2$ lies on a line perpendicular to the parabola $f_1$ and passing through the point $(x_0,y_0)$ on the parabola. The family of such circles is called the {\em family of curves\/} determined by the polynomial $f_1$ (see \cite{Cox}). There are of course two such points for any given point $(x_0,y_0):$ one on each side of the parabola. All these points $P$ belong to an affine variety 
$\bV=\bV(f_1,f_2,f_3)$ -- the {\em envelope\/} of the family of circles -- and define two parallel lines at the distance $r$ from the parabola. In what follows we will compute a reduced \Grobner basis for $\bV$ and by eliminating variables $x_0,y_0$ we will obtain one polynomial $g$ in variables $x$ and 
$y$ alone (also $p,r$ in a general case) that defines both parallel lines on each side of the parabola for the given value of the offset $r.$ Thus, in general, we have $f_1,f_2,f_3 \in \BR[x_0,y_0,x,y,p,r].$ However, first we illustrate computation for specific values of the parameters $p$ and $r,$ and later we show computation for unspecified parameters.

\subsection{Parallel lines for specific parameters $p$ and $r$}
\label{parabola1}

The parallel lines constitute a subvariety of the affine variety $\bV = \bV(I)$ where $I = \langle f_1, f_2, f_3\rangle.$ In order to find the polynomial $g$ for the second elimination ideal $I_2 = I \cap  \BR[x,y],$ we first find a \Grobner basis $G$ for $\bV$ for the lexicographic order 
$y_0 > x_0 >x > y$ that contains sixteen polynomials.\footnote{Computations have been performed with {\tt Groebner} package in {\tt Maple 8} \cite{Maple} as well as with a small package written by the Authors \cite{RJGrobner}, and the results have been confirmed with {\tt Singular} \cite{Singular} using \cite{Singularlink}.} Then, we reduce $G$ to a minimal \Grobner basis $G_m$ that contains six (or eight, depending on the case) polynomials and, further, to a unique reduced \Grobner basis $G_r$ with three (or five, depending on the case) polynomials only. Finally, we eliminate from $G_r$ all polynomials that contain variables $y_0$ and $x_0,$ and keep the single polynomial $g \in \BR[x,y]$ that gives the parallel lines and generates the second elimination ideal $I_2 = \langle g \rangle.$

We are also interested when the parallel lines form a smooth (non singular) variety $\bV(g)$ or a $C^{\infty}$ $1$-manifold or a {\it hypersurface defined by $g.$} Since $g$ can be viewed as a $C^{\infty}$ function $\BR^2 \rightarrow \BR$ and $\bV(g) = \{(x,y) \in \BR^2 \,|\, g(x,y) = 0 \},$ it is known that $\bV(g)$ is a smooth manifold at those points $p \in \bV(g)$ where $(\grad g)(p) \neq 0.$ \cite[Theorem 2.6]{Millman} We will find singular points of $\bV(g)$ by finding points $p \in \bV(g)$ where $(\grad g)(p) = 0.$  

Thus, $\bV(g)$ can be thought of as a union of points where the gradient is zero and those points where the gradient is not zero. Connected components of $\bV(g)$ where the gradient is not zero form the $1$-dimensional manifold of the parallel lines. However, we are also interested in describing a discrete finite set $S$ of singular points where the gradient is zero. For a formal definition of a singular point of an affine variety 
$\bV(g)$ see \cite[page 136]{Cox}.

\begin{figure}[htb]
\centerline{\scalebox{1.20}{\includegraphics{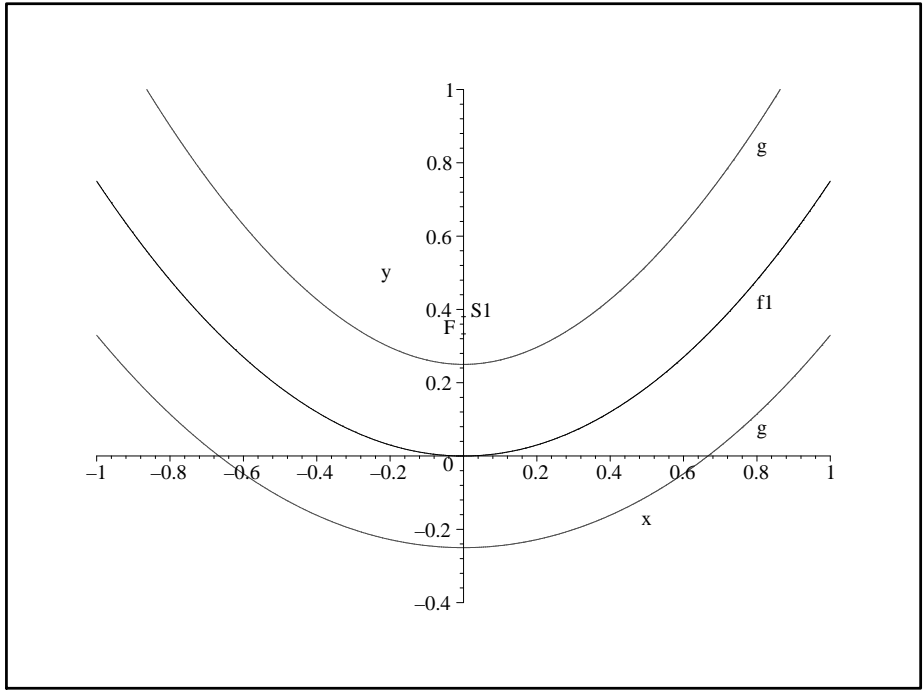}}}
\caption{Parabola with parallel lines $g,$ focus $F$ and one virtual singular point $S_1$ for $p=\frac{1}{3}$ and $r=\frac{1}{4}$}
\label{fig:parabola1}
\end{figure}

\begin{example}[$r < r_{crit}$]
\label{example1}
Let $p=\frac{1}{3}$ and $r=\frac{1}{4}.$ Then,
$$
\bV=\bV(4 y_0-3 x_0^2, 16 y^2-32 y y_0+16 y_0^2+16 x^2-32 x x_0+16 x_0^2-1, 2 x-2 x_0+3 x_0 y-3 x_0 y_0)
$$
We find that the reduced \Grobner basis $G_r$ contains only three polynomials including 
\begin{multline}
g = 331776 x^6-42192 x^2+48400 y^2-448512 y^3-84480 x^2 y+589824 y^4+28032 y \\[0.5ex]
-25344 x^4+1138176 x^2 y^2-884736 x^2 y^3-1105920 x^4 y-5329+331776 x^4 y^2
\label{eq:gp1}
\end{multline}
whose gradient is
\begin{setlength}{\multlinegap}{0pt}
\begin{multline}
\grad g = (-879 x-1760 y x-1056 x^3+23712 x y^2-46080 x^3 y-18432 y^3 x+13824 x^3 y^2+20736 x^5) \bi \\[0.5ex]
         +(3025 y-2640 x^2-42048 y^2+876+71136 x^2 y-34560 x^4-82944 x^2 y^2+73728 y^3+20736 x^4 y) \bj\,.
\end{multline}
\end{setlength}

\noindent
Thus, the only real point $S_1$ on $\bV(g)$ where $(\grad g)(S_1) = 0$ is $S_1=(0,\frac{73}{192}).$ We will refer to this point as a {\rm virtual singular point} for the following reason:  $\bV(g)$ is a disjoint union of $\{S_1\}$ and a non-singular smooth subvariety $\Pi$ where the gradient $\grad g$ does not vanish. The subvariety $\Pi$ constitutes the set of two parallel lines that we have been seeking. We show all in 
Figure~\ref{fig:parabola1}.
\end{example}

We have verified by direct computation that as the value of the offset $r$ increases towards some critical value $r_{crit},$ the virtual singular point approaches the smooth component of the variety $\bV(g),$ that is, it approaches the parallel lines $\Pi.$\footnote{In the next subsection we will derive a formula for the critical value $r_{crit}$ of $r$ for a general parabola.} As the next example shows, when $r = r_{crit},$ then the virtual singular point is actually located on the smooth component $\Pi$ of the variety $\bV(g).$

\begin{figure}[htb]
\centerline{\scalebox{1.20}{\includegraphics{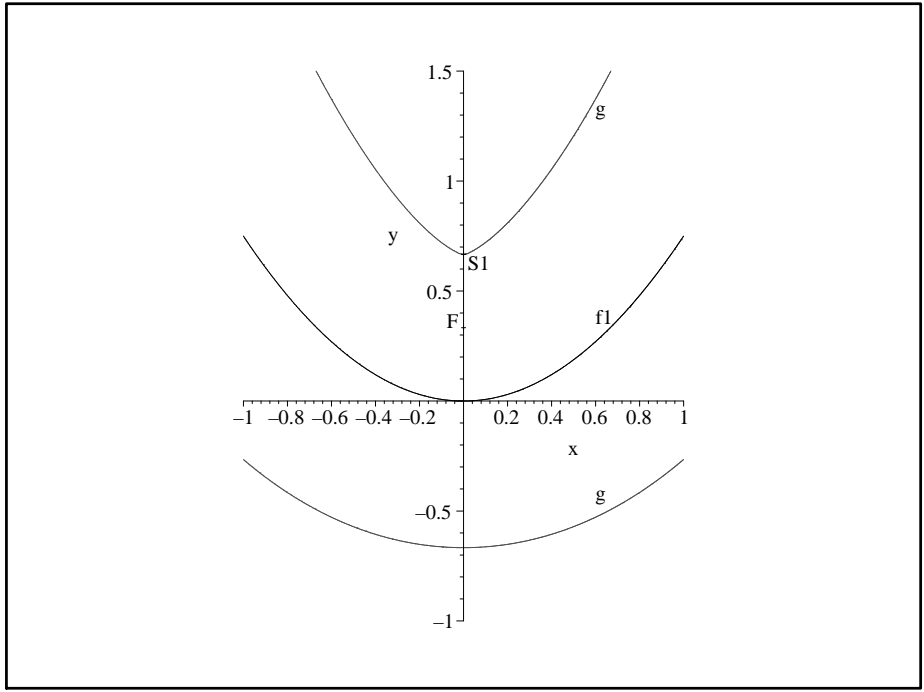}}}
\caption{Parabola with parallel lines $g,$ focus $F$ and one singular point $S_1$ for $p=\frac{1}{3}$ and $r=\frac{2}{3}$}
\label{fig:parabola2}
\end{figure}

\begin{example}[$r = r_{crit}$]
\label{example2}
Let $p=\frac{1}{3}$ and $r=\frac{2}{3}.$ Then,
$$
\bV=\bV(4 y_0-3 x_0^2, 9 y^2-18 y y_0+9 y_0^2+9 x^2-18 x x_0+9 x_0^2-4, 2 x-2 x_0+3 x_0 y-3 x_0 y_0)
$$
In this case, we have found that the reduced \Grobner basis $G_r$ contains five polynomials including
\begin{multline}
g = 768 y-288 x^2-1728 y^3+1944 x^2 y^2-891 x^4-2430 x^4 y \\[0.5ex]
    -1944 x^2 y^3+1296 y^4+729 x^4 y^2+729 x^6-256
\label{eq:gp2}
\end{multline}
whose gradient is
\begin{setlength}{\multlinegap}{0pt}
\begin{multline}
\grad g = (-32 x+216 x y^2-198 x^3-540 x^3 y-216 y^3 x+162 x^3 y^2+243 x^5) \bi \\[0.5ex]
         +(128-864 y^2+648 x^2 y-405 x^4-972 x^2 y^2+864 y^3+243 x^4 y) \bj\,.
\end{multline}
\end{setlength}

\noindent
Thus, like in Example~\ref{example1}, there is only one real point $S_1$ on $\bV(g)$ where $(\grad g)(S_1) = 0:$ it is $S_1=(0,\frac{2}{3}).$ However, unlike in the first example, this time point $S_1$ is not disjoint from the rest of the variety $\bV(g)$ as it can be seen in 
Figure~\ref{fig:parabola2}. We will refer to this point as a {\rm singular point} of the variety $\bV(g).$ The latter can still be thought of as union of $\{S_1\}$ and the non-singular smooth subvariety $\Pi$ where the gradient $\grad g$ does not vanish. We will say that in this case the parallel lines possess one singular point.
\end{example}

Probably the most interesting case occurs when the value of the offset $r > r_{crit}$ as shown in the next example. In this case, the single singular point of the variety $\bV(g)$ splits into three singular points.

\begin{figure}[htb]
\centerline{\scalebox{1.40}{\includegraphics{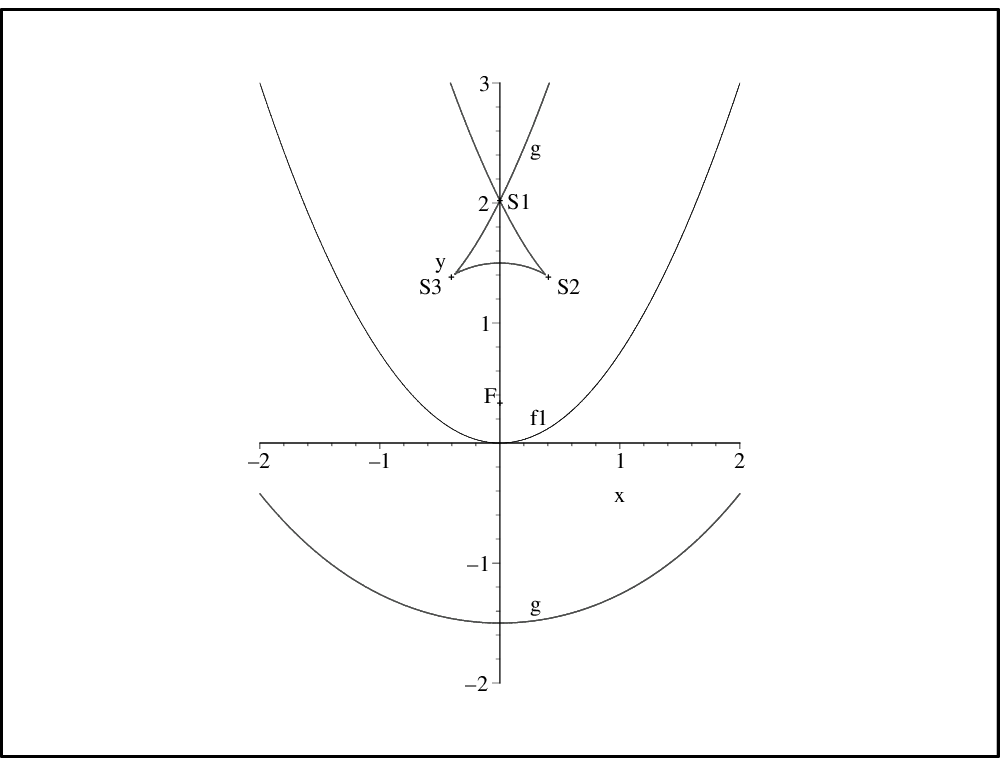}}}
\caption{Parabola with parallel lines $g,$ focus $F$ and three singular points $S_1,S_2,S_3$ for $p=\frac{1}{3}$ and $r=\frac32$}
\label{fig:parabola3}
\end{figure}

\begin{example}[$r > r_{crit}$]
\label{example3}
Let $p=\frac{1}{3}$ and $r=\frac{3}{2}.$ Then,
$$
\bV=\bV(4 y_0-3 x_0^2, 4 y^2-8 y y_0+4 y_0^2+4 x^2-8 x x_0+4 x_0^2-9, 2 x-2 x_0+3 x_0 y-3 x_0 y_0)
$$
In this case, we have found that the reduced \Grobner basis $G_r$ again contains only three polynomials including
\begin{multline}
g = 83808 y+52812 x^2+16900 y^2-37248 y^3-4896 x^2 y^2-34416 x^4-17280 x^4 y\\[0.5ex]
    -13824 x^2 y^3+9216 y^4+5184 x^4 y^2+5184 x^6-84681+6240 x^2 y 
\label{eq:gp3} 
\end{multline}
whose gradient is
\begin{setlength}{\multlinegap}{0pt}
\begin{multline}
\grad g = (4401 x-408 x y^2-5736 x^3-2880 x^3 y-1152 y^3 x+864 x^3 y^2+1296 x^5+520 y x) \bi \\[0.5ex]
         +(10476+4225 y-13968 y^2-1224 x^2 y-2160 x^4-5184 x^2 y^2+4608 y^3+1296 x^4 y+780 x^2) \bj\,.
\end{multline}
\end{setlength}

\noindent
Thus, unlike in the above two examples, this time there are three singular points: $S_1=(0,\frac{97}{48})$ located on the $y$-axis (the symmetry axis of the parabola), and two additional points $S_2$ and $S_3$
\beq
S_2 =
\left(
{\displaystyle \frac {\sqrt{65 + 36 \cdot 12^{\frac13} - 27 \cdot 12^{\frac23}}}{6}} , \,
{\displaystyle  \frac {3\cdot 12^{\frac13}}{4}}  - {\displaystyle  \frac {1}{3}} 
\right),\\[0.5ex] 
\label{eq:Sp2}
\eeq
\beq
S_3 =
\left( 
 - {\displaystyle \frac {\sqrt{65 + 36\cdot 12^{\frac13} - 27\cdot 12^{\frac23}}}{6}} , \,
   {\displaystyle \frac {3\cdot 12^{\frac13}}{4}}  - {\displaystyle \frac {1}{3}} 
\right) 
\label{eq:Sp3}
\eeq 
located symmetrically on each side of the $y$-axis as shown in Figure~\ref{fig:parabola3}. We will say that in this case the parallel lines possess three singular points. A similar figure can be found in \cite[page 141]{Cox} 
\end{example}

\subsection{Parallel lines for arbitrary parameters $p$ and $r$}
\label{parabola2}

In this section we will find a general polynomial $g$ that will generate the second elimination ideal $I_2 = I \cap \BR[x,y,p,r]$ where 
$I = \langle f_1,f_2,f_3 \rangle$ and the generating polynomials are defined in equations (\ref{eq:f1}), (\ref{eq:f2}), 
and~(\ref{eq:f3}). We will also find a formula for the critical value $r_{crit}$ of $r$ that determines whether $\bV(I_2)$ has one or three singular points, and, we will find general formulas for the coordinates of the said singular points for arbitrary values of the parameters $p$ and $r.$

We begin with a computation of a \Grobner basis for the ideal $I$ for the lexicographic order $y_0 > x_0 > x > y > r > p.$
$$
I = \langle 4 p y_0-x_0^2, y^2-2 y y_0+y_0^2+x^2-2 x x_0+x_0^2-r^2, 2 x p-2 x_0 p+x_0 y-x_0 y_0 \rangle \subset\BR[y_0, x_0, x, y, r, p]
$$

Using {\tt Maple} and {\tt Singular} \cite{Singularlink},  \cite{Singular}  we found that a reduced \Grobner basis for $I$ consists of fourteen homogeneous polynomials of degrees: 
$7$ (two polynomials),
$6$ (two polynomials),
$5$ (two polynomials),
$4$ (two polynomials),
$3$ (three polynomials),
and $2$ (three polynomials). These polynomials are displayed in Appendix~A. 

Thus, in this general case we observe that $I_2 = \langle g \rangle$ where $g = g_1/p$ is the following homogeneous polynomial in $\BR[x,y,r,p]$ of degree six: 
\begin{setlength}{\multlinegap}{0pt}
\begin{multline}
g = -2 p r^2 y x^2+8 p r^2 y^3+8 p^2 r^2 y^2-32 y p^3 r^2+16 p^4 r^2-16 y^4 p^2+ 32 y^3 p^3\\[0.5ex]
    -16 p^4 y^2+3 r^2 x^4+8 p^2 r^4+20 p^2 r^2 x^2-y^2 x^4 +10 y p x^4 -x^6-x^4 p^2+8 p y^3 x^2\\[0.5ex]
    -32 x^2 y^2 p^2+8 x^2 y p^3-3 r^4 x^2+2 r^2 x^2 y^2 +r^6-r^4 y^2-8 p r^4 y
\label{eq:gparabola}  
\end{multline}
\end{setlength}

\noindent
We want to find now the singular points on the variety $\bV=\bV(g)$ where the gradient $\grad g$ is zero. Thus, $\bV$ can be thought of as union of points where the gradient is zero and those points where the gradient is not zero. Connected components of $\bV$ where the gradient is not zero constitute the $1$-dimensional manifold that we encountered above. We are interested in the finite set $S$ of the singular points of the variety. We will show next that the critical value of $r$ is $r_{crit} = |2p| = 2 |p|.$ Therefore, we compute $\grad g= \grad g(x,y,r,p):$
\begin{setlength}{\multlinegap}{0pt}
\begin{multline}
\grad g = \\[0.5ex]
<-x (2 y p r^2-6 x^2 r^2-20 p^2 r^2+2 x^2 y^2-20 y p x^2+3 x^4+2 x^2 p^2-8 p y^3+32 y^2 p^2-8 y p^3+3 r^4-2 r^2 y^2),\\[0.5ex]
 -x^2 p r^2+12 y^2 p r^2+8 y p^2 r^2-16 p^3 r^2-32 y^3 p^2+48 p^3 y^2-16 y p^4-y x^4+5 x^4 p+12 x^2 p y^2-32 x^2 p^2 y\\[0.5ex]
 +4 x^2 p^3+2 r^2 y x^2-r^4 y-4 p r^4 ,\\[0.5ex]
 r (-2 y p x^2+8 p y^3+8 y^2 p^2-32 y p^3+16 p^4+3 x^4+16 p^2 r^2+20 x^2 p^2-6 x^2 r^2+2 x^2 y^2+3 r^4\\[0.5ex]
 -2 r^2 y^2-16 y p r^2),\\[0.5ex]
 -r^2 y x^2+4 r^2 y^3+8 y^2 p r^2-48 y p^2 r^2+32 p^3 r^2-16 y^4 p+48 y^3 p^2-32 p^3 y^2+8 p r^4+20 x^2 p r^2\\[0.5ex]
 +5 y x^4-x^4 p +4 y^3 x^2-32 x^2 p y^2+12 x^2 p^2 y-4 r^4 y
>
\end{multline}
\end{setlength}

\noindent
In order to find solutions to the vector equation $\grad g = 0,$ we need to solve in fact a resulting system of four polynomial equations in variables 
$x,y,p,r.$ The system could be simplified by assuming $r>0,$ yet to efficiently analyze its solutions one can again employ the \Grobner basis approach and find a reduced \Grobner basis for the ideal $I_S$ generated by the four equations. Since we want to make sure that the solutions belong to $\bV(g)$, we add equation $g=0$ to our system.  Let $F$ be a list consisting of the five polynomial equations, that is, $F = [\grad g=0,g=0  ]$ and let 
$I_S = \langle \grad g,g \rangle.$ Thus, we find that the affine variety $\bV(I_S)$ of singular points is generated by a reduced \Grobner basis consisting of eleven homogeneous polynomials $s_i, i=1,\ldots,11,$ in $y, x, r,$ and $p$ with total degrees varying from $5$ through $17.$ The third polynomial $s_3$ in the basis has the following form:
\begin{multline}
s_3 = x r (640 p^6 y-320 p^7+3 r^4 y x^2-3 r^2 y x^4-12 p r^2 x^4+48 x^2 p^3 r^2-r^6 y+14 r^6 p \\[0.5ex]
           +84 p^3 r^4-480 p^5 r^2+y x^6+13 x^6 p-24 x^4 p^3-384 x^2 p^5-96 y p^2 r^2 x^2-15 p r^4 x^2 \\[0.5ex]
           -60 y p^2 r^4+96 y p^4 r^2-60 y x^4 p^2-96 x^2 p^4 y)
\end{multline}

\noindent
It appears that since $r$ and $p$ are non-zero, one possible solution to the equation $s_3 = 0$ is $x = 0.$ Let's then see if the remaining polynomials in $F$ can be solved when $x = 0.$ Upon substituting $x=0$ into the list~$F,$ one obtains
\begin{multline}
[0, -(-4 p y+r^2+4 p^2) (-8 p y^2+y r^2+4 p r^2+4 y p^2), r (-4 p y+r^2+4 p^2) (-2 y^2+3 r^2-4 p y+4 p^2),\\[0.5ex]
 4 (-y+2 p) (r-y) (r+y) (-4 p y+r^2+4 p^2), (r-y) (r+y) (-4 p y+r^2+4 p^2)^2]
\label{eq:F0}
\end{multline}

\noindent
Since $r>0,$ one possible solution to (\ref{eq:F0}) is $r=y$ which gives the following system:
\beq
F_1 = [0, -(r-2 p)^4 r, (r-2 p)^4 r, 0, 0]
\label{eq:F1}
\eeq
This leads to the critical value of $r_{crit} = 2p.$ Of course, this makes sense only when $p > 0,$ that is, in general we have $r_{crit} = |2p| = 2|p|.$ When $p < 0,$ we can still get the same critical value of $r$ when we select $r = - y$ in (\ref{eq:F0}) which leads to
\beq
F_2 = [0, (r+2 p)^4 r, (r+2 p)^4 r, 0, 0].
\eeq
Thus, the critical value of $r$ in this case is  $r_{crit} = -2p.$ Combining both cases we get $r_{crit} = |2p| = 2|p|.$

Another solution of the system (\ref{eq:F0}) can be found by setting $y = \dfrac{4 p^2+r^2}{4p} = p + \dfrac{r^2}{4p}\; (p \neq 0).$ This gives the singular point 
\beq
S_1 = \left(0, \frac{r^2+4 p^2)}{4 p}\right)= \left(0,p + \dfrac{r^2}{4p}\right) \quad \mbox{when} \quad p \neq 0   \quad \mbox{and} \quad   r>0
\label{eq:Sp1}
\eeq
found in Examples~\ref{example1}, \ref{example2} and \ref{example3} above. When $0 < r < r_{crit} = 2|p|,$ then this is the only singular point $S_1$ as shown in Figure~\ref{fig:parabola1} that is not located on the continuous part of $\bV(g):$ It is located at the distance $\dfrac{r^2}{4|p|}$ above the focus (when $p>0)$ and below the focus (when $p<0)$ on the symmetry line of the parabola (the $y$-axis in our case). Then, the parallel line is smooth at the points $(0, r)$ and $(0, -r).$ These points are located on the $y$ axis above and below the vertex of the parabola.  As $r$ approaches 
$r_{crit}$ from the left, then $S_1$ approaches one of these points: It approaches the point located on the parallel line on the concave side of the parabola (see Figure~\ref{fig:parabola1}). When $r = r_{crit}$ then the singular point $S_1$ is located on the parallel line on the concave side of the parabola (see Figure~\ref{fig:parabola2}). 

Now we look for those singular points on $\bV(g)$ for which $x \neq 0.$ These points will be present only when $r > r_{crit}.$ Thus, we now proceed to solve the system of five polynomial equations contained in the list $F,$ or, equivalently, the system of eleven polynomial equations $s_i = 0, i=1,\ldots,11,$ in the \Grobner basis for $I_S$ under the condition that $x \neq 0.$ Using approach described in more details in Appendix~C for the ellipse, in the case when $r > r_{crit}$ one obtains the following two additional singular points
\beq
S_2 = \left(\alpha,\frac{\beta}{\gamma}\right) \quad \mbox{and} \quad S_3 = \left(-\alpha,\frac{\beta}{\gamma}\right)
       \quad \mbox{when} \quad p \neq 0  \quad \mbox{and} \quad  r>r_{crit}=2|p|.
\label{eq:Sp23}
\eeq
The expressions $\alpha, \beta$ and $\gamma$ are as follows:
\beq
\alpha = \dfrac{\sqrt{(p r^2)^{\frac23} r^2+6 r^2 p^2 2^{\frac13}-3 (p r^2)^{\frac13} p r^2 2^{\frac23}-4 (p r^2)^{\frac23} p^2}}{\sqrt[3]{p r^2}},
\eeq
\begin{multline}
\beta = p r^2 (22 r^6 (p r^2)^{\frac23}+1452 r^6 p^2 2^{\frac13}+7456 p^3 r^4 2^{\frac23} (p r^2)^{\frac13}-6560 p^2 r^4 (p r^2)^{\frac23}\\[0.5ex]
-15488 p^4 r^4 2^{\frac13}-39680 p^5 r^2 2^{\frac23} (p r^2)^{\frac13}+37600 p^4 r^2 (p r^2)^{\frac23}+7936 p^8 2^{\frac13}\\[0.5ex]
-15872 p^7 2^{\frac23} (p r^2)^{\frac13}+45760 p^6 r^2 2^{\frac13}-33 r^8 2^{\frac13}-5376 p^6 (p r^2)^{\frac23}),
\end{multline}
\begin{multline}
\gamma = 2 (p r^2)^{\frac13} (8640 p^6 r^2 (p r^2)^{\frac13}-3968 p^7 2^{\frac13} (p r^2)^{\frac23}+1984 p^7 r^2 2^{\frac23}\\[0.5ex]
-3920 r^4 p^4 (p r^2)^{\frac13}-9920 p^5 r^2 2^{\frac13} (p r^2)^{\frac23}+4960 p^5 r^4 2^{\frac23}+484 r^6 p^2 (p r^2)^{\frac13}\\[0.5ex]
+1864 p^3 r^4 2^{\frac13} (p r^2)^{\frac23}-932 p^3 r^6 2^{\frac23}-11 r^8 (p r^2)^{\frac13}).
\end{multline}

\noindent
Thus, in summary, we have these three cases:
\begin{itemize}
\item[(i)] One (virtual) singular point $S_1 = (0,p+\frac{r^2}{4p})$ when $r < r_{crit}=2|p|.$ In particular, when $p=\frac13$ and $r = \frac14,$ we get 
from (\ref{eq:Sp1}) that $S_1=(0,\frac{73}{192}):$ That is the only (discrete) singular point of $\bV(g)$ found in Example~\ref{example1}. Formula 
(\ref{eq:gparabola}) then gives
\begin{multline*}
g = \frac{5329}{331776}-\frac{3025}{20736} y^2-\frac{73}{864}y+\frac{55}{216} y x^2+\frac{73}{54} y^3+\frac{293}{2304} x^2\\[0.5ex]
    + \frac{11}{144} x^4-\frac{247}{72} x^2 y^2 -\frac{16}{9} y^4+\frac{10}{3} y x^4+\frac{8}{3} y^3 x^2-x^6-y^2 x^4
\end{multline*}
which is, up to the form and sign, the same as polynomial (\ref{eq:gp1}).

\item[(ii)] One (true) singular point $S_1 = (0,p+\frac{r^2}{4p})$ when $r = r_{crit}=2|p|.$ In particular, when $p=\frac13$ and $r=\frac23,$ we get
from (\ref{eq:Sp1}) that $S_1=(0,\frac23):$ That is the only singular point of located on the continuous part of $\bV(g)$ found in 
Example~\ref{example2}. Formula (\ref{eq:gparabola}) then gives
\begin{gather*}
g = \frac{256}{729}-\frac{256}{243} y+\frac{64}{27} y^3+\frac{32}{81} x^2+\frac{11}{9} x^4-\frac{8}{3} x^2 y^2-
    \frac{16}{9} y^4 + \frac{10}{3} y x^4+\frac{8}{3} y^3 x^2-x^6-y^2 x^4
\end{gather*}
which is, up to the form and sign, the same as polynomial (\ref{eq:gp2}).

\item[(iii)] Three singular points: $S_1 = (0,p+\frac{r^2}{4p}),\, S_2=(\alpha,\frac{\beta}{\gamma}),$ and $S_3=(-\alpha,\frac{\beta}{\gamma})$ 
when $r > r_{crit}=2|p|.$  In particular, when $p=\frac13$ and $r = \frac32,$ we get from (\ref{eq:Sp1}) that  $S_1=(0,\frac{97}{48})$ while from 
(\ref{eq:Sp23}) we get the two additional singular points $S_2$ and $S_3$ displayed in (\ref{eq:Sp2}) and (\ref{eq:Sp3}) found in Example~3. Formula 
(\ref{eq:gparabola}) then gives
\begin{multline*}
g = \frac{9409}{576}-\frac{4225}{1296} y^2-\frac{97}{6} y-\frac{65}{54} y x^2+\frac{194}{27} y^3-\frac{163}{16} x^2+\frac{239}{36} x^4\\[0.5ex]
+\frac{17}{18} x^2 y^2-\frac{16}{9} y^4+\frac{10}{3} y x^4+\frac{8}{3} y^3 x^2-x^6-y^2 x^4
\end{multline*}
which is, up to the form and sign, the same as polynomial (\ref{eq:gp3}).
\end{itemize}

As a final comment, let's notice that for the parabola $4p(y-k) = (x-h)^2$ the critical value of the offset 
$$
r_{crit} = 2|p| = \frac{1}{\kappa_{max}} = \rho_{min}
$$ 
is the reciprocal of $\kappa_{max},$ the maximum curvature of the parabola at its vertex, or, equivalently, it equals $\rho_{min},$ the minimum radius of the osculating circle.\footnote{Recall that when $\beta(t)$ is  a regular curve in $\BR^3$ then $\kappa =|\Dot{\beta} \times \Ddot{\beta}|/|\Dot{\beta}|^3.$ 
\cite[page 46]{Millman}.} Furthermore, the quantity  $2|p|$ is known as the {\em semi-latus rectum} and represents the distance from the parabola focus to the parabola measured along a line parallel to the parabola directrix. \cite{wiki}

\section{Parallel Lines to an Ellipse}
\label{ellipse}

In this section we compute various parallel lines to an ellipse defined by a polynomial
\beq
f_1 = b^2 y_0^2+a^2 x_0^2-a^2 b^2 = 0
\label{eq:fe1}
\eeq
where we assume that $a>b$ and the center of the ellipse is at the origin. Polynomial $f_2$ defines a circle of radius (offset) $r$ centered at a point $P=(x_0,y_0)$ on the ellipse,
\beq
f_2 = (y-y_0)^2+(x-x_0)^2-r^2 = 0  
\label{eq:fe2}
\eeq
while polynomial $f_3$
\beq
f_3 = b^2 (x-x_0) y_0 - a^2 (y-y_0) x_0 = 0 
\label{eq:fe3}
\eeq
gives a condition that a point $P(x,y)$ on the circle $f_2$ lies on a line perpendicular to the ellipse $f_1$ and passing through the point 
$(x_0,y_0)$ on the ellipse. All these points $P$ belong to the affine variety $\bV = \bV(f_1,f_2,f_3)$ -- the envelope -- and, like for the parabola, define two parallel lines at the distance $r$ from the ellipse. We proceed like in Section~\ref{parabola}: We will compute a reduced \Grobner basis for $\bV$ and by eliminating variables $x_0,y_0$ we will obtain one polynomial $g$ in variables $x$ and $y$ alone (also $a,b,r$ in a general case) that defines both parallel lines on each side of the ellipse for the given value of the offset $r.$ Thus, in general, we have 
$f_1,f_2,f_3 \in \BR[y_0,x_0,x,y,a,b,r].$ However, first we illustrate computations for specific values of the parameters $a,b$ and $r,$ and later we show computation for unspecified parameters.

\subsection{Parallel lines for specific parameters $a,b,r$}
\label{ellipse1}

The parallel lines constitute a subvariety of the affine variety $\bV = \bV(I)$ where, like before, $I \langle f_1, f_2, f_3\rangle.$ In order to find the polynomial $g$ for the second elimination ideal $I_2 = I \cap \BR[x,y]$, we first find a \Grobner basis~$G$ for $\bV$ for the lexicographic order 
$y_0> x_0 > x >y$ and reduce it eventually to a unique reduced basis~$G_r$ of nine (or eleven, depending on the case) polynomials only. Finally, we eliminate from that basis all polynomials that contain variables $y_0$ and $x_0,$ and keep the single polynomial $g \in \BR[x,y]$ that gives the parallel lines and generates the second elimination ideal $I_2= \langle g \rangle.$ We will also proceed to find all singular points $p$ of $\bV(g)$ where $(\grad g)(p) = 0.$ Thus, let's consider three cases when $0<r<r_{crit},\,r = r_{crit},$ and $r>r_{crit}.$ In the next section we will consider a general case for arbitrary values of $a,b\; (a>b>0)$ and $r\; (r>0).$

\begin{figure}[htb]
\centerline{\scalebox{1.30}{\includegraphics{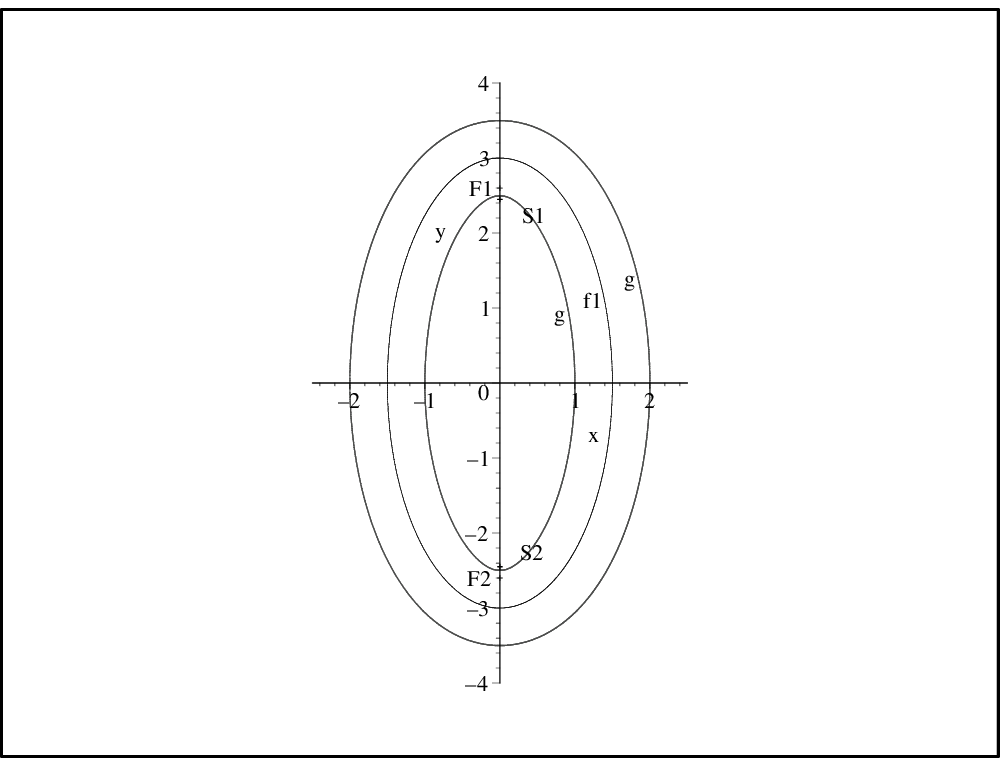}}}
\caption{Ellipse with parallel lines $g,$ foci $F_1,F_2$ and two singular points $S_1,S_2$ for $a=3,b=\frac32$ and $r=\frac12$}
\label{fig:ellipse1}
\end{figure}

\begin{example}[$r < r_{crit}$]
\label{example4}
Let $a=3,$ $b=\frac32,$ and $r=\frac12.$ Then,
$$
\bV=\bV(y_0^2+4 x_0^2-9, 4 y^2-8 y y_0+4 y_0^2+4 x^2-8 x x_0+4 x_0^2-1, y_0 x+3 y_0 x_0-4 x_0 y)
$$
We find that the reduced \Grobner basis $G_r$ contains nine polynomials including
\begin{setlength}{\multlinegap}{0pt}
\begin{multline}
g = 256 x^8+2080 x^6-41685 x^2+44100-25356 y^2+5353 y^4-4751 x^4-488 y^6-3360 y^4 x^2 \\[0.5ex]
    -4680 y^2 x^4+528 y^4 x^4+640 x^6 y^2+160 y^6 x^2+16 y^8+21410 x^2 y^2
\label{eq:ge4}
\end{multline}
\end{setlength}
whose gradient is
\begin{align}
\begin{split}
\grad g &= x (1024 x^6+6240 x^4-41685-9502 x^2-3360 y^4-9360 x^2 y^2+1056 y^4 x^2\\[0.5ex]
        &\hspace*{2in}  +1920 y^2 x^4+160 y^6+21410 y^2) \bi \\[0.5ex]
        &\hspace*{0.18in} + y (-12678+5353 y^2-732 y^4-3360 x^2 y^2-2340 x^4+528 y^2 x^4+320 x^6 \\[0.5ex]
        &\hspace*{2in}  +240 y^4 x^2+32 y^6+10705 x^2) \bj\,.
\end{split}
\end{align}
\noindent
There are exactly two solutions to the equation $(\grad g)(p) = 0$ on $\bV(g):$ 
$$
S_1 = (0, \sqrt{6}) \quad \mbox{and} \quad  S_2 = (0, -\sqrt{6})
$$ 
We will refer to these points as {\rm virtual singular points} for the following reason:  $\bV(g)$ is a disjoint union of $\{S_1,S_2\}$ and a non-singular smooth subvariety $\Pi$ where the gradient $\grad g$ does not vanish. The subvariety $\Pi$ constitutes the set of two parallel lines that we have been seeking. Note that the foci $F_1 = (0,\frac{3\sqrt{3}}{2})$ and $F_2 = (0,-\frac{3\sqrt{3}}{2})$ are different from the singular points. We show all in Figure~\ref{fig:ellipse1}.
\end{example}

We have verified by direct computation that as the value of the offset $r$ increases towards some critical value $r_{crit},$ the virtual singular points approach the smooth component of the variety $\bV(g),$ that is, they approach the parallel lines $\Pi.$\footnote{In the next subsection we will derive a formula for the critical value $r_{crit}$ of $r$ for a general ellipse.} As the next example shows, when $r = r_{crit},$ then the virtual singular points are actually located on the smooth component $\Pi$ of the variety $\bV(g).$

\begin{figure}[htb]
\centerline{\scalebox{1.30}{\includegraphics{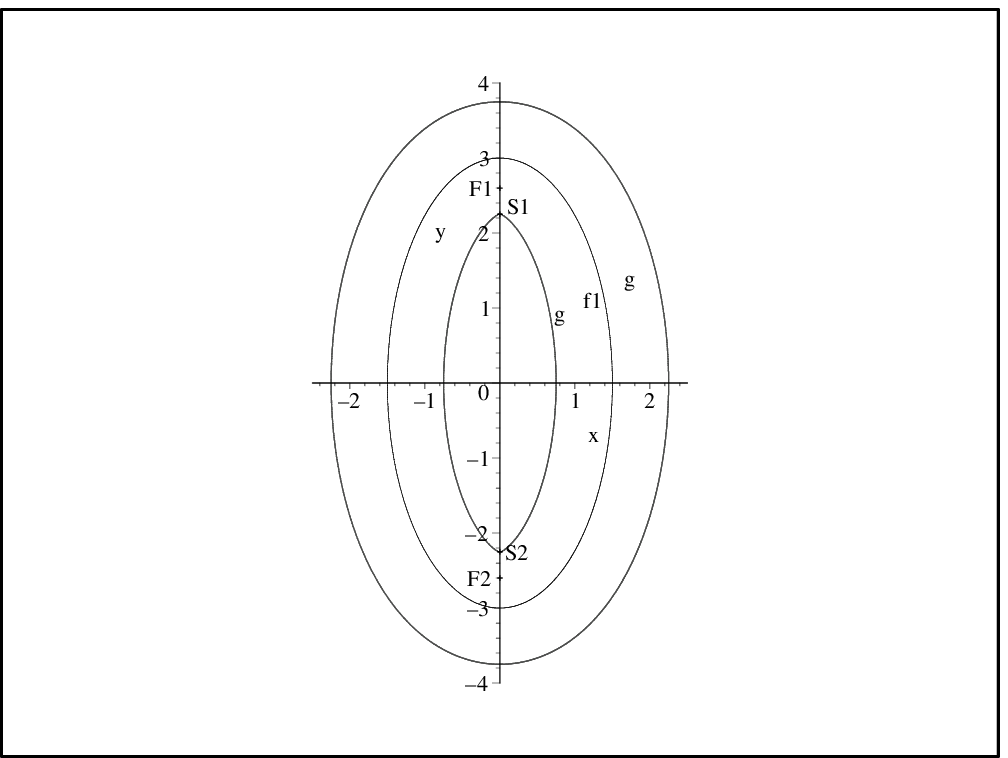}}}
\caption{Ellipse with parallel lines $g,$ foci $F_1,F_2$ and two singular points $S_1,S_2$ for $a=3,$ $b=\frac32$ and $r=\frac34$}
\label{fig:ellipse2}
\end{figure}

\begin{example}[$r = r_{crit}$]
\label{example5}
Let $a=3,$ $b=\frac32,$ and $r=\frac34.$ Then,
$$
\bV=\bV(y_0^2+4 x_0^2-9, 16 y^2-32 y y_0+16 y_0^2+16 x^2-32 x x_0+16 x_0^2-9, y_0 x+3 y_0 x_0-4 x_0 y)
$$
We find that the reduced \Grobner basis $G_r$ contains eleven polynomials including
\begin{setlength}{\multlinegap}{0pt}
\begin{multline}
g = -29673216 x^4+19035648 y^4+95178240 x^2 y^2-1916928 y^6-14376960 y^4 x^2 \\[0.5ex]
-21012480 y^2 x^4+7372800 x^6+65536 y^8+2162688 y^4 x^4+655360 y^6 x^2\\[0.5ex]
+2621440 x^6 y^2+1048576 x^8-79361856 y^2+119574225-198404640 x^2 
\label{eq:ge2}
\end{multline}
\end{setlength}

\noindent
whose gradient is
\begin{align}
\begin{split}
\grad g &= x (-1854576 x^2+2974320 y^2-449280 y^4-1313280 x^2 y^2+691200 x^4+135168 y^4 x^2\\[0.5ex]
        &\hspace*{2in}   +20480 y^6+245760 y^2 x^4+131072 x^6-6200145)\bi \\[0.5ex]
        &\hspace*{0.18in} + y (594864 y^2+1487160 x^2-89856 y^4-449280 x^2 y^2-328320 x^4+4096 y^6\\[0.5ex]
        &\hspace*{2in}   +67584 y^2 x^4+30720 y^4 x^2+40960 x^6-1240029)\bj\,.
\end{split}
\end{align}

\noindent
There are exactly two solutions to the equation $(\grad g)(p) = 0$ on $\bV(g):$ 
$$
S_1 = \left(0, \frac94\right) \quad \mbox{and} \quad  S_2 = \left(0, -\frac94\right)
$$ 
These points are the {\rm singular points} of the variety $\bV(g):$ At every other point of this variety the gradient $\grad g$ does not vanish. We will collect these points as before into a non-singular subvariety $\Pi.$ The variety $\Pi$ constitutes the set of two parallel lines that we have been seeking. Note that the foci $F_1 = (0,\frac{3\sqrt{3}}{2})$ and $F_2 = (0,-\frac{3\sqrt{3}}{2})$ are still different from the singular points. We show all in Figure~\ref{fig:ellipse2}.
\end{example}

As we can see from the last figure, when $r = r_{crit}= \frac34$ for the given values of $a$ and $b,$ the singular points are located on the 
continuous component of $\Pi.$

Probably the most interesting parallel lines we obtain when $r >  r_{crit}$ as we show in the next example. There, we will find six singular points as each of the two singular points found in Example~\ref{example5} splits into three points.

\begin{figure}[htb]
\centerline{\scalebox{1.30}{\includegraphics{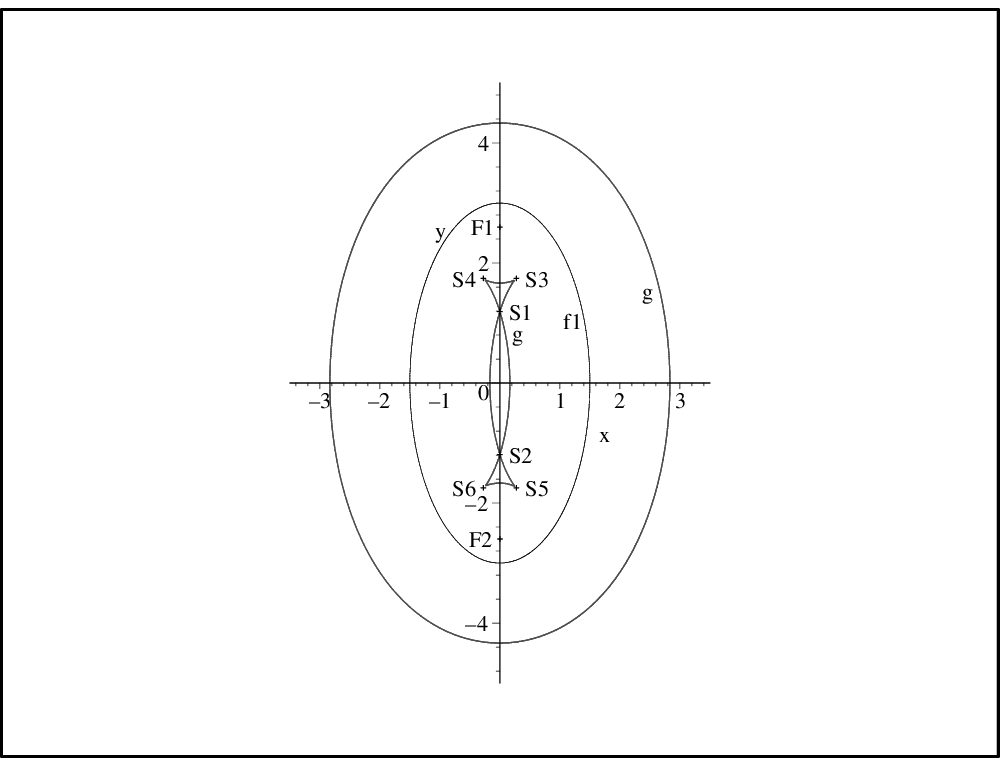}}}
\caption{Ellipse with parallel lines $g,$ foci $F_1,F_2$ and six singular points $S_i,i=1,\ldots,6,$ for $a=3,$ $b=\frac32$ and $r=\frac43$}
\label{fig:ellipse3}
\end{figure}

\begin{example}[$r > r_{crit}$]
\label{example6}
Let $a=3,$ $b=\frac32,$ and $r=\frac32.$ Then,
$$
\bV=\bV(y_0^2+4 x_0^2-9, 9 y^2-18 y y_0+9 y_0^2+9 x^2-18 x x_0+9 x_0^2-16, y_0 x+3 y_0 x_0-4 x_0 y)
$$
We find that again the reduced \Grobner basis $G_r$ contains nine polynomials including
\begin{setlength}{\multlinegap}{0pt}
\begin{multline}
g = 23461560 x^2 y^2+1221025+186624 x^8-2984040 y^4 x^2-5015520 y^2 x^4\\[0.5ex]
    -284472 y^6+384912 y^4 x^4-43658160 x^2+518400 x^6+466560 x^6 y^2\\[0.5ex]
    +116640 y^6 x^2+11664 y^8-2228394 y^2-10769184 x^4+1344177 y^4
\label{eq:ge3}
\end{multline}
\end{setlength}

\noindent
whose gradient is
\begin{align}
\begin{split}
\grad g &= x (977565 y^2+31104 x^6-124335 y^4-417960 x^2 y^2+32076 y^4 x^2-1819090+64800 x^4\\[0.5ex]
        &\hspace*{2in}  +58320 y^2 x^4+4860 y^6-897432 x^2) \bi \\[0.5ex]
        &\hspace*{0.18in} + y (3910260 x^2-994680 x^2 y^2-835920 x^4-142236 y^4+128304 y^2 x^4+77760 x^6\\[0.5ex]
        &\hspace*{2in}  +58320 y^4 x^2+7776 y^6-371399+448059 y^2) \bj\,.
\end{split}
\end{align}

\noindent
There are exactly six solutions to the equation $(\grad g)(p) = 0$ on $\bV(g),$ or, in another words, six singular points of the variety:
\begin{gather*} 
S_1 = \left(0, \frac{\sqrt{51}}{6}\right) \quad \mbox{and} \quad  S_2 = \left(0, -\frac{\sqrt{51}}{6}\right)
\end{gather*}
\begin{align*} 
S_3 &= \left(\frac{ \sqrt{ 525+324 \cdot 6^{\frac23}-864 \cdot 6^{\frac13} }}{18}, 
             \frac{2 \sqrt{231 - 81 \cdot 6^{\frac23} + 54 \cdot 6^{\frac13}}}{9}\right),\\[0.5ex]
S_4 &= \left(-\frac{ \sqrt{ 525+324 \cdot 6^{\frac23}-864 \cdot 6^{\frac13} }}{18}, 
              \frac{2 \sqrt{231 - 81 \cdot 6^{\frac23} + 54 \cdot 6^{\frac13}}}{9}\right),\\[0.5ex]    
S_5 &= \left( \frac{ \sqrt{ 525+324 \cdot 6^{\frac23}-864 \cdot 6^{\frac13} }}{18}, 
             -\frac{2 \sqrt{231 - 81 \cdot 6^{\frac23} + 54 \cdot 6^{\frac13}}}{9}\right),\\[0.5ex]     
S_6 &= \left(-\frac{ \sqrt{ 525+324 \cdot 6^{\frac23}-864 \cdot 6^{\frac13} }}{18}, 
             -\frac{2 \sqrt{231 - 81 \cdot 6^{\frac23} + 54 \cdot 6^{\frac13}}}{9}\right).    
\end{align*}
We display these six points in Figure~\ref{fig:ellipse3}.
\end{example}

\subsection{Parallel lines for arbitrary parameters $a,b,r$}
\label{ellipse2}

In this section we will find a general polynomial $g$ that will generate the second elimination ideal $I_2 = I \cap \BR[x,y,a,b,r]$ where 
$I = \langle f_1, f_2, f_3\rangle$ and the generating polynomials are defined in equations (\ref{eq:fe1}), (\ref{eq:fe2}), and (\ref{eq:fe3}). We will also find a formula for the critical value $r_{crit}$ of $r$ that determines whether $\bV(I_2)$ has two or six singular points, and, we will find general formulas for the coordinates of the said singular points for arbitrary values of the parameters $a,b$ and $r.$  

We begin with a computation of a \Grobner basis for the ideal $I$ for the lexicographic order $y_0 > x_0 > x > y > a > b >r.$
$$
I = \langle b^2 y_0^2+a^2 x_0^2-a^2 b^2, y^2-2 y y_0+y_0^2+x^2-2 x x_0+x_0^2-r^2, b^2 y_0 x-b^2 y_0 x_0-a^2 x_0 y+a^2 x_0 y_0 \rangle
$$
We found that a reduced \Grobner basis for $I$ consists of fifteen homogeneous polynomials of degrees: 
$16$ (two polynomials),
$14$ (one polynomial),
$12$ (one polynomial),
$11$ (two polynomials),
$10$ (two polynomials),
$9$ (one polynomial),
$8$ (one polynomial),
$7$ (two polynomials),
$5$ (one polynomial),
$4$ (one polynomial),
and $2$ (one polynomial). These polynomials are displayed in Appendix~B.

Thus, in this general case we observe that $I_2 = \langle g \rangle$ where $g = g_1/(a^2 b^2)$ is the following homogeneous polynomial in
$\BR[x,y,a,b,r]$ of degree twelve\footnote{We assume $a>b>0.$}:

\beq
\begin{split}
g &= 6\,b^{4}\,x^{4}\,a^{4} - 4\,y^{6}\,b^{4}\,a^{2} - 6\,a^{6}\,x^{4}\,b^{2} - 4\,b^{6}\,x^{2}\,a^{4} + 6\,b^{4}\,x^{2}\,a^{6} 
  + a^{4}\,b^{8} - 2\,y^{2}\,b^{8}\,a^{2} 
\\
& \mytabb + 6\,y^{2}\,b^{6}\,x^{2}\,a^{2} + 6\,a^{6}\,y^{2}\,x^{2}\,b^{2} - 10\,b^{4}\,x^{2}\,a^{4}\,y^{2} + y^{8}\,b^{4} + b^{8}\,
    r^{4} + b^{4}\,r^{8} - 2\,b^{6}\,r^{6} 
\\
& \mytabb  - 2\,b^{4}\,y^{2}\,x^{4}\,r^{2} - 6\,x^{2}\,y^{4}\,b^{4}\,r^{2} + 6\,x^{2}\,b^{4}\,r^{4}\,y^{2} - 4\,b^{4}\,y^{6}\,r^{2}
  + 2\,b^{4}\,y^{6}\,x^{2} + y^{4}\,b^{4}\,x^{4} 
\\
& \mytabb  + b^{4}\,x^{4}\,r^{4} - 2\,b^{4}\,r^{6}\,x^{2} - 2\,b^{6}\,x^{2}\,r^{4} + 6\,y^{4}\,b^{4}\,r^{4} - 2\,b^{8}\,y^{2}\,r^{2}
  + 6\,b^{6}\,r^{4}\,y^{2} - 4\,b^{4}\,r^{6}\,y^{2} 
\\
& \mytabb  + 6\,a^{4}\,y^{2}\,b^{6} + 6\,a^{6}\,y^{2}\,b^{2}\,r^{2} + 2\,a^{4}\,b^{2}\,x^{4}\,y^{2} - 4\,a^{4}\,x^{6}\,b^{2} + 4\,a^{6}\,x^{2}\,b^{2}\,r^{2} - 8\,b^{4}\,x^{2}\,a^{4}\,r^{2} 
\\
& \mytabb  + 10\,a^{4}\,x^{4}\,b^{2}\,r^{2} - 6\,b^{4}\,r^{2}\,a^{2}\,x^{2}\,y^{2} - 6\,a^{4}\,x^{2}\,b^{2}\,r^{2}\,y^{2} - 8\,a^{4}
\,x^{2}\,b^{2}\,r^{4} - 6\,b^{4}\,r^{2}\,a^{2}\,x^{4} 
\\
& \mytabb  + 4\,b^{4}\,r^{4}\,a^{2}\,x^{2} - 6\,y^{4}\,a^{4}\,b^{2}\,x^{2} - 6\,b^{4}\,x^{4}\,a^{2}\,y^{2} + 6\,b^{6}\,x^{2}\,r^{2}
\,a^{2} + 2\,b^{4}\,x^{2}\,y^{4}\,a^{2} 
\\
& \mytabb  + 10\,b^{4}\,y^{4}\,a^{2}\,r^{2} + 4\,b^{6}\,y^{2}\,r^{2}\,a^{2} - 8\,b^{4}\,y^{2}\,r^{4}\,a^{2} - 8\,b^{4}\,y^{2}\,a^{4}
\,r^{2} + 2\,r^{6}\,a^{2}\,b^{4} + a^{8}\,b^{4} 
\\
& \mytabb  - 6\,a^{4}\,y^{4}\,r^{2}\,b^{2} + 4\,a^{4}\,y^{2}\,r^{4}\,b^{2} - 2\,a^{8}\,b^{2}\,x^{2} + x^{4}\,a^{4}\,y^{4} + 4\,x^{2}
\,b^{6}\,y^{2}\,r^{2} - 6\,y^{4}\,b^{6}\,a^{2} 
\\
& \mytabb  - 4\,a^{6}\,y^{2}\,b^{4} - 2\,y^{4}\,a^{4}\,x^{2}\,r^{2} + 2\,b^{2}\,r^{2}\,a^{2}\,x^{4}\,y^{2} - 2\,r^{8}\,a^{2}\,b^{2}
 + r^{8}\,a^{4} + a^{8}\,x^{4} + a^{4}\,x^{8} 
\\
& \mytabb  - 4\,a^{4}\,x^{6}\,r^{2} + 6\,a^{4}\,x^{4}\,r^{4} + 2\,b^{2}\,r^{2}\,a^{2}\,x^{6} - 6\,b^{2}\,r^{4}\,a^{2}\,x^{4} - 2\,a
^{6}\,r^{6} + 2\,a^{6}\,x^{6} + a^{8}\,r^{4} 
\\
& \mytabb  + 2\,y^{6}\,b^{2}\,a^{2}\,x^{2} - 2\,r^{6}\,a^{4}\,y^{2} + 6\,r^{6}\,b^{2}\,y^{2}\,a^{2} - 2\,a^{8}\,x^{2}\,r^{2} + y^{4}
\,a^{4}\,r^{4} - 6\,y^{4}\,b^{2}\,r^{4}\,a^{2} 
\\
& \mytabb  - 2\,a^{6}\,y^{2}\,r^{4} + 4\,x^{2}\,a^{6}\,y^{2}\,r^{2} + 6\,x^{2}\,a^{4}\,y^{2}\,r^{4} + 2\,y^{4}\,b^{2}\,x^{2}\,a^{2}
\,r^{2} - 10\,x^{2}\,b^{2}\,y^{2}\,r^{4}\,a^{2} 
\\
& \mytabb  + 2\,y^{6}\,b^{2}\,a^{2}\,r^{2} + 6\,r^{6}\,b^{2}\,x^{2}\,a^{2} + 4\,b^{2}\,x^{4}\,y^{4}\,a^{2} + 2\,a^{4}\,x^{6}\,y^{2}
 - 6\,a^{6}\,x^{4}\,r^{2} - 4\,r^{6}\,a^{4}\,x^{2} 
\\
& \mytabb  + 2\,b^{2}\,x^{6}\,a^{2}\,y^{2} + 6\,r^{4}\,a^{6}\,x^{2} - 6\,a^{4}\,x^{4}\,r^{2}\,y^{2} + 2\,b^{6}\,y^{6} - 2\,a^{8}\,b
^{2}\,r^{2} - 2\,x^{2}\,y^{4}\,b^{6} 
\\
& \mytabb  + 6\,a^{4}\,y^{4}\,b^{4} + y^{4}\,b^{8} - 6\,y^{4}\,b^{6}\,r^{2} - 2\,x^{4}\,a^{6}\,y^{2} - 2\,a^{6}\,b^{6} + 2\,b^{6}\,r
^{4}\,a^{2} - 2\,b^{8}\,r^{2}\,a^{2} 
\\
& \mytabb  + 2\,b^{4}\,a^{6}\,r^{2} - 6\,b^{4}\,a^{4}\,r^{4} + 2\,r^{2}\,a^{4}\,b^{6} + 2\,r^{4}\,a^{6}\,b^{2} + 2\,r^{6}\,a^{4}\,b^{2}
\end{split}
\label{eq:gellipse}
\eeq
We want to find now those points on the variety $\bV=\bV(g)$ where the gradient $\grad g$ is zero. Thus, $\bV$ can be thought of as union of points where the gradient is zero and those points where the gradient is not zero. Connected components of $\bV$ where the gradient is not zero constitute the $1$-dimensional manifold. We are interested in a discrete finite set $S$ of points where the gradient is zero. These points will give the singular points of the variety $V.$  We will show next that the critical value of $r$ is $r_{crit} = \frac{b^2}{a}.$ Therefore, we compute the gradient 
$\grad g = \grad g(x,y,a,b,r) = <h_1,h_2,h_3,h_4,h_5>$ where
\beq
\begin{split}
h_1 
&=  x(2\,a^{4}\,b^{6} - 2\,a^{4}\,b^{2}\,x^{2}\,y^{2} - 6\,b^{4}\,x^{2}\,a^{4} + 5\,b^{4}\,y^{2}\,a^{4} + b^{6}\,y^{4} + 6\,b^{4}\,y^{2}\,a^{2}\,x^{2} + b^{6}\,r^{4} 
\\
& \hspace*{0.5in} - 2\,b^{6}\,y^{2}\,r^{2} - a^{8}\,x^{2} - 2\,a^{6}\,y^{2}\,r^{2} + a^{4}\,y^{4}\,r^{2} - 3\,a^{4}\,y^{2}\,r^{4} + 2\,a^{4}\,r^{6} + a^{8}\,r^{2} - 3\,a^{6}\,r^{4} 
\\
& \hspace*{0.5in} - b^{4}\,r^{4}\,x^{2} + b^{4}\,r^{6} - 4\,x^{2}\,b^{2}\,a^{2}\,y^{4} - 3\,x^{4}\,b^{2}\,a^{2}\,y^{2} - x^{2}\,b^{4}\,y^{4} - a^{2}\,y^{6}\,b^{2} - 3\,x^{4}\,a^{4}\,y^{2} 
\\
& \hspace*{0.5in} + 6\,a^{4}\,x^{2}\,r^{2}\,y^{2} - b^{2}\,y^{4}\,a^{2}\,r^{2} + 5\,b^{2}\,y^{2}\,r^{4}\,a^{2} - 3\,b^{2}\,a^{2}\,x^{4}\,r^{2} + 6\,b^{2}\,a^{2}\,x^{2}\,r^{4} + 6\,a^{4}\,x^{4}\,r^{2} 
\\
& \hspace*{0.5in} + 6\,a^{6}\,x^{2}\,r^{2} - 6\,a^{4}\,x^{2}\,r^{4} + 2\,b^{4}\,y^{2}\,x^{2}\,r^{2} - 3\,b^{4}\,r^{4}\,y^{2} + 3\,y^{4}\,b^{4}\,r^{2}   - 3\,b^{2}\,r^{6}\,a^{2} 
\\
& \hspace*{0.5in} - 2\,b^{2}\,a^{2}\,x^{2}\,r^{2}\,y^{2} + 3\,y^{4}\,a^{4}\,b^{2} + a^{8}\,b^{2} - 2\,a^{6}\,b^{2}\,r^{2} + 4\,a^{4}\,b^{2}\,r^{4} + 4\,a^{4}\,b^{4}\,r^{2} - 3\,a^{6}\,x^{4} 
\\
& \hspace*{0.5in} - b^{4}\,y^{6} + 3\,b^{4}\,y^{2}\,a^{2}\,r^{2} - 10\,a^{4}\,b^{2}\,x^{2}\,r^{2} + 3\,a^{4}\,b^{2}\,y^{2}\,r^{2} + 6\,b^{4}\,x^{2}\,a^{2}\,r^{2} + 6\,x^{2}\,a^{6}\,b^{2} 
\\
& \hspace*{0.5in} - 3\,a^{6}\,b^{2}\,y^{2} - 3\,b^{6}\,y^{2}\,a^{2} - y^{4}\,b^{4}\,a^{2} - 3\,b^{6}\,r^{2}\,a^{2} - 2\,b^{4}\,r^{4}\,a^{2} + 2\,x^{2}\,y^{2}\,a^{6} + 6\,x^{4}\,a^{4}\,b^{2} 
\\
& \hspace*{0.5in} - 3\,a^{6}\,b^{4} - a^{4}\,y^{4}\,x^{2} - 2\,x^{6}\,a^{4})
\end{split}
\eeq
\beq
\begin{split}
h_2 
&=   y(3\,a^{2}\,b^{4}\,x^{4} - 3\,a^{4}\,b^{6} + 6\,a^{4}\,b^{2}\,x^{2}\,y^{2} + 5\,b^{4}\,x^{2}\,a^{4} - 6\,b^{4}\,y^{2}\,a^{4} + 2\,b^{6}\,y^{2}\,x^{2} - 3\,b^{6}\,y^{4} 
\\
& \hspace*{0.5in}  - 2\,b^{4}\,y^{2}\,a^{2}\,x^{2} - 3\,b^{6}\,r^{4} + 6\,b^{6}\,y^{2}\,r^{2} - a^{4}\,y^{2}\,r^{4} + a^{4}\,r^{6} + a^{6}\,r^{4} + x^{4}\,b^{4}\,r^{2} - 3\,b^{4}\,r^{4}\,x^{2} 
\\
& \hspace*{0.5in}  - 2\,b^{6}\,x^{2}\,r^{2} + b^{8}\,r^{2} + 2\,b^{4}\,r^{6} - 3\,x^{2}\,b^{2}\,a^{2}\,y^{4} - 4\,x^{4}\,b^{2}\,a^{2}\,y^{2}
 - x^{6}\,b^{2}\,a^{2} - 3\,x^{2}\,b^{4}\,y^{4} 
\\
& \hspace*{0.5in}  - x^{4}\,b^{4}\,y^{2} - x^{4}\,a^{4}\,y^{2} + 2\,a^{4}\,x^{2}\,r^{2}\,y^{2} - 3\,b^{2}\,y^{4}\,a^{2}\,r^{2} + 6\,b^{2}\,y^{2}\,r^{4}\,a^{2} - b^{2}\,a^{2}\,x^{4}\,r^{2} 
\\
& \hspace*{0.5in}  + 5\,b^{2}\,a^{2}\,x^{2}\,r^{4} + 3\,a^{4}\,x^{4}\,r^{2} - 2\,a^{6}\,x^{2}\,r^{2} - 3\,a^{4}\,x^{2}\,r^{4} + 6\,b^{4}\,y^{2}\,x^{2}\,r^{2} - 6\,b^{4}\,r^{4}\,y^{2} 
\\
& \hspace*{0.5in}  + 6\,y^{4}\,b^{4}\,r^{2} - 3\,b^{2}\,r^{6}\,a^{2} - 2\,b^{2}\,a^{2}\,x^{2}\,r^{2}\,y^{2} + b^{8}\,a^{2} - 3\,a^{6}\,b^{2}
\,r^{2} - 2\,a^{4}\,b^{2}\,r^{4} 
\\
& \hspace*{0.5in}  + 4\,a^{4}\,b^{4}\,r^{2} + a^{6}\,x^{4} - b^{8}\,y^{2} - 2\,b^{4}\,y^{6} - 10\,b^{4}\,y^{2}\,a^{2}\,r^{2} + 3\,a^{4}\,b^{2
}\,x^{2}\,r^{2} + 6\,a^{4}\,b^{2}\,y^{2}\,r^{2} 
\\
& \hspace*{0.5in}  + 3\,b^{4}\,x^{2}\,a^{2}\,r^{2} - 3\,x^{2}\,a^{6}\,b^{2} - 3\,x^{2}\,b^{6}\,a^{2} + 6\,b^{6}\,y^{2}\,a^{2} + 6\,y^{4}\,b
^{4}\,a^{2} - 2\,b^{6}\,r^{2}\,a^{2} 
\\
& \hspace*{0.5in}  + 4\,b^{4}\,r^{4}\,a^{2} - x^{4}\,a^{4}\,b^{2} + 2\,a^{6}\,b^{4} - x^{6}\,a^{4})
\end{split}
\eeq
\beq
\begin{split}
h_3 
&= 6\,a^{2}\,b^{4}\,x^{4} + 3\,a^{4}\,b^{6} - 9\,a^{4}\,b^{2}\,x^{2}\,y^{2} - 9\,b^{4}\,x^{2}\,a^{4} + 6\,b^{4}\,y^{2}\,a^{4} - 3\,b^{6}\,y^{2}\,x^{2} + 3\,b^{6}\,y^{4} 
\\
& \hspace*{0.5in}  + 10\,b^{4}\,y^{2}\,a^{2}\,x^{2} - b^{6}\,r^{4} - 2\,b^{6}\,y^{2}\,r^{2} + 3\,a^{4}\,y^{2}\,r^{4} + 3\,a^{4}\,r^{6} - 2\,a^{6}\,r^{4} + 3\,x^{4}\,b^{4}\,r^{2} 
\\
& \hspace*{0.5in}  - 2\,b^{4}\,r^{4}\,x^{2} - 3\,b^{6}\,x^{2}\,r^{2} + b^{8}\,r^{2} - b^{4}\,r^{6} + 6\,x^{2}\,b^{2}\,a^{2}\,y^{4} - 2\,x^{4}
\,b^{2}\,a^{2}\,y^{2} + 4\,x^{6}\,b^{2}\,a^{2} 
\\
& \hspace*{0.5in}   - x^{2}\,b^{4}\,y^{4} + 3\,x^{4}\,b^{4}\,y^{2} - a^{2}\,r^{8} - 3\,r^{6}\,x^{2}\,b^{2} - 3\,r^{6}\,b^{2}\,y^{2} + 5\,x^{2}
\,b^{2}\,y^{2}\,r^{4} - x^{4}\,b^{2}\,r^{2}\,y^{2} 
\\
& \hspace*{0.5in}   - x^{2}\,b^{2}\,y^{4}\,r^{2} - a^{2}\,x^{8} - x^{6}\,b^{2}\,y^{2} - 2\,x^{4}\,b^{2}\,y^{4} - 6\,a^{2}\,x^{2}\,y^{2}\,r^{4}
 - y^{6}\,b^{2}\,r^{2} + 3\,y^{4}\,b^{2}\,r^{4} 
\\
& \hspace*{0.5in}   + 6\,a^{2}\,x^{4}\,r^{2}\,y^{2} + 2\,a^{2}\,y^{4}\,x^{2}\,r^{2} - b^{2}\,x^{2}\,y^{6} - x^{6}\,b^{2}\,r^{2} + 3\,x^{4}\,b
^{2}\,r^{4} + 3\,x^{4}\,a^{4}\,y^{2} 
\\
& \hspace*{0.5in}   - 6\,a^{4}\,x^{2}\,r^{2}\,y^{2} + 6\,b^{2}\,y^{4}\,a^{2}\,r^{2} - 4\,b^{2}\,y^{2}\,r^{4}\,a^{2} - 10\,b^{2}\,a^{2}\,x^{4}
\,r^{2} + 8\,b^{2}\,a^{2}\,x^{2}\,r^{4} 
\\
& \hspace*{0.5in}   + 9\,a^{4}\,x^{4}\,r^{2} + 4\,a^{6}\,x^{2}\,r^{2} - 9\,a^{4}\,x^{2}\,r^{4} + 3\,b^{4}\,y^{2}\,x^{2}\,r^{2} + 4\,b^{4}\,r
^{4}\,y^{2} - 5\,y^{4}\,b^{4}\,r^{2} 
\\
& \hspace*{0.5in}   - 2\,b^{2}\,r^{6}\,a^{2} + 6\,b^{2}\,a^{2}\,x^{2}\,r^{2}\,y^{2} - b^{8}\,a^{2} + 4\,a^{6}\,b^{2}\,r^{2} - 3\,a^{4}\,b^{2}
\,r^{4} - 3\,a^{4}\,b^{4}\,r^{2} - 2\,a^{6}\,x^{4} 
\\
& \hspace*{0.5in}   + b^{8}\,y^{2} + 2\,b^{4}\,y^{6} + 8\,b^{4}\,y^{2}\,a^{2}\,r^{2} - 6\,a^{4}\,b^{2}\,x^{2}\,r^{2} - 9\,a^{4}\,b^{2}\,y^{2}
\,r^{2} + 8\,b^{4}\,x^{2}\,a^{2}\,r^{2} 
\\
& \hspace*{0.5in}  + 4\,x^{2}\,a^{6}\,b^{2} + 4\,x^{2}\,b^{6}\,a^{2} - 6\,b^{6}\,y^{2}\,a^{2} - 6\,y^{4}\,b^{4}\,a^{2} - 2\,b^{6}\,r^{2}\,a
^{2} + 6\,b^{4}\,r^{4}\,a^{2} 
\\
& \hspace*{0.5in}   + 9\,x^{4}\,a^{4}\,b^{2} - 2\,a^{6}\,b^{4} - 3\,x^{6}\,a^{4} + 4\,a^{2}\,x^{6}\,r^{2} + 4\,a^{2}\,r^{6}\,x^{2} - a^{2}\,y
^{4}\,x^{4} - 6\,a^{2}\,x^{4}\,r^{4} 
\\
& \hspace*{0.5in}   - 2\,a^{2}\,x^{6}\,y^{2} - a^{2}\,y^{4}\,r^{4} + 2\,a^{2}\,r^{6}\,y^{2} + r^{8}\,b^{2}
\end{split}
\eeq
\beq
\begin{split}
h_4 
&= 2\,a^{4}\,b^{6} - 10\,a^{4}\,b^{2}\,x^{2}\,y^{2} - 6\,b^{4}\,x^{2}\,a^{4} + 9\,b^{4}\,y^{2}\,a^{4} + 2\,b^{6}\,y^{4} + 9\,b^{4}\,y^{2}\,a^{2}\,x^{2} + 2\,b^{6}\,r^{4} 
\\
& \hspace*{0.5in}   - 4\,b^{6}\,y^{2}\,r^{2} - a^{8}\,x^{2} + 3\,a^{6}\,y^{2}\,r^{2} - 3\,a^{4}\,y^{4}\,r^{2} + 2\,a^{4}\,y^{2}\,r^{4} + a^{4}
\,r^{6} - a^{8}\,r^{2} + a^{6}\,r^{4} 
\\
& \hspace*{0.5in}   - 3\,b^{4}\,r^{4}\,x^{2} - 3\,b^{4}\,r^{6} + 2\,x^{2}\,b^{2}\,a^{2}\,y^{4} - 6\,x^{4}\,b^{2}\,a^{2}\,y^{2} - 3\,x^{2}\,b
^{4}\,y^{4} - 4\,a^{2}\,y^{6}\,b^{2} - a^{2}\,r^{8} 
\\
& \hspace*{0.5in}    - 2\,r^{6}\,x^{2}\,b^{2} - 4\,r^{6}\,b^{2}\,y^{2} + 6\,x^{2}\,b^{2}\,y^{2}\,r^{4} - 2\,x^{4}\,b^{2}\,r^{2}\,y^{2} - 6\,x^{2}\,b^{2}\,y^{4}\,r^{2} + x^{4}\,b^{2}\,y^{4} 
\\
& \hspace*{0.5in}    - 5\,a^{2}\,x^{2}\,y^{2}\,r^{4} - 4\,y^{6}\,b^{2}\,r^{2} + 6\,y^{4}\,b^{2}\,r^{4} + a^{2}\,x^{4}\,r^{2}\,y^{2} + a^{2}\,y^{4}\,x^{2}\,r^{2} + 2\,b^{2}\,x^{2}\,y^{6} 
\\
& \hspace*{0.5in}    + x^{4}\,b^{2}\,r^{4} + x^{4}\,a^{4}\,y^{2} - 3\,a^{4}\,x^{2}\,r^{2}\,y^{2} + 10\,b^{2}\,y^{4}\,a^{2}\,r^{2} - 8\,b^{2}\,y
^{2}\,r^{4}\,a^{2} - 6\,b^{2}\,a^{2}\,x^{4}\,r^{2} 
\\
& \hspace*{0.5in}    + 4\,b^{2}\,a^{2}\,x^{2}\,r^{4} + 5\,a^{4}\,x^{4}\,r^{2} + 2\,a^{6}\,x^{2}\,r^{2} - 4\,a^{4}\,x^{2}\,r^{4} + 6\,b^{4}\,y
^{2}\,x^{2}\,r^{2} + 9\,b^{4}\,r^{4}\,y^{2} 
\\
& \hspace*{0.5in}    - 9\,y^{4}\,b^{4}\,r^{2} + 2\,b^{2}\,r^{6}\,a^{2} - 6\,b^{2}\,a^{2}\,x^{2}\,r^{2}\,y^{2} + 6\,y^{4}\,a^{4}\,b^{2} + a^{8}
\,b^{2} + 2\,a^{6}\,b^{2}\,r^{2} 
\\
& \hspace*{0.5in}   - 6\,a^{4}\,b^{2}\,r^{4} + 3\,a^{4}\,b^{4}\,r^{2} - 3\,a^{6}\,x^{4} + 3\,b^{4}\,y^{6} + 6\,b^{4}\,y^{2}\,a^{2}\,r^{2} - 8
\,a^{4}\,b^{2}\,x^{2}\,r^{2} 
\\
& \hspace*{0.5in}   - 8\,a^{4}\,b^{2}\,y^{2}\,r^{2} + 9\,b^{4}\,x^{2}\,a^{2}\,r^{2} + b^{2}\,y^{8} + a^{2}\,y^{6}\,r^{2} + x^{2}\,a^{2}\,y^{6
} + 6\,x^{2}\,a^{6}\,b^{2} - 4\,a^{6}\,b^{2}\,y^{2} 
\\
& \hspace*{0.5in}   - 4\,b^{6}\,y^{2}\,a^{2} - 9\,y^{4}\,b^{4}\,a^{2} - 4\,b^{6}\,r^{2}\,a^{2} + 3\,b^{4}\,r^{4}\,a^{2} + 3\,x^{2}\,y^{2}\,a
^{6} + 6\,x^{4}\,a^{4}\,b^{2} - 3\,a^{6}\,b^{4} 
\\
& \hspace*{0.5in}   - 3\,a^{4}\,y^{4}\,x^{2} - 2\,x^{6}\,a^{4} + a^{2}\,x^{6}\,r^{2} + 3\,a^{2}\,r^{6}\,x^{2} + 2\,a^{2}\,y^{4}\,x^{4} - 3\,a
^{2}\,x^{4}\,r^{4} + a^{2}\,x^{6}\,y^{2} 
\\
& \hspace*{0.5in}   - 3\,a^{2}\,y^{4}\,r^{4} + 3\,a^{2}\,r^{6}\,y^{2} + r^{8}\,b^{2}
\end{split}
\eeq
\beq
\begin{split}
h_5 
&=3\,a^{2}\,b^{4}\,x^{4} + a^{4}\,b^{6} - 3\,a^{4}\,b^{2}\,x^{2}\,y^{2} - 4\,b^{4}\,x^{2}\,a^{4} - 4\,b^{4}\,y^{2}\,a^{4} + 2\,b^{6}\,y^{2}\,x^{2} - 3\,b^{6}\,y^{4} 
\\
& \hspace*{0.5in}   - 3\,b^{4}\,y^{2}\,a^{2}\,x^{2} - 3\,b^{6}\,r^{4} + 6\,b^{6}\,y^{2}\,r^{2} - a^{8}\,x^{2} - 2\,a^{6}\,y^{2}\,r^{2} + a^{4
}\,y^{4}\,r^{2} - 3\,a^{4}\,y^{2}\,r^{4} 
\\
& \hspace*{0.5in}   + 2\,a^{4}\,r^{6} + a^{8}\,r^{2} - 3\,a^{6}\,r^{4} + x^{4}\,b^{4}\,r^{2} - 3\,b^{4}\,r^{4}\,x^{2} - 2\,b^{6}\,x^{2}\,r^{2}
 + b^{8}\,r^{2} + 2\,b^{4}\,r^{6} 
\\
& \hspace*{0.5in}   + x^{2}\,b^{2}\,a^{2}\,y^{4} + x^{4}\,b^{2}\,a^{2}\,y^{2} + x^{6}\,b^{2}\,a^{2} - 3\,x^{2}\,b^{4}\,y^{4} - x^{4}\,b^{4}\,y
^{2} + a^{2}\,y^{6}\,b^{2} - 3\,x^{4}\,a^{4}\,y^{2} 
\\
& \hspace*{0.5in}   + 6\,a^{4}\,x^{2}\,r^{2}\,y^{2} - 6\,b^{2}\,y^{4}\,a^{2}\,r^{2} + 9\,b^{2}\,y^{2}\,r^{4}\,a^{2} - 6\,b^{2}\,a^{2}\,x^{4}
\,r^{2} + 9\,b^{2}\,a^{2}\,x^{2}\,r^{4} 
\\
& \hspace*{0.5in}   + 6\,a^{4}\,x^{4}\,r^{2} + 6\,a^{6}\,x^{2}\,r^{2} - 6\,a^{4}\,x^{2}\,r^{4} + 6\,b^{4}\,y^{2}\,x^{2}\,r^{2} - 6\,b^{4}\,r
^{4}\,y^{2} + 6\,y^{4}\,b^{4}\,r^{2} 
\\
& \hspace*{0.5in}   - 4\,b^{2}\,r^{6}\,a^{2} - 10\,b^{2}\,a^{2}\,x^{2}\,r^{2}\,y^{2} - 3\,y^{4}\,a^{4}\,b^{2} - a^{8}\,b^{2} - b^{8}\,a^{2} + 
2\,a^{6}\,b^{2}\,r^{2} + 3\,a^{4}\,b^{2}\,r^{4} 
\\
& \hspace*{0.5in}   - 6\,a^{4}\,b^{4}\,r^{2} - 3\,a^{6}\,x^{4} - b^{8}\,y^{2} - 2\,b^{4}\,y^{6} - 8\,b^{4}\,y^{2}\,a^{2}\,r^{2} - 8\,a^{4}\,b
^{2}\,x^{2}\,r^{2} + 4\,a^{4}\,b^{2}\,y^{2}\,r^{2} 
\\
& \hspace*{0.5in}   + 4\,b^{4}\,x^{2}\,a^{2}\,r^{2} + 2\,x^{2}\,a^{6}\,b^{2} + 3\,x^{2}\,b^{6}\,a^{2} + 3\,a^{6}\,b^{2}\,y^{2} + 2\,b^{6}\,y
^{2}\,a^{2} + 5\,y^{4}\,b^{4}\,a^{2} 
\\
& \hspace*{0.5in}   + 2\,b^{6}\,r^{2}\,a^{2} + 3\,b^{4}\,r^{4}\,a^{2} + 2\,x^{2}\,y^{2}\,a^{6} + 5\,x^{4}\,a^{4}\,b^{2} + a^{6}\,b^{4} - a^{4
}\,y^{4}\,x^{2} - 2\,x^{6}\,a^{4}
\end{split}
\eeq
are polynomials in $\BR[a,y,a,b,r].$ Continuing our approach from Section~\ref{parabola}, in order to find solutions to the vector equation $\grad g=0$ in $\bV(g),$ we need to solve a system of six polynomial equations $h_i=0,i=1,\ldots,5,$ and $g=0$ in the variables $x,y,a,b,r$ under the assumption that $a>b>r>0.$ Let 
$I_S = \langle \grad g,g\rangle = \langle h_1, h_2, h_3, h_4, h_5, g \rangle$ and let $F$ be a list containing the six polynomial equations:
\beq
F = [h_1 = 0, h_2 = 0, h_3 = 0, h_4 = 0, h_5 = 0, g = 0]
\label{eq:Fellipse} 
\eeq
In order to solve system (\ref{eq:Fellipse}), we can either employ the \Grobner basis approach again (see Appendix~C) or simply use {\tt solve} command in {\tt Maple}. Thus we find that the affine variety $\bV(I_S)$ consists of the following eight singular points of which six are real and two are complex:
\beq
S_1 = \left(0, \,{\displaystyle \frac {\sqrt{ - (a^{2} - b^{2})\,(r^{2} - b^{2})}}{b}} \right), \qquad
S_2 = \left(0, \, - {\displaystyle \frac {\sqrt{ - (a^{2} - b^{2})\,(r^{2} - b^{2})}}{b}}  \right),
\label{eq:S12}
\eeq
\beq
S_3 = \left({\displaystyle \frac {\sqrt{(a^{2} - b^{2})\,( - a^{2} + r^{2})}}{a}} , \,0 \right), \qquad 
S_4 = \left( - {\displaystyle \frac {\sqrt{(a^{2} - b^{2})\,( - a^{2} + r^{2})}}{a}} , \,0 \right),
\eeq
\beq
S_5 = \left(\sqrt{ - {\displaystyle \frac {\delta }{( - a^{2} + b^{2})\,(a\,b\,r)^{\frac23}}} }, \, - {\displaystyle \frac {\sqrt{( - a^{2} 
              + b^{2})\,\varepsilon }}{(a\,b\,r)^{\frac13}\,( - a^{2} + b^{2})}}  \right),
\label{eq:S5}
\eeq
\beq
S_6 = \left( - \sqrt{ - {\displaystyle \frac {\delta }{( - a^{2} + b^{2})\,(a\,b\,r)^{\frac23}}} }, \, - {\displaystyle \frac {\sqrt{( - a^{2}
                 + b^{2})\,\varepsilon }}{(a\,b\,r)^{\frac13}\,( - a^{2} + b^{2})}}  \right),
\label{eq:S6}
\eeq
\beq
S_7 = \left(\sqrt{ - {\displaystyle \frac {\delta }{( - a^{2} + b^{2})\,(a\,b\,r)^{\frac23}}} }, \,{\displaystyle \frac {\sqrt{( - a^{2} 
              + b^{2})\,\varepsilon }}{(a\,b\,r)^{\frac13}\,( - a^{2} + b^{2})}}  \right),
\label{eq:S7}
\eeq
\beq
S_8 = \left( - \sqrt{ - {\displaystyle \frac {\delta }{( - a^{2} + b^{2})\,(a\,b\,r)^{\frac23}}} }, \,{\displaystyle \frac {\sqrt{( - a^{2} 
                 + b^{2})\,\varepsilon }}{(a\,b\,r)^{\frac13}\,( - a^{2} + b^{2})}}  \right),
\label{eq:S8}
\eeq
where
\begin{align}
\delta   &=  - 3\,b^{2}\,a^{2}\,r^{2} + r^{2}\,a^{2}\,(a\,b\,r)^{\frac23} + 3\,(a\,b\,r)^{\frac13}\,b^{3}\,a\,r - b^{4}\,(a\,b\,r)^{\frac23},
\label{eq:delta}\\[0.5ex]
\varepsilon &=  - a^{4}\,(a\,b\,r)^{\frac23} + 3\,(a\,b\,r)^{\frac13}\,b\,r\,a^{3} - 3\,b^{2}\,a^{2}\,r^{2} + (a\,b\,r)^{\frac23}\,b^{2}\,r^{2}.
\label{eq:epsilon}
\end{align}
Notice that $\delta$ becomes zero in three cases:
\begin{enumerate}
\item When $r = 0$ -- this is not possible since we expect $r > 0,$

\item When $r = - \frac{b^2}{a}$  -- this is not possible as $r > 0,$

\item When $r = \frac{b^2}{a}.$
\end{enumerate}
Only this last value is possible thus we set $r_{crit} = \frac{b^2}{a}.$ When $r = r_{crit},$ then under the assumptions that $a>b>r>0,$ we essentially obtain two real singular points:
\beq
S_1 = S_5 = S_6 = \left( 0,\frac{a^2-b^2}{a}\right) \quad \mbox{and} \quad S_2 = S_7 = S_8 = \left(0,-\frac{a^2-b^2}{a}\right),
\label{eq:S156}
\eeq
while the complex points $S_3$ and $S_4$ become:
\beq
S_3 = \left(\frac{\sqrt{a^2+b^2}(a^2-b^2)}{a^2}i,0 \right) \quad \mbox{and} \quad S_4 = \left(- \frac{\sqrt{a^2+b^2}(a^2-b^2)}{a^2}i,0 \right).
\eeq
For the ellipse with $a > b,$ the ratio $l = \frac{b^2}{a}$ is again called a {\it semi-latus rectum\/.} It is the distance between a focus and the ellipse itself measured along a line perpendicular to the major axis. The value of  $l$ is also the reciprocal of the maximum curvature $\kappa_{max}$ of the ellipse, or, it is equal to the radius of the smallest osculating circle.~\cite{wiki}

To summarize, the variety $\bV(I_S)$ contains the following real singular points:
\begin{enumerate}
\item Two virtual points $S_1$ and $S_2$ shown in (\ref{eq:S12}) when $0 < r < r_{crit}.$ In particular, when $a=3,b=\frac32,$ and $r=\frac12,$ upon substituting these values into (\ref{eq:gellipse}) we obtain the parallel lines (\ref{eq:ge4}) derived in Example~\ref{example4}. Furthermore, when these values are substituted into (\ref{eq:S12}), then we obtain the two virtual singular points found in the same example.
\item Two points $S_1$ and $S_2$ shown in (\ref{eq:S156}) when $r = r_{crit}.$ In particular, we can recover results of Example~\ref{example5} when we substitute values $a=3,$ $b=\frac32,$ and $r=\frac34$ into the general equation (\ref{eq:gellipse}) and the singular points (\ref{eq:S156}). 
\item Six points $S_1,S_2,S_5,S_6,S_7,S_8$ shown in (\ref{eq:S12}), (\ref{eq:S5}), (\ref{eq:S6}), (\ref{eq:S7}), and (\ref{eq:S8}) when $r>r_{crit}.$ In particular we can recover in the same manner six points and the equation for the parallel lines found in Example~\ref{example6} when we set  
$a=3,$ $b=\frac32,$ and $r=\frac32.$ These six singular points are displayed in Figure~\ref{fig:ellipse3}.
\end{enumerate}

\section{Parallel Lines to a Hyperbola}
\label{hyperbola}

Since a parabola, an ellipse, and a hyperbola belong to the same class of nondegenerate conic 
$$
C:\, (X^2 + Y^2 - Z^2 = 0)
$$ 
in $\BP^2(\BR),$ (cf. \cite{Reid}) we expect essentially similar results for a hyperbola. Thus, our approach will be similar to the one for the ellipse from Section~\ref{ellipse}. In particular, there is only one sign change needed in the polynomials $f_1$ and $f_2$ shown in (\ref{eq:fe1}) and 
(\ref{eq:fe3}) whereas the polynomial $f_2$ remains the same. Thus, we have now
\beq
f_1 = b^2 y_0^2 - a^2 x_0^2-a^2 b^2 = 0
\label{eq:fh1}
\eeq
where we assume that $a>b$ and the center of the hyperbola is at the origin. Polynomial $f_2$ defines a circle of radius (offset) $r$ centered at a point $P=(x_0,y_0)$ on the hyperbola,
\beq
f_2 = (y-y_0)^2+(x-x_0)^2-r^2 = 0  
\label{eq:fh2}
\eeq
while polynomial $f_3$
\beq
f_3 = b^2 (x-x_0) y_0 + a^2 (y-y_0) x_0 = 0 
\label{eq:fh3}
\eeq
gives a condition that a point $P(x,y)$ on the circle $f_2$ lies on a line perpendicular to the hyperbola $f_1$ and passing through the point 
$(x_0,y_0)$ on the  hyperbola. All these points $P$ belong to the affine variety $\bV = \bV(f_1,f_2,f_3)$ -- the  hyperbola envelope -- and, like for the parabola in Section~\ref{parabola} and the ellipse in Section~\ref{ellipse}, define two parallel lines at the distance $r$ from the hyperbola. We proceed like in Section~\ref{ellipse}: We will compute a reduced \Grobner basis for $\bV$ and by eliminating variables $x_0,y_0$ we will obtain one polynomial $g$ in variables $x$ and $y$ alone (also $a,b,r$ in the general case) that defines both parallel lines on each side of the  hyperbola for the given value of the offset $r.$ Thus, in general, we have $f_1,f_2,f_3 \in \BR[y_0,x_0,x,y,a,b,r].$ However, first we illustrate computations for specific values of the parameters $a,b$ and $r,$ and later we show computation for unspecified parameters.

\subsection{Parallel lines for specific parameters $a,b,r$}
\label{hyperbola1}

The parallel lines constitute a subvariety of the affine variety $\bV = \bV(I)$ where, like before, $I \langle f_1, f_2, f_3\rangle.$ In order to find the polynomial $g$ for the second elimination ideal $I_2 = I \cap \BR[x,y],$ we first find a \Grobner basis~$G$ for $\bV$ for the lexicographic order 
$y_0> x_0 > x >y$ and reduce it eventually to a unique reduced basis~$G_r$ of nine (or eleven, depending on the case) polynomials. Finally, we eliminate from that basis all polynomials that contain variables $y_0$ and $x_0,$ and keep the single polynomial $g \in \BR[x,y]$ that gives the parallel lines and generates the second elimination ideal $I_2= \langle g \rangle.$ We will also proceed to find all singular points of $\bV(g)$ where $(\grad g)(p) = 0.$ Thus, let's consider three cases when $0<r<r_{crit},\,r = r_{crit},$ and $r>r_{crit}.$ In the following section we will consider a general case for arbitrary values of $a,b\; (a>b>0)$ and $r\; (r>0).$

\begin{figure}[htb]
\centerline{\scalebox{1.30}{\includegraphics{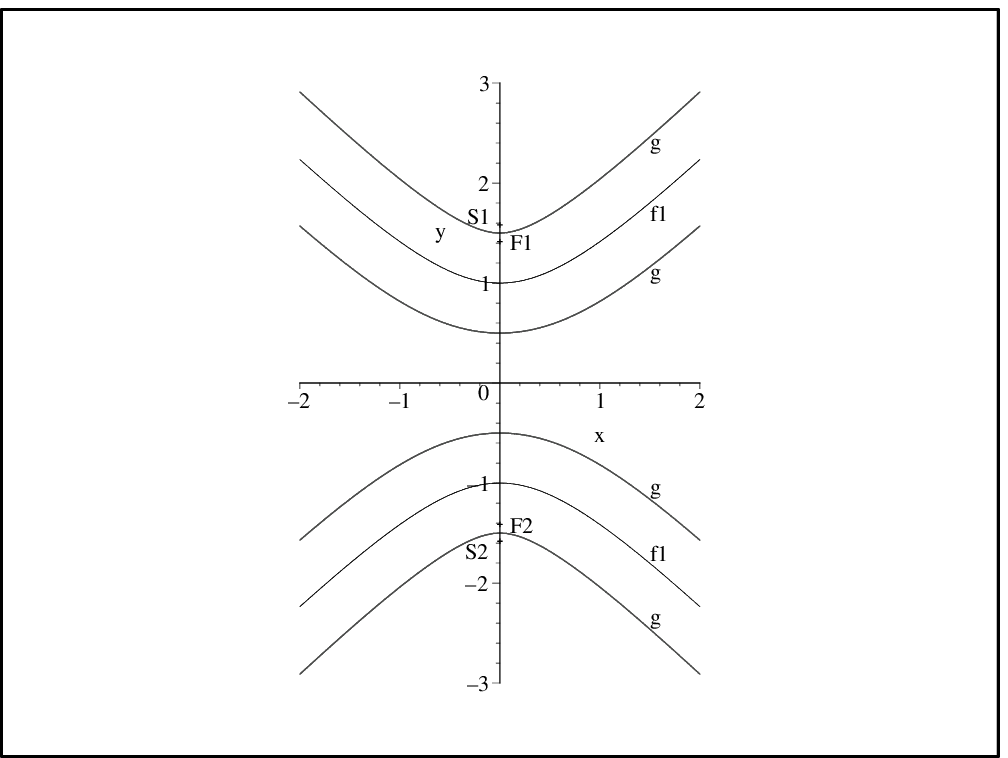}}}
\caption{Hyperbola with parallel lines $g,$ foci $F_1,F_2$ and two singular points $S_1,S_2$ for $a=1,b=1$ and $r=\frac12$}
\label{fig:hyperbola1}
\end{figure}

\begin{example}[$r < r_{crit}$]
\label{example7}
Let $a=3,$ $b=\frac32,$ and $r=\frac12.$ Then,
$$
\bV=\bV(y_0^2-x_0^2-1, 4 y^2-8 y y_0+4 y_0^2+4 x^2-8 x x_0+4 x_0^2-1, x y_0+y x_0-2 y_0 x_0 )
$$
We find that the reduced \Grobner basis $G_r$ contains nine polynomials including
\begin{multline}
g = 225+64 y^8+288 x^6+480 x^2 y^4-800 y^2 x^4-480 y^6-128 y^4 x^4\\[0.5ex]
    -1180 y^2+516 x^2-1320 x^2 y^2+532 x^4+1236 y^4+64 x^8 
\label{eq:gh7}
\end{multline}
whose gradient is
\begin{align}
\begin{split}
\grad g &= x (216 x^4+120 y^4-400 x^2 y^2-64 x^2 y^4 + 129-330 y^2+266 x^2+64 x^6)  \bi \\[0.5ex]
        &\hspace*{0.18in} -y (-64 y^6-240 x^2 y^2+200 x^4+360 y^4+64 y^2 x^4+295+330 x^2-618 y^2) \bj\,.
\end{split}
\end{align}
\noindent
There are exactly two solutions to the equation $(\grad g)(p) = 0$ on $\bV(g):$ 
$$
S_1 = \left(0, \frac{\sqrt{10}}{2}\right) \quad \mbox{and} \quad  S_2 = \left(0, -\frac{\sqrt{10}}{2}\right).
$$ 
We will refer to these points as {\rm virtual singular points} for the following reason:  $\bV(g)$ is a disjoint union of $\{S_1,S_2\}$ and a non-singular smooth subvariety $\Pi$ where the gradient $\grad g$ does not vanish. The subvariety $\Pi$ constitutes the set of two parallel lines that we have been seeking. Note that the foci $F_1 = (0,\sqrt{2})$ and $F_2 = (0,-\sqrt{2})$ are different from the singular points. We show all in 
Figure~\ref{fig:hyperbola1}.
\end{example}

We have verified by direct computation that as the value of the offset $r$ increases towards some critical value $r_{crit},$ the virtual singular points approach the smooth component of the variety $\bV(g),$ that is, they approach the parallel lines $\Pi.$\footnote{In the next subsection we will derive a formula for the critical value $r_{crit}$ of $r$ for a general hyperbola.} As the next example shows, when $r = r_{crit},$ then the virtual singular points are actually located on the smooth component $\Pi$ of the variety $\bV(g).$

\begin{figure}[htb]
\centerline{\scalebox{1.30}{\includegraphics{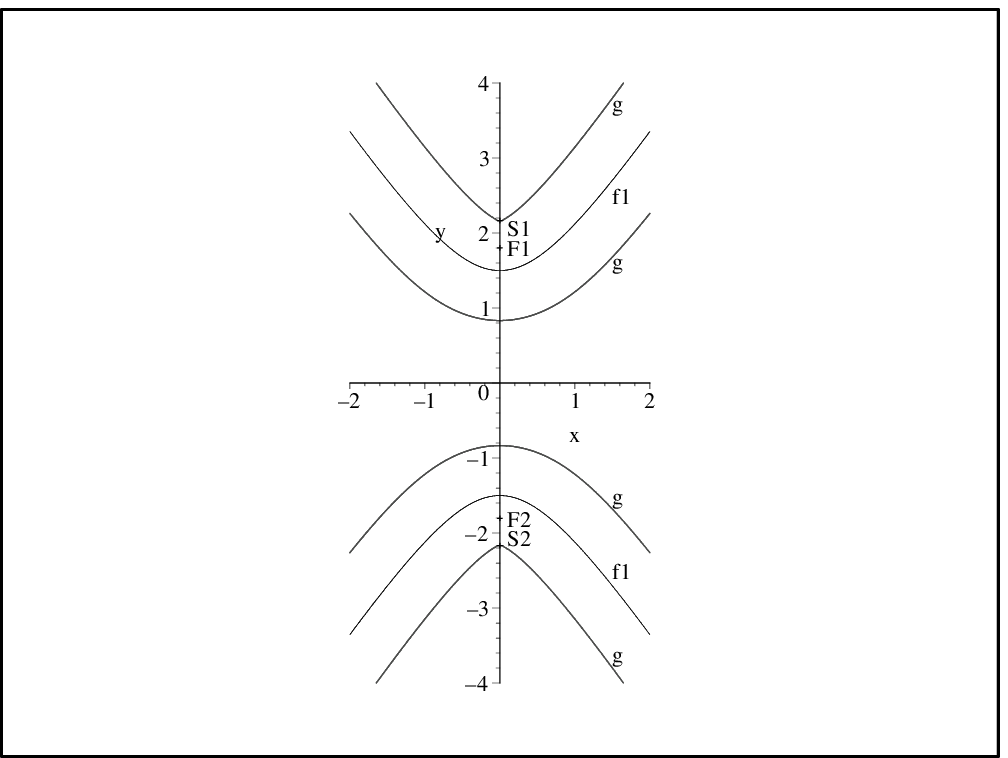}}}
\caption{Hyperbola with parallel lines $g,$ foci $F_1,F_2$ and singular points $S_1,S_2$ for $a=\frac32,$ $b=1$ and $r=\frac23$}
\label{fig:hyperbola2}
\end{figure}

\begin{example}[$r = r_{crit}$]
\label{example8}
Let $a=\frac32,$ $b=1,$ and $r=\frac23.$ Then,
$$
\bV=\bV(4 y_0^2-9 x_0^2-9, 9 y^2-18 y y_0+9 y_0^2+9 x^2-18 x x_0+9 x_0^2-4, 4 x y_0-13 x_0 y_0+9 x_0 y)
$$
We find that the reduced \Grobner basis $G_r$ contains eleven polynomials including 
\begin{setlength}{\multlinegap}{0pt}
\begin{multline}
g = 9447840 x^6 y^2-4199040 y^6 x^2+53800200 x^6+1679616 y^8+124857369 x^4\\[0.5ex]
    +127471968 y^4-159339960 x^2 y^2-24820992 y^6+46539360 y^4 x^2-105471720 y^2 x^4\\[0.5ex]
    -4933872 y^4 x^4+120670225+8503056 x^8-250879824 y^2+156799890 x^2 
\label{eq:gh2}
\end{multline}
\end{setlength}

\noindent
whose gradient is
\begin{align}
\begin{split}
\grad g &= x (1574640 y^2 x^4-233280 y^6+8966700 x^4+13873041 x^2-8852220 y^2+2585520 y^4\\[0.5ex]
        &\hspace {2in}  -11719080 x^2 y^2-548208 y^4 x^2+1889568 x^6+8711105) \bi \\[0.5ex]
        &\hspace*{0.18in} +y (131220 x^6-174960 y^4 x^2+93312 y^6+3540888 y^2-2213055 x^2-1034208 y^4\\[0.5ex]
        &\hspace {2in}  +1292760 x^2 y^2-1464885 x^4-137052 y^2 x^4-3484442)  \bj\,.
\end{split}
\end{align}
\noindent
There are exactly two solutions to the equation $(\grad g)(p) = 0$ on $\bV(g):$ 
$$
S_1 = \left(0, \frac{13}{6}\right) \quad \mbox{and} \quad  S_2 = \left(0, -\frac{13}{6}\right)
$$ 
These points are the {\rm singular points} of the variety $\bV(g):$ At every other point of this variety the gradient $\grad g$ does not vanish. We will collect these points as before into a non-singular subvariety $\Pi.$ The variety $\Pi$ constitutes the set of two parallel lines that we have been seeking. Note that the foci $F_1 = (0,\frac{\sqrt{13}}{2})$ and $F_2 = (0,-\frac{\sqrt{13}}{2})$ are still different from the singular points. We show all in Figure~\ref{fig:hyperbola2}.
\end{example}

As we can see from the last figure, when $r = r_{crit}= \frac23$ for the given values of $a$ and $b,$ the singular points are located on the 
continuous component of $\Pi.$

Probably the most interesting parallel lines we again obtain when $r >  r_{crit}$ as we show in the next example. There, we will find six singular points as each of the two singular points found in Example~\ref{example8} splits into three points.

\begin{figure}[htb]
\centerline{\scalebox{1.30}{\includegraphics{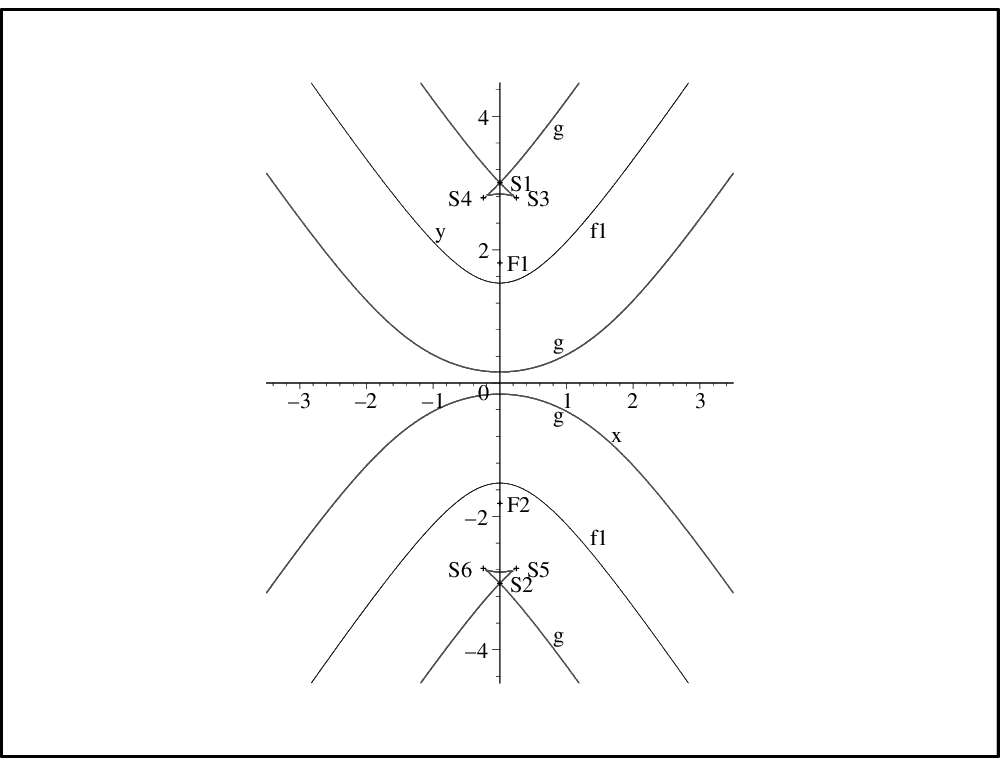}}}
\caption{Hyperbola with parallel lines $g,$ foci $F_1,F_2$ and six singular points $S_i,i=1,\ldots,6,$ for $a=\frac32,$ $b=1$ and $r=\frac43$}
\label{fig:hyperbola3}
\end{figure}

\begin{example}[$r > r_{crit}$]
\label{example9}
Let $a=\frac32,$ $b=1,$ and $r=\frac43.$ Then,
$$
\bV=\bV(4 y_0^2-9 x_0^2-9, 9 y^2-18 y y_0+9 y_0^2+9 x^2-18 x x_0+9 x_0^2-16, 4 x y_0-13 x_0 y_0+9 x_0 y)
$$
We find that the reduced \Grobner basis $G_r$ again contains nine polynomials including 
\begin{setlength}{\multlinegap}{0pt}
\begin{multline}
g = -1627128  x^6+8503056 x^8+30525625+83351034 x^2-1109487600 y^2+264529800 x^2 y^2\\[0.5ex]
    +51521913 x^4+381560544 y^4+349920 y^4 x^2-188052840 y^2 x^4-43856640 y^6-4933872 y^4 x^4\\[0.5ex]
    +9447840 x^6 y^2-4199040 y^6 x^2+1679616 y^8    
\label{eq:gh3}
\end{multline}
\end{setlength}
\noindent
whose gradient is
\begin{align}
\begin{split}
\grad g &=                  x (-271188 x^4+1889568 x^6+4630613+14696100 y^2+5724657 x^2+19440 y^4 \\[0.5ex]
        &\hspace*{1.0in}    -20894760 x^2 y^2-548208 y^4 x^2+1574640 y^2 x^4-233280 y^6) \bi \\[0.5ex]
        &\hspace*{0.18in}     +y (-15409550+3674025 x^2+10598904 y^2+9720 x^2 y^2-2611845 x^4-1827360 y^4 \\[0.5ex]
        &\hspace*{1.0in}    -137052 y^2 x^4+131220 x^6-174960 y^4 x^2+93312 y^6) \bj\,.
\end{split}
\end{align}


\noindent
There are exactly six solutions to the equation $(\grad g)(p) = 0$ on $\bV(g),$ or, in another words, six singular points of the variety:
\begin{gather*} 
S_1 = \left(0, \frac {5\,\sqrt{13}}{6}\right) \quad \mbox{and} \quad  S_2 = \left(0, -\frac {5\,\sqrt{13}}{6}\right)
\end{gather*}
\begin{align*} 
S_3 &= \left(\frac {2\,\sqrt{39 - 78 \cdot 2^{\frac13} + 39 \cdot 2^{\frac23}}}{13}, 
             \frac {\sqrt{12805 + 11232 \cdot 2^{\frac13} + 12636 \cdot 2^{\frac23}}}{78}\right),\\[0.5ex]
S_4 &= \left(-\frac {2\,\sqrt{39 - 78 \cdot 2^{\frac13} + 39 \cdot 2^{\frac23}}}{13}, 
             \frac {\sqrt{12805 + 11232 \cdot 2^{\frac13} + 12636 \cdot 2^{\frac23}}}{78}\right),\\[0.5ex]
S_5 &= \left(\frac {2\,\sqrt{39 - 78 \cdot 2^{\frac13} + 39 \cdot 2^{\frac23}}}{13}, 
             -\frac {\sqrt{12805 + 11232 \cdot 2^{\frac13} + 12636 \cdot 2^{\frac23}}}{78}\right),\\[0.5ex]
S_6 &= \left(-\frac {2\,\sqrt{39 - 78 \cdot 2^{\frac13} + 39 \cdot 2^{\frac23}}}{13}, 
             -\frac {\sqrt{12805 + 11232 \cdot 2^{\frac13} + 12636 \cdot 2^{\frac23}}}{78}\right).\\[0.5ex]
\end{align*}
We display these six points in Figure~\ref{fig:hyperbola3}.
\end{example}

\subsection{Parallel lines for arbitrary parameters $a,b,r$}
\label{hyperbola2}
In this section we will find a general polynomial $g$ that will generate the second elimination ideal $I_2 = I \cap \BR[x,y,a,b,r]$ where 
$I = \langle f_1, f_2, f_3\rangle$ and the generating polynomials are defined in equations (\ref{eq:fh1}), (\ref{eq:fh2}), and (\ref{eq:fh3}). We will also find a formula for the critical value $r_{crit}$ of $r$ that determines whether $\bV(I_2)$ has two or six singular points, and, we will find general formulas for the coordinates of the said singular points for arbitrary values of the parameters $a,b$ and $r.$  

We begin with a computation of a \Grobner basis for the ideal $I$ for the lexicographical order $y_0 > x_0 > x > y > a > b >r.$
$$
I = \langle b^2 y_0^2-a^2 x_0^2-a^2 b^2, y^2-2 y y_0+y_0^2+x^2-2 x x_0+x_0^2-r^2, b^2 y_0 x-b^2 y_0 x_0+a^2 x_0 y-a^2 x_0 y_0  \rangle
$$
We found that a reduced \Grobner basis for $I$ consists of fifteen homogeneous polynomials of degrees: 
$16$ (two polynomials),
$14$ (one polynomial),
$12$ (one polynomial),
$11$ (two polynomials),
$10$ (two polynomials),
$9$ (one polynomial),
$8$ (one polynomial),
$7$ (two polynomials),
$5$ (one polynomial),
$4$ (one polynomial),
and $2$ (one polynomial). We will not display these polynomials as they are similar to the polynomials for the ellipse shown in Appendix~B.

Thus, in this general case we observe that $I_2 = \langle g \rangle$ where $g = g_1/(a^2 b^2)$ is the following homogeneous polynomial in
$\BR[x,y,a,b,r]$ of degree twelve:

\beq
\begin{split}
g &= 6\,a^{4}\,x^{4}\,r^{4} - 2\,b^{2}\,x^{4}\,a^{2}\,r^{2}\,y^{2} + 4\,b^{6}\,x^{2}\,a^{4} + 2\,a^{8}\,b^{2}\,x^{2} + 2\,a^{6}\,x
^{6} + 6\,y^{4}\,a^{2}\,b^{6} - 2\,a^{4}\,b^{2}\,x^{4}\,y^{2} 
\\
& \mytabb + 4\,a^{4}\,x^{6}\,b^{2} - 6\,a^{6}\,x^{2}\,b^{2}\,y^{2} - 4\,a^{6}\,x^{2}\,b^{2}\,r^{2} - 10\,a^{4}\,x^{2}\,b^{4}\,y^{2}
 - 8\,a^{4}\,x^{2}\,b^{4}\,r^{2} 
\\
& \mytabb + 8\,a^{4}\,x^{2}\,b^{2}\,r^{4} - 6\,b^{4}\,r^{2}\,a^{2}\,x^{4} + 4\,b^{4}\,r^{4}\,a^{2}\,x^{2} - 10\,a^{4}\,x^{4}\,b^{2}
\,r^{2} - 6\,b^{4}\,r^{2}\,a^{2}\,x^{2}\,y^{2} 
\\
& \mytabb  + 6\,a^{4}\,x^{2}\,b^{2}\,r^{2}\,y^{2} - 6\,b^{4}\,x^{4}\,a^{2}\,y^{2} + 2\,b^{4}\,x^{2}\,y^{4}\,a^{2} - 6\,b^{6}\,x^{2}
\,r^{2}\,a^{2} + 10\,b^{4}\,y^{4}\,a^{2}\,r^{2} 
\\
& \mytabb  - 8\,b^{4}\,y^{2}\,a^{4}\,r^{2} - 4\,b^{6}\,y^{2}\,r^{2}\,a^{2} - 8\,b^{4}\,y^{2}\,r^{4}\,a^{2} + 6\,y^{4}\,a^{4}\,b^{2}
\,x^{2} - 2\,a^{2}\,y^{2}\,b^{8} + 6\,b^{4}\,y^{4}\,a^{4} 
\\
& \mytabb  - 4\,b^{4}\,y^{6}\,a^{2} - 6\,y^{2}\,a^{4}\,b^{6} + 6\,a^{4}\,x^{4}\,b^{4} - 4\,y^{2}\,a^{6}\,b^{4} - 6\,y^{2}\,a^{6}\,r
^{2}\,b^{2} + 6\,y^{4}\,a^{4}\,r^{2}\,b^{2} 
\\
& \mytabb  - 4\,y^{2}\,a^{4}\,r^{4}\,b^{2} + 2\,a^{6}\,b^{6} + a^{4}\,b^{8} + 6\,a^{6}\,b^{4}\,x^{2} + 6\,a^{6}\,x^{4}\,b^{2} + b^{8}
\,r^{4} + a^{8}\,b^{4} - 2\,a^{2}\,b^{8}\,r^{2} 
\\
& \mytabb  - 2\,a^{2}\,b^{6}\,r^{4} + 2\,b^{6}\,x^{2}\,r^{4} - 4\,b^{6}\,x^{2}\,y^{2}\,r^{2} + 6\,b^{2}\,x^{4}\,a^{2}\,r^{4} + 4\,a
^{6}\,x^{2}\,r^{2}\,y^{2} - 2\,b^{2}\,x^{2}\,a^{2}\,y^{6} 
\\
& \mytabb  + 2\,b^{4}\,x^{2}\,y^{6} - 2\,b^{2}\,x^{2}\,y^{4}\,a^{2}\,r^{2} + 6\,b^{4}\,x^{2}\,y^{2}\,r^{4} + 10\,b^{2}\,x^{2}\,y^{2}
\,a^{2}\,r^{4} + 2\,b^{6}\,x^{2}\,y^{4} 
\\
& \mytabb  - 6\,b^{4}\,x^{2}\,y^{4}\,r^{2} + a^{4}\,x^{8} + 2\,b^{6}\,r^{6} + b^{4}\,r^{8} + 6\,a^{6}\,x^{2}\,r^{4} - 6\,a^{6}\,x^{4}
\,r^{2} - 2\,a^{8}\,x^{2}\,r^{2} 
\\
& \mytabb  - 6\,a^{4}\,b^{4}\,r^{4} - 2\,a^{4}\,b^{6}\,r^{2} - 2\,a^{6}\,y^{2}\,r^{4} - 2\,a^{6}\,b^{2}\,r^{4} - 6\,a^{2}\,r^{6}\,b
^{2}\,y^{2} - 2\,a^{4}\,r^{6}\,y^{2} 
\\
& \mytabb  - 6\,a^{4}\,x^{4}\,r^{2}\,y^{2} + 6\,a^{2}\,y^{4}\,b^{2}\,r^{4} + 2\,a^{8}\,b^{2}\,r^{2} + a^{4}\,y^{4}\,r^{4} - 2\,a^{2}
\,y^{6}\,b^{2}\,r^{2} - 4\,a^{4}\,r^{6}\,x^{2} 
\\
& \mytabb  + a^{8}\,x^{4} + 6\,a^{4}\,x^{2}\,y^{2}\,r^{4} - 4\,a^{4}\,x^{6}\,r^{2} + a^{8}\,r^{4} + 2\,a^{4}\,x^{6}\,y^{2} - 2\,a^{6}
\,x^{4}\,y^{2} - 2\,a^{6}\,r^{6} 
\\
& \mytabb  - 2\,a^{4}\,r^{6}\,b^{2} + 2\,a^{6}\,r^{2}\,b^{4} + a^{4}\,y^{4}\,x^{4} + 2\,a^{2}\,r^{8}\,b^{2} - 2\,a^{4}\,y^{4}\,x^{2}
\,r^{2} + 2\,a^{2}\,b^{4}\,r^{6} + a^{4}\,r^{8} 
\\
& \mytabb  + 6\,y^{4}\,b^{4}\,r^{4} + 6\,y^{4}\,b^{6}\,r^{2} - 4\,y^{6}\,b^{4}\,r^{2} - 4\,y^{2}\,b^{4}\,r^{6} - 6\,y^{2}\,b^{6}\,r
^{4} + y^{4}\,b^{8} + y^{8}\,b^{4} - 2\,y^{6}\,b^{6} 
\\
& \mytabb  - 2\,b^{2}\,x^{6}\,a^{2}\,y^{2} - 2\,b^{4}\,x^{2}\,r^{6} - 6\,b^{2}\,x^{2}\,r^{6}\,a^{2} + b^{4}\,x^{4}\,y^{4} - 4\,b^{2}
\,x^{4}\,y^{4}\,a^{2} + b^{4}\,x^{4}\,r^{4} 
\\
& \mytabb  - 2\,b^{2}\,x^{6}\,a^{2}\,r^{2} - 2\,b^{4}\,x^{4}\,r^{2}\,y^{2} - 2\,y^{2}\,b^{8}\,r^{2} - 6\,y^{2}\,x^{2}\,a^{2}\,b^{6}
\end{split}
\label{eq:ghyperbola}
\eeq
We want to find now those points on the variety $\bV=\bV(g)$ where the gradient $\grad g$ is zero. Thus, $\bV$ can be thought of as union of points where the gradient is zero and those points where the gradient is not zero. Connected components of $\bV$ where the gradient is not zero constitute the $1$-dimensional manifold. We are interested in a discrete finite set $S$ of points where the gradient is zero. These points will give the singular points of the variety $V.$  We will show next that the critical value of $r$ is $r_{crit} = \frac{b^2}{a}.$ Therefore, we compute 
$\grad g = \grad g(x,y,a,b,r) = <h_1,h_2,h_3,h_4,h_5>$ where
\beq
\begin{split}
h_1 
&=  x(2\,b^{6}\,y^{2}\,r^{2} - a^{8}\,b^{2} + b^{4}\,r^{6} - x^{2}\,b^{4}\,r^{4} + 2\,x^{2}\,b^{4}\,r^{2}\,y^{2} - 3\,b^{4}\,y^{2}\,r^{4} + 3\,b^{4}\,y^{4}\,r^{2} - b^{6}\,y^{4} 
\\
& \hspace*{0.5in}  - 3\,a^{6}\,b^{4} - a^{8}\,x^{2} - 2\,a^{4}\,b^{6} - 3\,a^{4}\,y^{4}\,b^{2} + 3\,y^{2}\,a^{6}\,b^{2} + 3\,y^{2}\,b^{6}\,a
^{2} + 5\,b^{4}\,y^{2}\,a^{4} 
\\
& \hspace*{0.5in}  - b^{4}\,y^{4}\,a^{2} - 6\,a^{4}\,x^{2}\,b^{4} + r^{2}\,a^{8} - b^{4}\,x^{2}\,y^{4} + 2\,a^{6}\,x^{2}\,y^{2} + 6\,b^{4}\,x
^{2}\,a^{2}\,y^{2} + 2\,a^{4}\,b^{2}\,x^{2}\,y^{2} 
\\
& \hspace*{0.5in}  - b^{4}\,y^{6} - b^{6}\,r^{4} + 10\,a^{4}\,x^{2}\,b^{2}\,r^{2} - 6\,a^{6}\,x^{2}\,b^{2} - 6\,a^{4}\,x^{4}\,b^{2} - 3\,a^{6
}\,x^{4} + 3\,b^{6}\,r^{2}\,a^{2} 
\\
& \hspace*{0.5in}  - 2\,b^{4}\,r^{4}\,a^{2} + 2\,r^{6}\,a^{4} + 3\,r^{6}\,b^{2}\,a^{2} - 3\,r^{4}\,a^{6} + 4\,r^{2}\,a^{4}\,b^{4} + 2\,b^{2}
\,r^{2}\,a^{2}\,x^{2}\,y^{2} 
\\
& \hspace*{0.5in}  + 6\,a^{4}\,x^{2}\,r^{2}\,y^{2} + 6\,a^{6}\,x^{2}\,r^{2} + 3\,b^{2}\,r^{2}\,a^{2}\,x^{4} - 6\,a^{4}\,x^{2}\,r^{4} + b^{2}
\,y^{4}\,a^{2}\,r^{2} - 5\,b^{2}\,y^{2}\,r^{4}\,a^{2} 
\\
& \hspace*{0.5in}  - 3\,y^{2}\,a^{4}\,r^{4} - 6\,b^{2}\,r^{4}\,a^{2}\,x^{2} + 3\,b^{2}\,x^{4}\,a^{2}\,y^{2} + 4\,b^{2}\,x^{2}\,y^{4}\,a^{2}
 - y^{4}\,a^{4}\,x^{2} - 3\,a^{4}\,x^{4}\,y^{2} 
\\
& \hspace*{0.5in}  + 6\,a^{4}\,x^{4}\,r^{2} - 2\,y^{2}\,a^{6}\,r^{2} + y^{4}\,a^{4}\,r^{2} + b^{2}\,y^{6}\,a^{2} + 3\,b^{4}\,y^{2}\,a^{2}\,r
^{2} - 3\,b^{2}\,y^{2}\,a^{4}\,r^{2} 
\\
& \hspace*{0.5in}  + 2\,b^{2}\,a^{6}\,r^{2} + 6\,b^{4}\,x^{2}\,a^{2}\,r^{2} - 4\,b^{2}\,a^{4}\,r^{4} - 2\,a^{4}\,x^{6}) 
\end{split}
\eeq
\beq
\begin{split}
h_2 
&=  y( - 6\,b^{6}\,y^{2}\,r^{2} + 2\,b^{4}\,r^{6} - 3\,x^{2}\,b^{4}\,r^{4} + 6\,x^{2}\,b^{4}\,r^{2}\,y^{2} - 6\,b^{4}\,y^{2}\,r^{4} + 6\,b^{4}\,y^{4}\,r^{2} + r^{2}\,b^{8} 
\\
& \hspace*{0.5in}   + b^{4}\,x^{4}\,r^{2} + 2\,x^{2}\,b^{6}\,r^{2} + 3\,b^{6}\,y^{4} + 2\,a^{6}\,b^{4} + 3\,a^{4}\,b^{6} - 6\,y^{2}\,b^{6}\,a^{2} - 6\,b^{4}\,y^{2}\,a^{4} 
\\
& \hspace*{0.5in}   + 6\,b^{4}\,y^{4}\,a^{2} + b^{2}\,x^{6}\,a^{2} + 3\,b^{6}\,x^{2}\,a^{2} + 3\,b^{4}\,a^{2}\,x^{4} + 5\,a^{4}\,x^{2}\,b^{4} + a^{2}\,b^{8} - 2\,y^{2}\,x^{2}\,b^{6} 
\\
& \hspace*{0.5in}   - 3\,b^{4}\,x^{2}\,y^{4} - b^{4}\,x^{4}\,y^{2} - 2\,b^{4}\,x^{2}\,a^{2}\,y^{2} - 6\,a^{4}\,b^{2}\,x^{2}\,y^{2} - 2\,b^{4}
\,y^{6} + 3\,b^{6}\,r^{4} 
\\
& \hspace*{0.5in}   - 3\,a^{4}\,x^{2}\,b^{2}\,r^{2} - y^{2}\,b^{8} + 3\,a^{6}\,x^{2}\,b^{2} + a^{4}\,x^{4}\,b^{2} + a^{6}\,x^{4} + 2\,b^{6}\,r
^{2}\,a^{2} + 4\,b^{4}\,r^{4}\,a^{2} + r^{6}\,a^{4} 
\\
& \hspace*{0.5in}   + 3\,r^{6}\,b^{2}\,a^{2} + r^{4}\,a^{6} + 4\,r^{2}\,a^{4}\,b^{4} + 2\,b^{2}\,r^{2}\,a^{2}\,x^{2}\,y^{2} + 2\,a^{4}\,x^{2}
\,r^{2}\,y^{2} - 2\,a^{6}\,x^{2}\,r^{2} 
\\
& \hspace*{0.5in}   + b^{2}\,r^{2}\,a^{2}\,x^{4} - 3\,a^{4}\,x^{2}\,r^{4} + 3\,b^{2}\,y^{4}\,a^{2}\,r^{2} - 6\,b^{2}\,y^{2}\,r^{4}\,a^{2} - y
^{2}\,a^{4}\,r^{4} - 5\,b^{2}\,r^{4}\,a^{2}\,x^{2} 
\\
& \hspace*{0.5in}   + 4\,b^{2}\,x^{4}\,a^{2}\,y^{2} + 3\,b^{2}\,x^{2}\,y^{4}\,a^{2} - a^{4}\,x^{4}\,y^{2} + 3\,a^{4}\,x^{4}\,r^{2} - 10\,b^{4
}\,y^{2}\,a^{2}\,r^{2} - 6\,b^{2}\,y^{2}\,a^{4}\,r^{2} 
\\
& \hspace*{0.5in}   + 3\,b^{2}\,a^{6}\,r^{2} + 3\,b^{4}\,x^{2}\,a^{2}\,r^{2} + 2\,b^{2}\,a^{4}\,r^{4} - a^{4}\,x^{6})
\end{split}
\eeq
\beq
\begin{split}
h_3 
&= 2\,b^{2}\,x^{4}\,y^{4} + b^{2}\,x^{6}\,y^{2} + b^{2}\,x^{2}\,y^{6} + 2\,b^{6}\,y^{2}\,r^{2} - b^{4}\,r^{6} - 2\,x^{2}\,b^{4}\,r^{4} + 3\,x^{2}\,b^{4}\,r^{2}\,y^{2} 
\\
& \hspace*{0.5in}   + 4\,b^{4}\,y^{2}\,r^{4} - 5\,b^{4}\,y^{4}\,r^{2} + r^{2}\,b^{8} + 3\,b^{4}\,x^{4}\,r^{2} + 3\,x^{2}\,b^{6}\,r^{2} - 3\,b^{6}\,y^{4} - 2\,a^{6}\,b^{4} - a^{2}\,x^{8} 
\\
& \hspace*{0.5in}   + b^{2}\,x^{4}\,r^{2}\,y^{2} + b^{2}\,x^{2}\,y^{4}\,r^{2} - 5\,b^{2}\,x^{2}\,y^{2}\,r^{4} - a^{2}\,y^{4}\,x^{4} - r^{8}\,b^{2} + 2\,a^{2}\,r^{6}\,y^{2} - a^{2}\,y^{4}\,r^{4} 
\\
& \hspace*{0.5in}   + 4\,a^{2}\,r^{6}\,x^{2} - 6\,a^{2}\,x^{4}\,r^{4} - 3\,b^{2}\,x^{4}\,r^{4} + 6\,a^{2}\,x^{4}\,r^{2}\,y^{2} - 3\,y^{4}\,b^{2}\,r^{4} + 3\,r^{6}\,b^{2}\,y^{2} + y^{6}\,b^{2}\,r^{2} 
\\
& \hspace*{0.5in}   - 6\,a^{2}\,x^{2}\,y^{2}\,r^{4} + 2\,a^{2}\,y^{4}\,x^{2}\,r^{2} + 3\,b^{2}\,x^{2}\,r^{6} + b^{2}\,x^{6}\,r^{2} + 4\,a^{2}
\,x^{6}\,r^{2} - 2\,a^{2}\,x^{6}\,y^{2} - r^{8}\,a^{2} 
\\
& \hspace*{0.5in}   - 3\,a^{4}\,b^{6} + 6\,y^{2}\,b^{6}\,a^{2} + 6\,b^{4}\,y^{2}\,a^{4} - 6\,b^{4}\,y^{4}\,a^{2} - 4\,b^{2}\,x^{6}\,a^{2} - 4\,b^{6}\,x^{2}\,a^{2} - 6\,b^{4}\,a^{2}\,x^{4} 
\\
& \hspace*{0.5in}   - 9\,a^{4}\,x^{2}\,b^{4} - a^{2}\,b^{8} + 3\,y^{2}\,x^{2}\,b^{6} - b^{4}\,x^{2}\,y^{4} + 3\,b^{4}\,x^{4}\,y^{2} + 10\,b^{4}\,x^{2}\,a^{2}\,y^{2} \\
& \hspace*{0.5in}   + 9\,a^{4}\,b^{2}\,x^{2}\,y^{2} + 2\,b^{4}\,y^{6} + b^{6}\,r^{4} + 6\,a^{4}\,x^{2}\,b^{2}\,r^{2} + y^{2}\,b^{8} - 4\,a^{6}
\,x^{2}\,b^{2} - 9\,a^{4}\,x^{4}\,b^{2} 
\\
& \hspace*{0.5in}   - 2\,a^{6}\,x^{4} + 2\,b^{6}\,r^{2}\,a^{2} + 6\,b^{4}\,r^{4}\,a^{2} + 3\,r^{6}\,a^{4} + 2\,r^{6}\,b^{2}\,a^{2} - 2\,r^{4}
\,a^{6} - 3\,r^{2}\,a^{4}\,b^{4} 
\\
& \hspace*{0.5in}   - 6\,b^{2}\,r^{2}\,a^{2}\,x^{2}\,y^{2} - 6\,a^{4}\,x^{2}\,r^{2}\,y^{2} + 4\,a^{6}\,x^{2}\,r^{2} + 10\,b^{2}\,r^{2}\,a^{2}
\,x^{4} - 9\,a^{4}\,x^{2}\,r^{4} 
\\
& \hspace*{0.5in}   - 6\,b^{2}\,y^{4}\,a^{2}\,r^{2} + 4\,b^{2}\,y^{2}\,r^{4}\,a^{2} + 3\,y^{2}\,a^{4}\,r^{4} - 8\,b^{2}\,r^{4}\,a^{2}\,x^{2}
 + 2\,b^{2}\,x^{4}\,a^{2}\,y^{2} 
\\
& \hspace*{0.5in}   - 6\,b^{2}\,x^{2}\,y^{4}\,a^{2} + 3\,a^{4}\,x^{4}\,y^{2} + 9\,a^{4}\,x^{4}\,r^{2} + 8\,b^{4}\,y^{2}\,a^{2}\,r^{2} + 9\,b
^{2}\,y^{2}\,a^{4}\,r^{2} - 4\,b^{2}\,a^{6}\,r^{2} 
\\
& \hspace*{0.5in}   + 8\,b^{4}\,x^{2}\,a^{2}\,r^{2} + 3\,b^{2}\,a^{4}\,r^{4} - 3\,a^{4}\,x^{6} 
\end{split}
\eeq
\beq
\begin{split}
h_4 
&= b^{2}\,x^{4}\,y^{4} + 2\,b^{2}\,x^{2}\,y^{6} - 4\,b^{6}\,y^{2}\,r^{2} + a^{8}\,b^{2} + 3\,b^{4}\,r^{6} + 3\,x^{2}\,b^{4}\,r^{4} - 6\,x^{2}\,b^{4}\,r^{2}\,y^{2} 
\\
& \hspace*{0.5in} - 9\,b^{4}\,y^{2}\,r^{4} + 9\,b^{4}\,y^{4}\,r^{2} + 2\,b^{6}\,y^{4} + 3\,a^{6}\,b^{4} + b^{2}\,y^{8} + a^{8}\,x^{2} - 2\,b^{2}\,x^{4}\,r^{2}\,y^{2} 
\\
& \hspace*{0.5in} - 6\,b^{2}\,x^{2}\,y^{4}\,r^{2} + 6\,b^{2}\,x^{2}\,y^{2}\,r^{4} - 2\,a^{2}\,y^{4}\,x^{4} + r^{8}\,b^{2} - 3\,a^{2}\,r^{6}\,y^{2} + 3\,a^{2}\,y^{4}\,r^{4} - 3\,a^{2}\,r^{6}\,x^{2} 
\\
& \hspace*{0.5in} + 3\,a^{2}\,x^{4}\,r^{4} + b^{2}\,x^{4}\,r^{4} - a^{2}\,x^{4}\,r^{2}\,y^{2} + 6\,y^{4}\,b^{2}\,r^{4} - 4\,r^{6}\,b^{2}\,y^{2} - 4\,y^{6}\,b^{2}\,r^{2} + 5\,a^{2}\,x^{2}\,y^{2}\,r^{4} 
\\
& \hspace*{0.5in} - a^{2}\,y^{4}\,x^{2}\,r^{2} - 2\,b^{2}\,x^{2}\,r^{6} - a^{2}\,x^{6}\,r^{2} - a^{2}\,x^{6}\,y^{2} + r^{8}\,a^{2} - x^{2}\,
a^{2}\,y^{6} - a^{2}\,y^{6}\,r^{2} + 2\,a^{4}\,b^{6} 
\\
& \hspace*{0.5in} + 6\,a^{4}\,y^{4}\,b^{2} - 4\,y^{2}\,a^{6}\,b^{2} - 4\,y^{2}\,b^{6}\,a^{2} - 9\,b^{4}\,y^{2}\,a^{4} + 9\,b^{4}\,y^{4}\,a^{2} + 6\,a^{4}\,x^{2}\,b^{4} + r^{2}\,a^{8} 
\\
& \hspace*{0.5in} + 3\,b^{4}\,x^{2}\,y^{4} - 3\,a^{6}\,x^{2}\,y^{2} - 9\,b^{4}\,x^{2}\,a^{2}\,y^{2} - 10\,a^{4}\,b^{2}\,x^{2}\,y^{2} - 3\,b^{4}\,y^{6} + 2\,b^{6}\,r^{4} 
\\
& \hspace*{0.5in} - 8\,a^{4}\,x^{2}\,b^{2}\,r^{2} + 6\,a^{6}\,x^{2}\,b^{2} + 6\,a^{4}\,x^{4}\,b^{2} + 3\,a^{6}\,x^{4} - 4\,b^{6}\,r^{2}\,a^{2} - 3\,b^{4}\,r^{4}\,a^{2} - r^{6}\,a^{4} 
\\
& \hspace*{0.5in} + 2\,r^{6}\,b^{2}\,a^{2} - r^{4}\,a^{6} - 3\,r^{2}\,a^{4}\,b^{4} - 6\,b^{2}\,r^{2}\,a^{2}\,x^{2}\,y^{2} + 3\,a^{4}\,x^{2}\,r^{2}\,y^{2} - 2\,a^{6}\,x^{2}\,r^{2} 
\\
& \hspace*{0.5in}  - 6\,b^{2}\,r^{2}\,a^{2}\,x^{4} + 4\,a^{4}\,x^{2}\,r^{4} + 10\,b^{2}\,y^{4}\,a^{2}\,r^{2} - 8\,b^{2}\,y^{2}\,r^{4}\,a^{2} - 2\,y^{2}\,a^{4}\,r^{4} + 4\,b^{2}\,r^{4}\,a^{2}\,x^{2} 
\\
& \hspace*{0.5in}  - 6\,b^{2}\,x^{4}\,a^{2}\,y^{2} + 2\,b^{2}\,x^{2}\,y^{4}\,a^{2} + 3\,y^{4}\,a^{4}\,x^{2} - a^{4}\,x^{4}\,y^{2} - 5\,a^{4}\,x^{4}\,r^{2} - 3\,y^{2}\,a^{6}\,r^{2} 
\\
& \hspace*{0.5in}  + 3\,y^{4}\,a^{4}\,r^{2} - 4\,b^{2}\,y^{6}\,a^{2} - 6\,b^{4}\,y^{2}\,a^{2}\,r^{2} - 8\,b^{2}\,y^{2}\,a^{4}\,r^{2} + 2\,b^{2}\,a^{6}\,r^{2} - 9\,b^{4}\,x^{2}\,a^{2}\,r^{2} 
\\
& \hspace*{0.5in}  - 6\,b^{2}\,a^{4}\,r^{4} + 2\,a^{4}\,x^{6}  
\end{split}
\eeq
\beq
\begin{split}
h_5 
&= - 6\,b^{6}\,y^{2}\,r^{2} + a^{8}\,b^{2} + 2\,b^{4}\,r^{6} - 3\,x^{2}\,b^{4}\,r^{4} + 6\,x^{2}\,b^{4}\,r^{2}\,y^{2} - 6\,b^{4}\,y^{2}\,r^{4} + 6\,b^{4}\,y^{4}\,r^{2} + r^{2}\,b^{8} 
\\
& \hspace*{0.5in}   + b^{4}\,x^{4}\,r^{2} + 2\,x^{2}\,b^{6}\,r^{2} + 3\,b^{6}\,y^{4} + a^{6}\,b^{4} - a^{8}\,x^{2} - a^{4}\,b^{6} + 3\,a^{4}\,y^{4}\,b^{2} - 3\,y^{2}\,a^{6}\,b^{2} 
\\
& \hspace*{0.5in}   - 2\,y^{2}\,b^{6}\,a^{2} - 4\,b^{4}\,y^{2}\,a^{4} + 5\,b^{4}\,y^{4}\,a^{2} - b^{2}\,x^{6}\,a^{2} - 3\,b^{6}\,x^{2}\,a^{2} - 3\,b^{4}\,a^{2}\,x^{4} - 4\,a^{4}\,x^{2}\,b^{4} 
\\
& \hspace*{0.5in}   + r^{2}\,a^{8} - a^{2}\,b^{8} - 2\,y^{2}\,x^{2}\,b^{6} - 3\,b^{4}\,x^{2}\,y^{4} - b^{4}\,x^{4}\,y^{2} + 2\,a^{6}\,x^{2}\,y^{2} - 3\,b^{4}\,x^{2}\,a^{2}\,y^{2} 
\\
& \hspace*{0.5in}   + 3\,a^{4}\,b^{2}\,x^{2}\,y^{2} - 2\,b^{4}\,y^{6} + 3\,b^{6}\,r^{4} + 8\,a^{4}\,x^{2}\,b^{2}\,r^{2} - y^{2}\,b^{8} - 2\,a^{6}\,x^{2}\,b^{2} - 5\,a^{4}\,x^{4}\,b^{2} 
\\
& \hspace*{0.5in}   - 3\,a^{6}\,x^{4} - 2\,b^{6}\,r^{2}\,a^{2} + 3\,b^{4}\,r^{4}\,a^{2} + 2\,r^{6}\,a^{4} + 4\,r^{6}\,b^{2}\,a^{2} - 3\,r^{4}\,a^{6} - 6\,r^{2}\,a^{4}\,b^{4} 
\\
& \hspace*{0.5in}   + 10\,b^{2}\,r^{2}\,a^{2}\,x^{2}\,y^{2} + 6\,a^{4}\,x^{2}\,r^{2}\,y^{2} + 6\,a^{6}\,x^{2}\,r^{2} + 6\,b^{2}\,r^{2}\,a^{2}\,x^{4} - 6\,a^{4}\,x^{2}\,r^{4} 
\\
& \hspace*{0.5in}   + 6\,b^{2}\,y^{4}\,a^{2}\,r^{2} - 9\,b^{2}\,y^{2}\,r^{4}\,a^{2} - 3\,y^{2}\,a^{4}\,r^{4} - 9\,b^{2}\,r^{4}\,a^{2}\,x^{2} - b^{2}\,x^{4}\,a^{2}\,y^{2} - b^{2}\,x^{2}\,y^{4}\,a^{2} 
\\
& \hspace*{0.5in}   - y^{4}\,a^{4}\,x^{2} - 3\,a^{4}\,x^{4}\,y^{2} + 6\,a^{4}\,x^{4}\,r^{2} - 2\,y^{2}\,a^{6}\,r^{2} + y^{4}\,a^{4}\,r^{2} - b^{2}\,y^{6}\,a^{2} - 8\,b^{4}\,y^{2}\,a^{2}\,r^{2} 
\\
& \hspace*{0.5in}   - 4\,b^{2}\,y^{2}\,a^{4}\,r^{2} - 2\,b^{2}\,a^{6}\,r^{2} + 4\,b^{4}\,x^{2}\,a^{2}\,r^{2} - 3\,b^{2}\,a^{4}\,r^{4} - 2\,a^{4}\,x^{6}
\end{split}
\eeq
Continuing our approach from Section~\ref{ellipse}, in order to find the solutions to the vector equation $\grad g=0,$ we need to solve a system of six polynomial equations $h_i=0,i=1,\ldots,5$ and $g=0$ in the variables $x,y,a,b,r$ under the assumption that $a>b>r>0.$ Let 
$I_S = \langle \grad g,g\rangle = \langle h_1, h_2, h_3, h_4, h_5, g \rangle$ and let $F$ be a list containing the six polynomial equations:
\beq
F = [h_1 = 0, h_2 = 0, h_3 = 0, h_4 = 0, h_5 = 0, g = 0]
\label{eq:Fhyperbola} 
\eeq
In order to solve system (\ref{eq:Fhyperbola}), we can again either employ the \Grobner basis approach similar to the one shown in Appendix~C in case of the ellipse, or simply use {\tt solve} command in {\tt Maple}. Thus we find that the affine variety $\bV(I_S)$ consists of the following eight singular points of which six are real and two are complex:
\beq
S_1 = \left( 0, \,{\displaystyle \frac {\sqrt{(a^{2} + b^{2})\,(r^{2} + b^{2})}}{b}} \right), \qquad
S_2 = \left( 0, \,- {\displaystyle \frac {\sqrt{(a^{2} + b^{2})\,(r^{2} + b^{2})}}{b}}  \right),
\label{eq:SS12}
\eeq
\beq
S_3 = \left( {\displaystyle \frac {\sqrt{(a^{2} + b^{2})\,( - a^{2} + r^{2})}}{a}}  , \,0 \right), \qquad 
S_4 = \left( - {\displaystyle \frac {\sqrt{(a^{2} + b^{2})\,( - a^{2} + r^{2})}}{a}}  , \,0 \right),
\label{eq:SS34}
\eeq
\beq
S_5 = \left( \sqrt{-\frac{\delta}{(a^2+b^2)(r\,b)^{\frac23}}} ,\,-\frac{\sqrt{(r\,b)^{\frac13}(a^2+b^2)\varepsilon}}{(r\,b)^{\frac13}(a^2+b^2)}\right),
\label{eq:SS5}
\eeq
\beq
S_6 = \left( -\sqrt{-\frac{\delta}{(a^2+b^2)(r\,b)^{\frac23}}} ,\,-\frac{\sqrt{(r\,b)^{\frac13}(a^2+b^2)\varepsilon}}{(r\,b)^{\frac13}(a^2+b^2)}\right),
\label{eq:SS6}
\eeq
\beq
S_7 = \left( \sqrt{-\frac{\delta}{(a^2+b^2)(r\,b)^{\frac23}}} ,\, \frac{\sqrt{(r\,b)^{\frac13}(a^2+b^2)\varepsilon}}{(r\,b)^{\frac13}(a^2+b^2)}\right),
\label{eq:SS7}
\eeq
\beq
S_8 = \left( -\sqrt{-\frac{\delta}{(a^2+b^2)(r\,b)^{\frac23}}} ,\, \frac{\sqrt{(r\,b)^{\frac13}(a^2+b^2)\varepsilon}}{(r\,b)^{\frac13}(a^2+b^2)}\right),
\label{eq:SS8}
\eeq
where
\begin{align}
\delta   &=  3\,b^{2}\,a^{\frac43}\,r^{2} - r^{2}\,a^{2}\,(r\,b)^{\frac23} - 3\,a^{\frac23}\,(r\,b)^{(1/3)}\,b^{3}\,r + b^{4}\,(r\,b)^{\frac23},
\label{eq:deltah}\\[0.5ex]
\varepsilon &= a^{4}\,(r\,b)^{\frac13} + b^{2}\,r^{2}\,(r\,b)^{\frac13} + 3\,b\,a^{\frac43}\,r\,(r\,b)^{\frac23} + 3\,a^{\frac83}\,b\,r.
\label{eq:epsilonh}
\end{align}
Notice that $\delta$ becomes zero in three cases:
\begin{enumerate}
\item When $r = 0$ -- this is not possible since we expect $r > 0,$

\item When $r = - \frac{b^2}{a}$  -- this is not possible as $r > 0,$

\item When $r = \frac{b^2}{a}.$
\end{enumerate}
Only this last value is possible thus we set $r_{crit} = \frac{b^2}{a}.$ When $r = r_{crit},$ then under the assumptions that $a>b>r>0,$ we essentially obtain two real singular points:
\beq
S_1 = S_7 = S_8 = \left( 0,\frac{a^2+b^2}{a}\right) = \left(0,\frac{c^2}{a} \right)\quad \mbox{and} \quad 
S_2 = S_5 = S_6 = \left(0,-\frac{a^2+b^2}{a}\right) = \left(0,-\frac{c^2}{a}\right) ,
\label{eq:SS178}
\eeq
while the complex points $S_3$ and $S_4$ become:
\beq
S_3 = \left( \frac{a^2+b^2}{a^2}\sqrt{a^2-b^2}\,i,0 \right) \quad \mbox{and} \quad 
S_4 = \left(- \frac{a^2+b^2}{a^2}\sqrt{a^2-b^2}\,i,0 \right).
\eeq
For a hyperbola with $a > b,$ the ratio $l = \frac{b^2}{a}$ is called a {\it semi-latus rectum\/.} It is the distance between a focus and the hyperbola itself measured along a line perpendicular to the major axis. The value of  $l$ is also the reciprocal of the maximum curvature $\kappa_{max}$ of the hyperbola, or, it is equal to the radius of the smallest osculating circle.~\cite{wiki}

To summarize, the variety $\bV(I_S)$ contains the following real singular points:
\begin{enumerate}
\item Two virtual points $S_1$ and $S_2$ shown in (\ref{eq:SS12}) when $0 < r < r_{crit}.$ In particular, when $a=1,b=1,$ and $r=\frac12,$ upon substituting these values into (\ref{eq:ghyperbola}) we obtain the parallel lines (\ref{eq:gh7}) derived in Example~\ref{example7}. Furthermore, when these values are substituted into (\ref{eq:SS12}), then we obtain the two virtual singular points $S_1$ and $S_2$ found in the same example. The remaining six singular points are complex.
\item Two points $S_1$ and $S_2$ shown in (\ref{eq:SS12}) when $r = r_{crit}.$ In particular, we can recover results of Example~\ref{example8} when we substitute values $a=\frac32,$ $b=1,$ and $r=\frac23$ into the general equation (\ref{eq:ghyperbola}) and the singular points (\ref{eq:SS178}). To be precise, three of the singular points give $S_1,$ another three give $S_2$ whereas the remaining two gave  two distinct complex points. 
\item Six points $S_1,S_2,S_5,S_6,S_7,S_8$ shown in (\ref{eq:SS12}), (\ref{eq:SS5}), (\ref{eq:SS6}), (\ref{eq:SS7}), and (\ref{eq:SS8}) when 
$r>r_{crit}.$ In particular we can recover in the same manner six points and the equation for the parallel lines found in Example~\ref{example9} when we set $a=\frac32,$ $b=1,$ and $r=\frac43.$ These six singular points are displayed in Figure~\ref{fig:hyperbola3}. The remaining two points are complex and distinct.
\end{enumerate}

\section{Applications and Conclusions}
\label{applications}

There are several immediate applications for the results obtained in the previous sections as well as for the general method of solving systems of polynomial and rational equations using the \Grobner basis technique. In particular, we refer to the explicit form of the polynomial $g$ that defines the parallel lines to the conics given in (\ref{eq:gparabola}), (\ref{eq:gellipse}), and (\ref{eq:ghyperbola}). While finding $g$ amounts to computing the reduced \Grobner bases for the appropriate ideals and then applying the elimination theory -- both being standard techniques of algebraic geometry 
\cite{Cox}, \cite{Cox2}, \cite{Eisenbud}, \cite{Hartshorne}, \cite{Reid} -- these techniques are not used and probably not even known in the area of engineering applications where numeric, CAD generated curves, dominate. However, it does not mean that the fact that the critical value $r_{crit}$ of $r$ equals the reciprocal of the maximum curvature of the curve -- as derived above for the non-generate conics -- is not known in engineering: In fact, this fact is very well known in the theory of mechanisms when designing cams with reciprocating roller follower \cite[Chapter 8]{Norton},  
\cite[Section 5.10]{Uicker}. 

In the theory of machines, the ``inner'' parallel line (on the concave side of the curve) is referred to as the {\em cam profile}, while the path followed by a center of a roller -- in our case it was always given above as the polynomial $f_1$ -- is referred to as the 
{\em pitch curve} \cite{Norton}, \cite{Uicker}. The pitch curve is usually a smooth curve defined piece-wise by the pre-determined motion of the follower that could be harmonic, or determined by a polynomial of degree eight (the so called {\em displacement function}). This results in more complicated -- than in our examples above -- cam contours that may have portions which are concave, convex or flat \cite[page 404]{Norton}.  The problem of the inner curve being ``non smooth'', or, in our language, having singularity points on the smooth portion of the variety $\bV(g),$ is referred to as {\em undercutting} -- when the radius $r$ of the roller follower is greater than $r_{crit}.$ When $r = r_{crit},$ that is when the radius of the roller follower equals the minimum positive (convex) radius of the pitch curve, a cutter cutting out the cam will create a {\em cusp} point on the cam surface. In this latter case, while the cam contour is a continuous curve, it is not smooth at one or more points that results in a cam not running well at high speeds. In the first case when $r > r_{crit},$ the cutter cutting the cam out {\em (...) undercuts or removes material needed for cam contours in different locations and also creates a sharp point or cusp on the cam surface. This cam no longer has the same displacement function you so carefully designed} \cite[page 404]{Norton}. Such cams are not acceptable as the radius of the roller follower is too large. As a rule of thumb, the absolute value of the minimum radius of curvature of the cam pitch curve is preferably at least $2$ or $3$ time as large as the radius of the roller follower $r.$  These situations that correspond to our cases discussed in examples~\ref{example2} and \ref{example3} for the parabola,  \ref{example5} and \ref{example6} for the ellipse, and \ref{example8} and \ref{example9} for the hyperbola, are depicted very well in \cite[Figure 8-45]{Norton} and in \cite[Figure 5-29]{Uicker}.  

What differs in the approaches to the problem of parallel lines presented above and designing cam profiles is that (i) in the former approach we obtain explicit formulas for the parallel lines given by the polynomial $g$ rather than obtain a numeric representation for it using some CAD software, (ii) in reality, cam profiles are more complicated as they are not determined by a single curve but are determined piecewise and, preferably, smooth-wise. Our approach allows for an analytic, hence exact, representation of the parallel lines, and it therefore affords an exact differential calculus to be applied to $g.$ There is no reason why in principle this approach could not be extended to piece-wise defined pitch curves. In particular, it allows for the exact computation of the coordinates of the singular points.

In the following, we briefly show one application of the result presented in Example~\ref{example4}. Finally, we remark that these methods extend to non conics such as a \Bezier cubic \cite{Bezier}.

\subsection{Parallel lines to an ellipse in the finite element method}

Let us consider an ellipse determined by a polynomial $f_1:$
\beq
f_1=4\,y_0^{2} + 16\,x_0^{2} - 64
\eeq
for the values of $a = 4$ and $b= 2.$ Thus, $c = 2\sqrt{3}$ and $r_{crit} = \frac{b^2}{a} = 1.$ We will generate a finite element mesh by choosing these values for the offset $r = 0.2,\, 0.4,\, 0.6.$ These numbers determine thicknesses of layers, or, the distances of the parallel lines to the ellipse. Number of rows in the mesh matrix is twice the number of values of $r$ values plus $1.$

\begin{remark}
The reason we are interested in generating several smooth parallel lines to the ellipse is that we want to apply the finite element analysis to composite materials that are built precisely from smooth layers of various materials. As long as the parallel lines are smooth, and this is guaranteed in the range $0 < r < r_{crit},$ the layers are smooth and touch at each point. This guarantees, in turn, an optimal performance of the composite. 
\end{remark}

\begin{figure}[htb]
\centerline{\scalebox{1.40}{\includegraphics{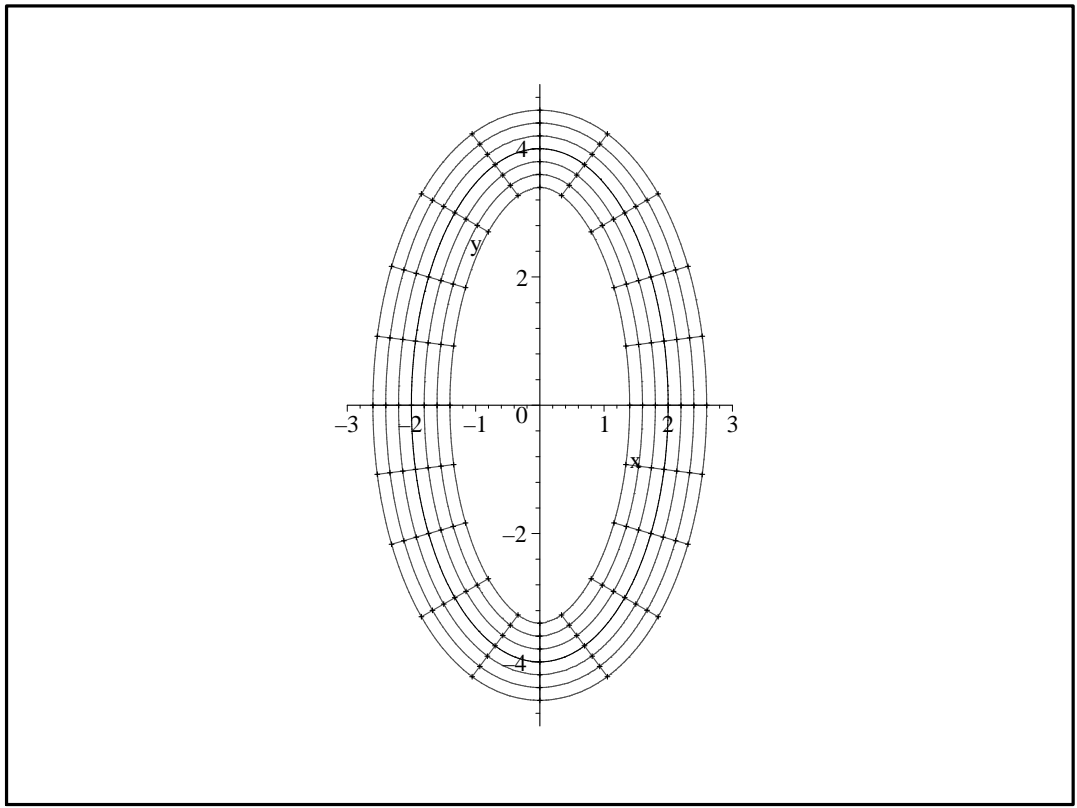}}}
\caption{Ellipse with parallel lines $g_1, g_2, g_3$ and nodes for $a=4,$ $b=2$ and $r=0.2,\,0.4,\,0.6$}
\label{fig:ellipsemesh}
\end{figure}

For each of the chosen values of $r,$ we can now compute the corresponding polynomial $g$ by substituting values of $a,b,r$ into (\ref{eq:gellipse}). Thus, we get:
\begin{align}
g_1 &=  - 0.1083369141 \cdot 10^{8}\,y^{2} - 0.1281707488 \cdot 10^{8}\,x^{2} + 0.4086271181 \cdot 10^{7}\,y^{2}\,x^{2} \notag \\
    & \mytabb + 33792\,y^{4}\,x^{4} + 10240\,x^{2}\,y^{6} + 40960\,x^{6}\,y^{2} - 495206.40\,x^{4}\,y^{2} - 369868.80\,x^{2}\,y^{4} \notag  \\
    & \mytabb + 259850.24\,x^{6} - 57180.16\,y^{6} + 0.1185870643 \cdot 10^{7}\,y^{4} - 564800.7168\,x^{4} + 1024\,y^{8} \notag  \\
    & \mytabb + 16384\,x^{8} + 0.3681277971 \cdot 10^{8} \qquad \mbox{when} \quad r = 0.2,\\[0.5ex]
g_2 &=  - 0.1031174671 \cdot 10^{8}\,y^{2} - 0.1349849501 \cdot 10^{8}\,x^{2} + 0.4156607693 \cdot 10^{7}\,y^{2}\,x^{2} \notag \\
    & \mytabb + 33792\,y^{4}\,x^{4} + 10240\,x^{2}\,y^{6} + 40960\,x^{6}\,y^{2} - 506265.60\,x^{4}\,y^{2} - 373555.20\,x^{2}\,y^{4} \notag  \\
    & \mytabb + 252968.96\,x^{6} - 56688.64\,y^{6} + 0.1155347251 \cdot 10^{7}\,y^{4} - 684903.6288\,x^{4} + 1024\,y^{8}  \notag \\
    & \mytabb + 16384\,x^{8} + 0.3409692932 \cdot 10^{8} \qquad \mbox{when} \quad r = 0.4,\\[0.5ex]
g_3 &=  - 0.9465872249 \cdot 10^{7}\,y^{2} - 0.1456495775 \cdot 10^{8}\,x^{2} + 0.4277898445 \cdot 10^{7}\,y^{2}\,x^{2} \notag \\
    & \mytabb + 33792\,y^{4}\,x^{4} + 10240\,x^{2}\,y^{6} + 40960\,x^{6}\,y^{2} - 524697.60\,x^{4}\,y^{2} - 379699.20\,x^{2}\,y^{4} \notag  \\
    & \mytabb + 241500.16\,x^{6} - 55869.44\,y^{6} + 0.1104343859 \cdot 10^{7}\,y^{4} - 880291.0208\,x^{4} \notag  \\
    & \mytabb + 0.2986886575 \cdot 10^{8} + 1024\,y^{8} + 16384\,x^{8}  \qquad \mbox{when} \quad r = 0.6.
\end{align}
Thus, so far we have rows in a mesh matrix. In order to create columns, we select some $y$ values that determine the $y$-coordinates of points on the ellipse through which we draw perpendicular lines to the ellipse. The intersections  of these lines with the ellipse itself as well as with the parallel lines are the mesh nodes. The number of these $y$-values gives the number of mesh columns minus $1.$ The last column is obtained by including intersection points between $g_1,\,g_2,\,g_3,$ and $f_1$ with the $y$-axis. For the sake of this example, we select these values: $y = 3.75, 3, 2, 1, 0, -1, -2, -3, -3.75$ (the $y$-values must satisfy $-a < y < a).$ We display the ellipse, the parallel lines, and the nodes in Figure~\ref{fig:ellipsemesh}.

It is easy now to produce a $7 \times 20$ mesh matrix with $140$ nodes and to create a list of $120$ $4$-node elements and a list of $30$ $9$-node elements. Once the elements have been created in {\tt Maple}, output has been sent to a finite element software for analysis of natural frequencies. It has been found, that there is a significant sensitivity in the higher mode of natural frequencies. For instance, when the mesh generated by exact parallel lines has been used, the higher mode of natural frequencies are much higher than the ones obtained from a mesh generated based on approximately parallel lines -- several confocal ellipses.

\appendix
\section{Appendix: Reduced \Grobner basis for the ideal\\ $I = \langle f_1,f_2,f_3 \rangle \subset\BR[y_0, x_0, x, y, r, p]$}
\label{AppendA}
For the ideal
$$
I = \langle 4 p y_0-x_0^2, y^2-2 y y_0+y_0^2+x^2-2 x x_0+x_0^2-r^2, 2 x p-2 x_0 p+x_0 y-x_0 y_0 \rangle \subset\BR[y_0, x_0, x, y, r, p]
$$
using lex order with $y_0 > x_0 > x > y > r > p,$ a reduced \Grobner basis for $I$ is given by the following fourteen polynomials:\footnote{We display these polynomials in a form returned by {\tt Maple}. Polynomials $g_1,g_2,g_3,g_4,g_6$ can be further simplified by dividing them by $p$ under the assumption that $p \neq 0.$}\\

\begin{smallmaplegroup}
\begin{mapleinput}
\mapleinline{active}{1d}{g1:=factor(reducedG[1]);}{%
}
\end{mapleinput}

\mapleresult
\begin{maplelatex}
\mapleinline{inert}{2d}{g1 :=
p*(-2*p*r^2*y*x^2+8*p*r^2*y^3+8*p^2*r^2*y^2-32*y*p^3*r^2+16*p^4*r^2-1
6*y^4*p^2+32*y^3*p^3-16*p^4*y^2+3*r^2*x^4+8*p^2*r^4+20*p^2*r^2*x^2-y^2
*x^4+10*y*p*x^4-x^6-x^4*p^2+8*p*y^3*x^2-32*x^2*y^2*p^2+8*x^2*y*p^3-3*r
^4*x^2+2*r^2*x^2*y^2+r^6-r^4*y^2-8*p*r^4*y);}{%
\maplemultiline{
\mathit{g1} :=  p( - 2\,p\,r^{2}\,y\,x^{2} + 8\,p\,r^{2}\,y^{3} + 8\,p^{2}\,r^{2}\,y^{2} - 32\,y\,p^{3}\,r^{2} + 16\,p^{4}\,r^{2} - 16\,y^{4}\,p^{2} + 32\,y^{3}\,p^{3} \\
\mytaba - 16\,p^{4}\,y^{2} + 3\,r^{2}\,x^{4} + 8\,p^{2}\,r^{4} + 20\,p^{2}\,r^{2}\,x^{2} - y^{2}\,x^{4} + 10\,y\,p\,x^{4} - x^{6}
 - x^{4}\,p^{2} + 8\,p\,y^{3}\,x^{2} \\
\mytaba  - 32\,x^{2}\,y^{2}\,p^{2} + 8\,x^{2}\,y\,p^{3} - 3\,r^{4}\,x^{2} + 2\,r^{2}\,x^{2}\,y^{2} + r^{6} - r^{4}\,y^{2} - 8\,p\,r^{4}\,y) }
}
\end{maplelatex}

\end{smallmaplegroup}
\begin{smallmaplegroup}
\begin{mapleinput}
\mapleinline{active}{1d}{g2:=factor(reducedG[2]);}{%
}
\end{mapleinput}

\mapleresult
\begin{maplelatex}
\mapleinline{inert}{2d}{g2 :=
p*(-114*x^3*p^2*y+232*x*p^3*y^2+40*x*p^4*y-316*x*p^3*r^2-8*x*y^3*p^2-2
*x^3*p^3+16*x*p^5-16*p*y^2*x*r^2-12*p*y^2*x0*r^2+54*x*y*p^2*r^2-120*y*
p^2*x0*r^2-32*y*p^4*x0+104*p^3*x0*r^2-14*p*y^2*x^3+32*p^2*y^3*x0+32*x*
p*r^4-46*x^3*p*r^2-16*p^5*x0+14*x^5*p+27*r^4*x0*p+2*y*x^5+2*y^3*x^3-2*
y^3*x*r^2-4*y^3*x0*r^2+2*y*x*r^4-4*y*x^3*r^2-24*p*y^4*x+16*p*y^4*x0);}
{%
\maplemultiline{
\mathit{g2} := p( - 114\,x^{3}\,p^{2}\,y + 232\,x\,p^{3}\,y^{2} + 40\,x\,p^{4}\,y - 316\,x\,p^{3}\,r^{2} - 8\,x\,y^{3}\,p^{2} - 
2\,x^{3}\,p^{3} + 16\,x\,p^{5} \\
\mytaba - 16\,p\,y^{2}\,x\,r^{2} - 12\,p\,y^{2}\,\mathit{x0}\,r^{2} + 54\,x\,y\,p^{2}\,r^{2} - 120\,y\,p^{2}\,\mathit{x0}\,r^{2}
 - 32\,y\,p^{4}\,\mathit{x0} \\
\mytaba + 104\,p^{3}\,\mathit{x0}\,r^{2} - 14\,p\,y^{2}\,x^{3} + 32\,p^{2}\,y^{3}\,\mathit{x0} + 32\,x\,p\,r^{4} - 46\,x^{3}\,p\,r
^{2} - 16\,p^{5}\,\mathit{x0} \\
\mytaba + 14\,x^{5}\,p + 27\,r^{4}\,\mathit{x0}\,p + 2\,y\,x^{5} + 2\,y^{3}\,x^{3} - 2\,y^{3}\,x\,r^{2} - 4\,y^{3}\,\mathit{x0}\,
r^{2} + 2\,y\,x\,r^{4} - 4\,y\,x^{3}\,r^{2} \\
\mytaba - 24\,p\,y^{4}\,x + 16\,p\,y^{4}\,\mathit{x0}) }
}
\end{maplelatex}

\end{smallmaplegroup}
\begin{smallmaplegroup}
\begin{mapleinput}
\mapleinline{active}{1d}{g3:=factor(reducedG[3]);}{%
}
\end{mapleinput}

\mapleresult
\begin{maplelatex}
\mapleinline{inert}{2d}{g3 :=
p*(4*x^2*p^3+56*p^3*r^2-4*x*p^3*x0-56*p^3*y^2+48*y^3*p^2+42*x^2*p^2*y
-12*x*p^2*y*x0-48*y*p^2*r^2-4*x^4*p+8*y^4*p+14*p*r^4-10*x^2*p*r^2-12*x
*p*y^2*x0-8*x^2*p*y^2+27*x*p*x0*r^2-22*y^2*p*r^2+2*r^4*y+2*y^3*x^2+2*y
*x^4-4*r^2*y*x^2-4*x*y^3*x0-2*r^2*y^3);}{%
\maplemultiline{
\mathit{g3} :=  p(4\,x^{2}\,p^{3} + 56\,p^{3}\,r^{2} - 4\,x\,p^{3}\,\mathit{x0} - 56\,p^{3}\,y^{2} + 48\,y^{3}\,p^{2} + 42\,x^{
2}\,p^{2}\,y - 12\,x\,p^{2}\,y\,\mathit{x0} \\
\mytaba - 48\,y\,p^{2}\,r^{2} - 4\,x^{4}\,p + 8\,y^{4}\,p + 14\,p\,r^{4} - 10\,x^{2}\,p\,r^{2} - 12\,x\,p\,y^{2}\,\mathit{x0} - 8
\,x^{2}\,p\,y^{2} \\
\mytaba + 27\,x\,p\,\mathit{x0}\,r^{2} - 22\,y^{2}\,p\,r^{2} + 2\,r^{4}\,y + 2\,y^{3}\,x^{2} + 2\,y\,x^{4} - 4\,r^{2}\,y\,x^{2}
 - 4\,x\,y^{3}\,\mathit{x0} - 2\,r^{2}\,y^{3}) }
}
\end{maplelatex}

\end{smallmaplegroup}
\begin{smallmaplegroup}
\begin{mapleinput}
\mapleinline{active}{1d}{g4:=factor(reducedG[4]);}{%
}
\end{mapleinput}

\mapleresult
\begin{maplelatex}
\mapleinline{inert}{2d}{g4 :=
p*(112*x*p^4-112*x0*p^4-120*y*x*p^3+176*y*x0*p^3-124*x*p^2*r^2-80*y^2
*x0*p^2-16*x0*p^2*r^2+160*x*p^2*y^2-96*x^2*p^2*x0+82*x^3*p^2-28*x^3*y*
p+4*y*x0*p*r^2-24*p*x*y^3+16*p*y^3*x0+22*x*y*p*r^2-4*r^2*y^2*x0+2*y^2*
x^3-4*r^2*x^3-2*r^2*x*y^2+2*r^4*x+3*r^4*x0-3*x^2*x0*r^2+2*x^5);}{%
\maplemultiline{
\mathit{g4} :=  p(112\,x\,p^{4} - 112\,\mathit{x0}\,p^{4} - 120
\,y\,x\,p^{3} + 176\,y\,\mathit{x0}\,p^{3} - 124\,x\,p^{2}\,r^{2} - 80\,y^{2}\,\mathit{x0}\,p^{2} \\
\mytaba - 16\,\mathit{x0}\,p^{2}\,r^{2} + 160\,x\,p^{2}\,y^{2} - 
96\,x^{2}\,p^{2}\,\mathit{x0} + 82\,x^{3}\,p^{2} - 28\,x^{3}\,y\,p + 4\,y\,\mathit{x0}\,p\,r^{2} \\
\mytaba - 24\,p\,x\,y^{3} + 16\,p\,y^{3}\,\mathit{x0} + 22\,x\,y\,p\,r^{2} - 4\,r^{2}\,y^{2}\,\mathit{x0} + 2\,y^{2}\,x^{3} - 4\,
r^{2}\,x^{3} - 2\,r^{2}\,x\,y^{2} \\
\mytaba + 2\,r^{4}\,x + 3\,r^{4}\,\mathit{x0} - 3\,x^{2}\,\mathit{x0}\,r^{2} + 2\,x^{5}) }
}
\end{maplelatex}

\end{smallmaplegroup}
\begin{smallmaplegroup}
\begin{mapleinput}
\mapleinline{active}{1d}{g5:=factor(reducedG[5]);}{%
}
\end{mapleinput}

\mapleresult
\begin{maplelatex}
\mapleinline{inert}{2d}{g5 :=
8*x*p*y^2-6*x^3*p+7*x0*p*x^2+6*x*p*r^2+4*p*x0*y^2+y*x0*x^2-y*x0*r^2+8
*y*x*p^2-12*y*x0*p^2+2*x0*p*r^2-8*x*p^3+8*x0*p^3;}{%
\maplemultiline{
\mathit{g5} :=  8\,x\,p\,y^{2} - 6\,x^{3}\,p + 7\,\mathit{x0}\,p\,x^{2} + 6\,x\,p\,r^{2} + 4\,p\,\mathit{x0}\,y^{2} + y\,
\mathit{x0}\,x^{2} - y\,\mathit{x0}\,r^{2} + 8\,y\,x\,p^{2} \\
\mytaba - 12\,y\,\mathit{x0}\,p^{2} + 2\,\mathit{x0}\,p\,r^{2} - 8\,x\,p^{3} + 8\,\mathit{x0}\,p^{3} }
}
\end{maplelatex}

\end{smallmaplegroup}
\begin{smallmaplegroup}
\begin{mapleinput}
\mapleinline{active}{1d}{g6:=factor(reducedG[6]);}{%
}
\end{mapleinput}

\mapleresult
\begin{maplelatex}
\mapleinline{inert}{2d}{g6 :=
p*(8*p^2*r^2-8*y^2*p^2-4*x0*x*p^2+4*x^2*p^2+4*x*y*x0*p+8*p*y^3+2*y*p*
x^2-8*y*p*r^2+3*x0*x*r^2-2*r^2*y^2+2*x^4+2*r^4+2*x^2*y^2-4*x^2*r^2-3*x
0*x^3-4*x*y^2*x0);}{%
\maplemultiline{
\mathit{g6} :=  p(8\,p^{2}\,r^{2} - 8\,y^{2}\,p^{2} - 4\,\mathit{x0}\,x\,p^{2} + 4\,x^{2}\,p^{2} + 4\,x\,y\,\mathit{x0}\,p
 + 8\,p\,y^{3} + 2\,y\,p\,x^{2} - 8\,y\,p\,r^{2} \\
\mytaba + 3\,\mathit{x0}\,x\,r^{2} - 2\,r^{2}\,y^{2} + 2\,x^{4}
 + 2\,r^{4} + 2\,x^{2}\,y^{2} - 4\,x^{2}\,r^{2} - 3\,\mathit{x0}\,x^{3} - 4\,x\,y^{2}\,\mathit{x0}) }
}
\end{maplelatex}

\end{smallmaplegroup}
\begin{smallmaplegroup}
\begin{mapleinput}
\mapleinline{active}{1d}{g7:=factor(reducedG[7]);}{%
}
\end{mapleinput}

\mapleresult
\begin{maplelatex}
\mapleinline{inert}{2d}{g7 :=
x0*x^4-2*x^2*x0*r^2+44*x^2*p^2*x0+32*y^2*x0*p^2-12*y*x0*p*r^2-80*y*x0*
p^3+r^4*x0+16*x0*p^2*r^2+48*x0*p^4-2*x^3*y*p-38*x^3*p^2-40*x*p^2*y^2+2
*x*y*p*r^2+56*y*x*p^3+20*x*p^2*r^2-48*x*p^4;}{%
\maplemultiline{
\mathit{g7} := \mathit{x0}\,x^{4} - 2\,x^{2}\,\mathit{x0}\,r^{2} + 44\,x^{2}\,p^{2}\,\mathit{x0} + 32\,y^{2}\,\mathit{x0}\,p^{2}
 - 12\,y\,\mathit{x0}\,p\,r^{2} - 80\,y\,\mathit{x0}\,p^{3} + r^{4}\,\mathit{x0} \\
\mytaba + 16\,\mathit{x0}\,p^{2}\,r^{2} + 48\,\mathit{x0}\,p^{4} - 2\,x^{3}\,y\,p - 38\,x^{3}\,p^{2} - 40\,x\,p^{2}\,y^{2} + 2\,x
\,y\,p\,r^{2} + 56\,y\,x\,p^{3} \\
\mytaba + 20\,x\,p^{2}\,r^{2} - 48\,x\,p^{4} }
}
\end{maplelatex}

\end{smallmaplegroup}
\begin{smallmaplegroup}
\begin{mapleinput}
\mapleinline{active}{1d}{g8:=factor(reducedG[8]);}{%
}
\end{mapleinput}

\mapleresult
\begin{maplelatex}
\mapleinline{inert}{2d}{g8 :=
3*x0^2*r^2-12*x0^2*p^2+16*x*y*x0*p+32*x0*x*p^2-20*x^2*p^2-48*y^2*p^2+4
8*p^2*r^2-2*x0*x^3+2*x0*x*r^2+4*y*p*x^2;}{%
\maplemultiline{
\mathit{g8} := 3\,\mathit{x0}^{2}\,r^{2} - 12\,\mathit{x0}^{2}\,p^{2} + 16\,x\,y\,\mathit{x0}\,p + 32\,\mathit{x0}\,x\,p^{2} - 20
\,x^{2}\,p^{2} - 48\,y^{2}\,p^{2} + 48\,p^{2}\,r^{2} \\
\mytaba - 2\,\mathit{x0}\,x^{3} + 2\,\mathit{x0}\,x\,r^{2} + 4\,y\,p\,x^{2} }
}
\end{maplelatex}

\end{smallmaplegroup}
\begin{smallmaplegroup}
\begin{mapleinput}
\mapleinline{active}{1d}{g9:=factor(reducedG[9]);}{%
}
\end{mapleinput}

\mapleresult
\begin{maplelatex}
\mapleinline{inert}{2d}{g9 := -4*p*y^2-4*p*x^2+6*p*x*x0-2*p*x0^2+4*p*r^2+x0^2*y;}{%
\[
\mathit{g9} :=  - 4\,p\,y^{2} - 4\,p\,x^{2} + 6\,p\,x\,\mathit{x0} - 2\,p\,\mathit{x0}^{2} + 4\,p\,r^{2} + \mathit{x0}^{2}\,y
\]
}
\end{maplelatex}

\end{smallmaplegroup}
\begin{smallmaplegroup}
\begin{mapleinput}
\mapleinline{active}{1d}{g10:=factor(reducedG[10]);}{%
}
\end{mapleinput}

\mapleresult
\begin{maplelatex}
\mapleinline{inert}{2d}{g10 :=
-2*x0*x^2+3*x0^2*x+2*x0*r^2+4*y*x*p-8*p*x0*y-8*x*p^2+8*x0*p^2;}{%
\[
\mathit{g10} :=  - 2\,\mathit{x0}\,x^{2} + 3\,\mathit{x0}^{2}\,x + 2\,\mathit{x0}\,r^{2} + 4\,y\,x\,p - 8\,p\,\mathit{x0}\,y - 8
\,x\,p^{2} + 8\,\mathit{x0}\,p^{2}
\]
}
\end{maplelatex}

\end{smallmaplegroup}
\begin{smallmaplegroup}
\begin{mapleinput}
\mapleinline{active}{1d}{g11:=factor(reducedG[11]);}{%
}
\end{mapleinput}

\mapleresult
\begin{maplelatex}
\mapleinline{inert}{2d}{g11 := x0^3-8*x*p^2+8*x0*p^2-4*p*x0*y;}{%
\[
\mathit{g11} := \mathit{x0}^{3} - 8\,x\,p^{2} + 8\,\mathit{x0}\,p^{2} - 4\,p\,\mathit{x0}\,y
\]
}
\end{maplelatex}

\end{smallmaplegroup}
\begin{smallmaplegroup}
\begin{mapleinput}
\mapleinline{active}{1d}{g12:=factor(reducedG[12]);}{%
}
\end{mapleinput}

\mapleresult
\begin{maplelatex}
\mapleinline{inert}{2d}{g12 := 4*p*y0-x0^2;}{%
\[
\mathit{g12} := 4\,p\,\mathit{y0} - \mathit{x0}^{2}
\]
}
\end{maplelatex}

\end{smallmaplegroup}
\begin{smallmaplegroup}
\begin{mapleinput}
\mapleinline{active}{1d}{g13:=factor(reducedG[13]);}{%
}
\end{mapleinput}

\mapleresult
\begin{maplelatex}
\mapleinline{inert}{2d}{g13 := -2*x*p+2*x0*p-x0*y+x0*y0;}{%
\[
\mathit{g13} :=  - 2\,x\,p + 2\,\mathit{x0}\,p - \mathit{x0}\,y + \mathit{x0}\,\mathit{y0}
\]
}
\end{maplelatex}

\end{smallmaplegroup}
\begin{smallmaplegroup}
\begin{mapleinput}
\mapleinline{active}{1d}{g14:=factor(reducedG[14]);}{%
}
\end{mapleinput}

\mapleresult
\begin{maplelatex}
\mapleinline{inert}{2d}{g14 := y^2-2*y*y0+y0^2+x^2-2*x*x0+x0^2-r^2;}{%
\[
\mathit{g14} := y^{2} - 2\,y\,\mathit{y0} + \mathit{y0}^{2} + x^{2} - 2\,x\,\mathit{x0} + \mathit{x0}^{2} - r^{2}
\]
}
\end{maplelatex}

\end{smallmaplegroup}

\section{Appendix: Reduced \Grobner basis for the ideal\\ $I = \langle f_1,f_2,f_3 \rangle \subset\BR[y_0, x_0, x, y, a, b, r]$}
\label{AppendB}
For the ideal
$$
I = \langle b^2 y_0^2+a^2 x_0^2-a^2 b^2, y^2-2 y y_0+y_0^2+x^2-2 x x_0+x_0^2-r^2, b^2 y_0 x-b^2 y_0 x_0-a^2 x_0 y+a^2 x_0 y_0 \rangle
\subset\BR[y_0, x_0, x, y, a, b, r]
$$
using lex order with $y_0 > x_0 > x > y > a> b> r,$ a reduced \Grobner basis for $I$ is given by the following fifteen polynomials:\footnote{We display these polynomials in a form returned by {\tt Maple}. All polynomials except $g_{13}$ and $g_{14}$ can be further simplified by dividing them by 
$a^2,b^2 $ under the assumption that $a,b \neq 0.$}\\

\begin{smallmaplegroup}
\begin{mapleinput}
\mapleinline{active}{1d}{g1:=factor(reducedG[1]);}{%
}
\end{mapleinput}

\mapleresult
\begin{maplelatex}
\mapleinline{inert}{2d}{g1 :=
a^2*b^2*(6*b^4*x^4*a^4-4*y^6*b^4*a^2-6*a^6*x^4*b^2-4*b^6*x^2*a^4+6*b^4
*x^2*a^6+a^4*b^8-2*y^2*b^8*a^2+6*y^2*b^6*x^2*a^2+6*a^6*y^2*x^2*b^2-10*
b^4*x^2*a^4*y^2+y^8*b^4+b^8*r^4+b^4*r^8-2*b^6*r^6-2*b^4*y^2*x^4*r^2-6*
x^2*y^4*b^4*r^2+6*x^2*b^4*r^4*y^2-4*b^4*y^6*r^2+2*b^4*y^6*x^2+y^4*b^4*
x^4+b^4*x^4*r^4-2*b^4*r^6*x^2-2*b^6*x^2*r^4+6*y^4*b^4*r^4-2*b^8*y^2*r^
2+6*b^6*r^4*y^2-4*b^4*r^6*y^2+6*a^4*y^2*b^6+6*a^6*y^2*b^2*r^2+2*a^4*b^
2*x^4*y^2-4*a^4*x^6*b^2+4*a^6*x^2*b^2*r^2-8*b^4*x^2*a^4*r^2+10*a^4*x^4
*b^2*r^2-6*b^4*r^2*a^2*x^2*y^2-6*a^4*x^2*b^2*r^2*y^2-8*a^4*x^2*b^2*r^4
-6*b^4*r^2*a^2*x^4+4*b^4*r^4*a^2*x^2-6*y^4*a^4*b^2*x^2-6*b^4*x^4*a^2*y
^2+6*b^6*x^2*r^2*a^2+2*b^4*x^2*y^4*a^2+10*b^4*y^4*a^2*r^2+4*b^6*y^2*r^
2*a^2-8*b^4*y^2*r^4*a^2-8*b^4*y^2*a^4*r^2+2*r^6*a^2*b^4+a^8*b^4-6*a^4*
y^4*r^2*b^2+4*a^4*y^2*r^4*b^2-2*a^8*b^2*x^2+x^4*a^4*y^4+4*x^2*b^6*y^2*
r^2-6*y^4*b^6*a^2-4*a^6*y^2*b^4-2*y^4*a^4*x^2*r^2+2*b^2*r^2*a^2*x^4*y^
2-2*r^8*a^2*b^2+r^8*a^4+a^8*x^4+a^4*x^8-4*a^4*x^6*r^2+6*a^4*x^4*r^4+2*
b^2*r^2*a^2*x^6-6*b^2*r^4*a^2*x^4-2*a^6*r^6+2*a^6*x^6+a^8*r^4+2*y^6*b^
2*a^2*x^2-2*r^6*a^4*y^2+6*r^6*b^2*y^2*a^2-2*a^8*x^2*r^2+y^4*a^4*r^4-6*
y^4*b^2*r^4*a^2-2*a^6*y^2*r^4+4*x^2*a^6*y^2*r^2+6*x^2*a^4*y^2*r^4+2*y^
4*b^2*x^2*a^2*r^2-10*x^2*b^2*y^2*r^4*a^2+2*y^6*b^2*a^2*r^2+6*r^6*b^2*x
^2*a^2+4*b^2*x^4*y^4*a^2+2*a^4*x^6*y^2-6*a^6*x^4*r^2-4*r^6*a^4*x^2+2*b
^2*x^6*a^2*y^2+6*r^4*a^6*x^2-6*a^4*x^4*r^2*y^2+2*b^6*y^6-2*a^8*b^2*r^2
-2*x^2*y^4*b^6+6*a^4*y^4*b^4+y^4*b^8-6*y^4*b^6*r^2-2*x^4*a^6*y^2-2*a^6
*b^6+2*b^6*r^4*a^2-2*b^8*r^2*a^2+2*b^4*a^6*r^2-6*b^4*a^4*r^4+2*r^2*a^4
*b^6+2*r^4*a^6*b^2+2*r^6*a^4*b^2);}{%
\maplemultiline{
\mathit{g1} := a^{2}\,b^{2}(6\,b^{4}\,x^{4}\,a^{4} - 4\,y^{6}\,b
^{4}\,a^{2} - 6\,a^{6}\,x^{4}\,b^{2} - 4\,b^{6}\,x^{2}\,a^{4} + 6
\,b^{4}\,x^{2}\,a^{6} + a^{4}\,b^{8} - 2\,y^{2}\,b^{8}\,a^{2} \\
\mytaba + 6\,y^{2}\,b^{6}\,x^{2}\,a^{2} + 6\,a^{6}\,y^{2}\,x^{2}
\,b^{2} - 10\,b^{4}\,x^{2}\,a^{4}\,y^{2} + y^{8}\,b^{4} + b^{8}\,
r^{4} + b^{4}\,r^{8} - 2\,b^{6}\,r^{6} \\
\mytaba  - 2\,b^{4}\,y^{2}\,x^{4}\,r^{2} - 6\,x^{2}\,y^{4}\,b^{4}
\,r^{2} + 6\,x^{2}\,b^{4}\,r^{4}\,y^{2} - 4\,b^{4}\,y^{6}\,r^{2}
 + 2\,b^{4}\,y^{6}\,x^{2} + y^{4}\,b^{4}\,x^{4} \\
\mytaba + b^{4}\,x^{4}\,r^{4} - 2\,b^{4}\,r^{6}\,x^{2} - 2\,b^{6}
\,x^{2}\,r^{4} + 6\,y^{4}\,b^{4}\,r^{4} - 2\,b^{8}\,y^{2}\,r^{2}
 + 6\,b^{6}\,r^{4}\,y^{2} - 4\,b^{4}\,r^{6}\,y^{2} \\
\mytaba + 6\,a^{4}\,y^{2}\,b^{6} + 6\,a^{6}\,y^{2}\,b^{2}\,r^{2}
 + 2\,a^{4}\,b^{2}\,x^{4}\,y^{2} - 4\,a^{4}\,x^{6}\,b^{2} + 4\,a
^{6}\,x^{2}\,b^{2}\,r^{2} - 8\,b^{4}\,x^{2}\,a^{4}\,r^{2} \\
\mytaba + 10\,a^{4}\,x^{4}\,b^{2}\,r^{2} - 6\,b^{4}\,r^{2}\,a^{2}
\,x^{2}\,y^{2} - 6\,a^{4}\,x^{2}\,b^{2}\,r^{2}\,y^{2} - 8\,a^{4}
\,x^{2}\,b^{2}\,r^{4} - 6\,b^{4}\,r^{2}\,a^{2}\,x^{4} \\
\mytaba + 4\,b^{4}\,r^{4}\,a^{2}\,x^{2} - 6\,y^{4}\,a^{4}\,b^{2}
\,x^{2} - 6\,b^{4}\,x^{4}\,a^{2}\,y^{2} + 6\,b^{6}\,x^{2}\,r^{2}
\,a^{2} + 2\,b^{4}\,x^{2}\,y^{4}\,a^{2} \\
\mytaba + 10\,b^{4}\,y^{4}\,a^{2}\,r^{2} + 4\,b^{6}\,y^{2}\,r^{2}
\,a^{2} - 8\,b^{4}\,y^{2}\,r^{4}\,a^{2} - 8\,b^{4}\,y^{2}\,a^{4}
\,r^{2} + 2\,r^{6}\,a^{2}\,b^{4} + a^{8}\,b^{4} \\
\mytaba - 6\,a^{4}\,y^{4}\,r^{2}\,b^{2} + 4\,a^{4}\,y^{2}\,r^{4}
\,b^{2} - 2\,a^{8}\,b^{2}\,x^{2} + x^{4}\,a^{4}\,y^{4} + 4\,x^{2}
\,b^{6}\,y^{2}\,r^{2} - 6\,y^{4}\,b^{6}\,a^{2} \\
\mytaba - 4\,a^{6}\,y^{2}\,b^{4} - 2\,y^{4}\,a^{4}\,x^{2}\,r^{2}
 + 2\,b^{2}\,r^{2}\,a^{2}\,x^{4}\,y^{2} - 2\,r^{8}\,a^{2}\,b^{2}
 + r^{8}\,a^{4} + a^{8}\,x^{4} + a^{4}\,x^{8} \\
\mytaba - 4\,a^{4}\,x^{6}\,r^{2} + 6\,a^{4}\,x^{4}\,r^{4} + 2\,b
^{2}\,r^{2}\,a^{2}\,x^{6} - 6\,b^{2}\,r^{4}\,a^{2}\,x^{4} - 2\,a
^{6}\,r^{6} + 2\,a^{6}\,x^{6} + a^{8}\,r^{4} \\
\mytaba + 2\,y^{6}\,b^{2}\,a^{2}\,x^{2} - 2\,r^{6}\,a^{4}\,y^{2}
 + 6\,r^{6}\,b^{2}\,y^{2}\,a^{2} - 2\,a^{8}\,x^{2}\,r^{2} + y^{4}
\,a^{4}\,r^{4} - 6\,y^{4}\,b^{2}\,r^{4}\,a^{2} \\
\mytaba - 2\,a^{6}\,y^{2}\,r^{4} + 4\,x^{2}\,a^{6}\,y^{2}\,r^{2}
 + 6\,x^{2}\,a^{4}\,y^{2}\,r^{4} + 2\,y^{4}\,b^{2}\,x^{2}\,a^{2}
\,r^{2} - 10\,x^{2}\,b^{2}\,y^{2}\,r^{4}\,a^{2} \\
\mytaba + 2\,y^{6}\,b^{2}\,a^{2}\,r^{2} + 6\,r^{6}\,b^{2}\,x^{2}
\,a^{2} + 4\,b^{2}\,x^{4}\,y^{4}\,a^{2} + 2\,a^{4}\,x^{6}\,y^{2}
 - 6\,a^{6}\,x^{4}\,r^{2} - 4\,r^{6}\,a^{4}\,x^{2} \\
\mytaba + 2\,b^{2}\,x^{6}\,a^{2}\,y^{2} + 6\,r^{4}\,a^{6}\,x^{2}
 - 6\,a^{4}\,x^{4}\,r^{2}\,y^{2} + 2\,b^{6}\,y^{6} - 2\,a^{8}\,b
^{2}\,r^{2} - 2\,x^{2}\,y^{4}\,b^{6} \\
\mytaba + 6\,a^{4}\,y^{4}\,b^{4} + y^{4}\,b^{8} - 6\,y^{4}\,b^{6}
\,r^{2} - 2\,x^{4}\,a^{6}\,y^{2} - 2\,a^{6}\,b^{6} + 2\,b^{6}\,r
^{4}\,a^{2} - 2\,b^{8}\,r^{2}\,a^{2} \\
\mytaba + 2\,b^{4}\,a^{6}\,r^{2} - 6\,b^{4}\,a^{4}\,r^{4} + 2\,r
^{2}\,a^{4}\,b^{6} + 2\,r^{4}\,a^{6}\,b^{2} + 2\,r^{6}\,a^{4}\,b
^{2}) }
}
\end{maplelatex}

\end{smallmaplegroup}
\begin{smallmaplegroup}
\begin{mapleinput}
\mapleinline{active}{1d}{g2:=factor(reducedG[2]);}{%
}
\end{mapleinput}

\mapleresult
\begin{maplelatex}
\mapleinline{inert}{2d}{g2 :=
b^2*(-2*y*a^2*b^2+b^2*y0*y^2+y0*a^2*b^2-b^2*y0*r^2+3*a^2*y*x*x0-a^2*y
0*x^2+a^2*y^2*y0-a^4*y0+a^2*y0*r^2+a^4*y-a^2*y^3-a^2*y*x^2+a^2*y*r^2);
}{%
\maplemultiline{
\mathit{g2} :=  b^{2}( - 2\,y\,a^{2}\,b^{2} + b^{2}\,\mathit{y0}\,y^{2} + \mathit{y0}\,a^{2}\,b^{2} - b^{2}\,\mathit{y0}\,r^{2}
 + 3\,a^{2}\,y\,x\,\mathit{x0} - a^{2}\,\mathit{y0}\,x^{2} + a^{2}\,y^{2}\,\mathit{y0} \\
\mytaba - a^{4}\,\mathit{y0} + a^{2}\,\mathit{y0}\,r^{2} + a^{4}\,y - a^{2}\,y^{3} - a^{2}\,y\,x^{2} + a^{2}\,y\,r^{2}) }
}
\end{maplelatex}

\end{smallmaplegroup}
\begin{smallmaplegroup}
\begin{mapleinput}
\mapleinline{active}{1d}{g3:=factor(reducedG[3]);}{%
}
\end{mapleinput}

\mapleresult
\begin{maplelatex}
\mapleinline{inert}{2d}{g3 :=
a^2*(-x0*b^2*y^2+x0*a^2*x^2+x0*a^2*b^2-x0*a^2*r^2+3*y*b^2*y0*x-2*x*a
^2*b^2+b^2*x^2*x0-b^4*x0+x0*b^2*r^2+b^4*x-b^2*x*y^2-b^2*x^3+x*b^2*r^2)
;}{%
\maplemultiline{
\mathit{g3} := a^{2}( - \mathit{x0}\,b^{2}\,y^{2} + \mathit{x0}\,a^{2}\,x^{2} + \mathit{x0}\,a^{2}\,b^{2} - \mathit{x0}\,a^{2}
\,r^{2} + 3\,y\,b^{2}\,\mathit{y0}\,x - 2\,x\,a^{2}\,b^{2} + b^{2}\,x^{2}\,\mathit{x0} \\
\mytaba - b^{4}\,\mathit{x0} + \mathit{x0}\,b^{2}\,r^{2} + b^{4}\,x - b^{2}\,x\,y^{2} - b^{2}\,x^{3} + x\,b^{2}\,r^{2}) }
}
\end{maplelatex}

\end{smallmaplegroup}
\begin{smallmaplegroup}
\begin{mapleinput}
\mapleinline{active}{1d}{g4:=factor(reducedG[4]);}{%
}
\end{mapleinput}

\mapleresult
\begin{maplelatex}
\mapleinline{inert}{2d}{g4 :=
b^2*(3*x0*a^2*x^3+6*a^4*x*x0-3*b^2*x*a^2*x0-3*x*x0*a^2*r^2+4*a^2*y*y
0*x^2+2*a^2*y^3*y0-2*a^4*y*y0-b^2*y0*y^3+5*y*y0*a^2*b^2-4*a^2*y*y0*r^2
+y*b^2*y0*r^2-3*x^4*a^2-5*x^2*y^2*a^2-2*y^4*a^2-3*a^4*x^2+2*a^4*y^2+3*
b^2*x^2*a^2-y^2*a^2*b^2-3*a^4*b^2+6*x^2*a^2*r^2+5*y^2*a^2*r^2+3*a^4*r^
2+3*b^2*a^2*r^2-3*r^4*a^2);}{%
\maplemultiline{
\mathit{g4} := b^{2}(3\,\mathit{x0}\,a^{2}\,x^{3} + 6\,a^{4}\,x\,\mathit{x0} - 3\,b^{2}\,x\,a^{2}\,\mathit{x0} - 3\,x\,\mathit{
x0}\,a^{2}\,r^{2} + 4\,a^{2}\,y\,\mathit{y0}\,x^{2} + 2\,a^{2}\,y^{3}\,\mathit{y0} \\
\mytaba - 2\,a^{4}\,y\,\mathit{y0} - b^{2}\,\mathit{y0}\,y^{3} + 5\,y\,\mathit{y0}\,a^{2}\,b^{2} - 4\,a^{2}\,y\,\mathit{y0}\,r^{2}
 + y\,b^{2}\,\mathit{y0}\,r^{2} - 3\,x^{4}\,a^{2} - 5\,x^{2}\,y^{2}\,a^{2} \\
\mytaba - 2\,y^{4}\,a^{2} - 3\,a^{4}\,x^{2} + 2\,a^{4}\,y^{2} + 3\,b^{2}\,x^{2}\,a^{2} - y^{2}\,a^{2}\,b^{2} - 3\,a^{4}\,b^{2} + 6
\,x^{2}\,a^{2}\,r^{2} + 5\,y^{2}\,a^{2}\,r^{2} \\
\mytaba + 3\,a^{4}\,r^{2} + 3\,b^{2}\,a^{2}\,r^{2} - 3\,r^{4}\,a^{2}) }
}
\end{maplelatex}

\end{smallmaplegroup}
\begin{smallmaplegroup}
\begin{mapleinput}
\mapleinline{active}{1d}{g5:=factor(reducedG[5]);}{%
}
\end{mapleinput}

\mapleresult
\begin{maplelatex}
\mapleinline{inert}{2d}{g5 :=
b^2*(-b^4*x*y0*r^2+a^2*b^4*y0*x+3*x0*b^4*y*a^2-5*b^4*x*y*a^2+b^4*y^2*
y0*x+2*b^2*x^3*a^2*y-7*b^2*x*a^2*y^2*y0-3*x0*b^2*y*a^2*r^2-a^2*b^2*y0*
x^3-3*x0*a^4*y*b^2-2*b^2*x*a^2*y*r^2+2*b^2*x*a^2*y^3+5*x*a^4*y*b^2+3*x
0*b^2*y^3*a^2+a^4*x*y^2*y0-x^3*y*a^4+a^4*x*y0*r^2+3*a^4*x0*y*r^2-a^6*x
*y0-a^4*x*y^3-a^4*y0*x^3+y*a^4*x*r^2+a^6*x*y);}{%
\maplemultiline{
\mathit{g5} :=  b^{2}( - b^{4}\,x\,\mathit{y0}\,r^{2} + a^{2}\,b^{4}\,\mathit{y0}\,x + 3\,\mathit{x0}\,b^{4}\,y\,a^{2} - 5\,b^{4
}\,x\,y\,a^{2} + b^{4}\,y^{2}\,\mathit{y0}\,x + 2\,b^{2}\,x^{3}\,a^{2}\,y \\
\mytaba - 7\,b^{2}\,x\,a^{2}\,y^{2}\,\mathit{y0} - 3\,\mathit{x0}\,b^{2}\,y\,a^{2}\,r^{2} - a^{2}\,b^{2}\,\mathit{y0}\,x^{3} - 3\,
\mathit{x0}\,a^{4}\,y\,b^{2} - 2\,b^{2}\,x\,a^{2}\,y\,r^{2} \\
\mytaba + 2\,b^{2}\,x\,a^{2}\,y^{3} + 5\,x\,a^{4}\,y\,b^{2} + 3\,\mathit{x0}\,b^{2}\,y^{3}\,a^{2} + a^{4}\,x\,y^{2}\,\mathit{y0}
 - x^{3}\,y\,a^{4} + a^{4}\,x\,\mathit{y0}\,r^{2} \\
\mytaba + 3\,a^{4}\,\mathit{x0}\,y\,r^{2} - a^{6}\,x\,\mathit{y0} - a^{4}\,x\,y^{3} - a^{4}\,\mathit{y0}\,x^{3} + y\,a^{4}\,x\,r^{2} + a^{6}\,x\,y) }
}
\end{maplelatex}

\end{smallmaplegroup}
\begin{smallmaplegroup}
\begin{mapleinput}
\mapleinline{active}{1d}{g6:=factor(reducedG[6]);}{%
}
\end{mapleinput}

\mapleresult
\begin{maplelatex}
\mapleinline{inert}{2d}{g6 :=
b^2*(-2*a^4*y^4-6*b^2*r^2*a^2*x*x0+b^4*y^2*a^2+4*a^4*y*y0*b^2+4*y0*b
^4*y*a^2-b^2*x^2*a^2*y^2+3*b^2*x^2*a^2*r^2+4*b^2*y^2*a^2*r^2-4*a^2*y*b
^2*y0*x^2-4*a^2*y*b^2*y0*r^2-2*b^4*y0*y^3+3*b^4*r^2*a^2-3*b^2*r^4*a^2-
3*a^4*b^4+2*y*b^4*y0*r^2-b^2*y^4*a^2+6*a^4*b^2*x^2+6*a^4*b^2*r^2-3*a^6
*b^2+4*a^4*y*y0*x^2-4*a^4*y*y0*r^2+3*a^6*r^2-3*a^4*r^4+2*a^4*y^3*y0-2*
a^6*y*y0+6*a^4*x^2*r^2+5*a^4*y^2*r^2+6*a^6*x*x0-3*x^4*a^4-3*x^2*a^6-5*
a^4*x^2*y^2+2*a^6*y^2);}{%
\maplemultiline{
\mathit{g6} := b^{2}( - 2\,a^{4}\,y^{4} - 6\,b^{2}\,r^{2}\,a^{2}\,x\,\mathit{x0} + b^{4}\,y^{2}\,a^{2} + 4\,a^{4}\,y\,\mathit{
y0}\,b^{2} + 4\,\mathit{y0}\,b^{4}\,y\,a^{2} - b^{2}\,x^{2}\,a^{2}\,y^{2} \\
\mytaba + 3\,b^{2}\,x^{2}\,a^{2}\,r^{2} + 4\,b^{2}\,y^{2}\,a^{2}\,r^{2} - 4\,a^{2}\,y\,b^{2}\,\mathit{y0}\,x^{2} - 4\,a^{2}\,y\,b
^{2}\,\mathit{y0}\,r^{2} - 2\,b^{4}\,\mathit{y0}\,y^{3} \\
\mytaba + 3\,b^{4}\,r^{2}\,a^{2} - 3\,b^{2}\,r^{4}\,a^{2} - 3\,a^{4}\,b^{4} + 2\,y\,b^{4}\,\mathit{y0}\,r^{2} - b^{2}\,y^{4}\,a^{
2} + 6\,a^{4}\,b^{2}\,x^{2} + 6\,a^{4}\,b^{2}\,r^{2} \\
\mytaba - 3\,a^{6}\,b^{2} + 4\,a^{4}\,y\,\mathit{y0}\,x^{2} - 4\,a^{4}\,y\,\mathit{y0}\,r^{2} + 3\,a^{6}\,r^{2} - 3\,a^{4}\,r^{4}
 + 2\,a^{4}\,y^{3}\,\mathit{y0} - 2\,a^{6}\,y\,\mathit{y0} \\
\mytaba + 6\,a^{4}\,x^{2}\,r^{2} + 5\,a^{4}\,y^{2}\,r^{2} + 6\,a^{6}\,x\,\mathit{x0} - 3\,x^{4}\,a^{4} - 3\,x^{2}\,a^{6} - 5\,a^{
4}\,x^{2}\,y^{2} + 2\,a^{6}\,y^{2}) }
}
\end{maplelatex}

\end{smallmaplegroup}
\begin{smallmaplegroup}
\begin{mapleinput}
\mapleinline{active}{1d}{g7:=factor(reducedG[7]);}{%
}
\end{mapleinput}

\mapleresult
\begin{maplelatex}
\mapleinline{inert}{2d}{g7 :=
b^2*(-21*a^4*x*b^4-6*x^3*b^4*a^2+2*a^6*y^2*x-3*x^5*a^4+6*b^6*x*a^2+9
*a^6*b^2*x+6*a^4*y^2*b^2*x-5*b^4*y^2*x*a^2-5*y^2*a^4*x^3-2*y^4*a^4*x+1
2*a^4*b^2*x^3-3*a^6*x^3-6*a^6*b^2*x0-6*b^6*x0*a^2+4*y*a^4*y0*x^3+2*y^3
*a^4*x*y0-2*y*a^6*x*y0-14*y*a^4*b^2*y0*x-4*y*a^4*x*y0*r^2+5*y^2*a^4*x*
r^2-4*y0*b^2*x^3*a^2*y+22*y0*b^4*x*y*a^2+6*a^4*x^3*r^2+3*a^6*x*r^2-2*b
^4*x*y0*y^3-3*a^4*x*r^4-b^2*x^3*a^2*y^2-b^2*y^4*x*a^2-6*b^4*y^2*a^2*x0
-12*b^2*a^4*x0*r^2+6*b^4*x^2*a^2*x0+6*x0*b^4*r^2*a^2+6*x0*a^6*r^2+6*x0
*a^4*b^2*y^2+2*b^4*x*y*y0*r^2+12*x0*a^4*b^4+3*b^2*x^3*a^2*r^2-4*y*a^2*
b^2*x*y0*r^2-3*b^2*x*r^4*a^2-6*b^2*x^2*r^2*a^2*x0+4*b^2*y^2*x*a^2*r^2+
9*b^4*r^2*a^2*x);}{%
\maplemultiline{
\mathit{g7} := b^{2}( - 21\,a^{4}\,x\,b^{4} - 6\,x^{3}\,b^{4}\,a^{2} + 2\,a^{6}\,y^{2}\,x - 3\,x^{5}\,a^{4} + 6\,b^{6}\,x\,a^{
2} + 9\,a^{6}\,b^{2}\,x + 6\,a^{4}\,y^{2}\,b^{2}\,x \\
\mytaba - 5\,b^{4}\,y^{2}\,x\,a^{2} - 5\,y^{2}\,a^{4}\,x^{3} - 2\,y^{4}\,a^{4}\,x + 12\,a^{4}\,b^{2}\,x^{3} - 3\,a^{6}\,x^{3} - 6
\,a^{6}\,b^{2}\,\mathit{x0} \\
\mytaba - 6\,b^{6}\,\mathit{x0}\,a^{2} + 4\,y\,a^{4}\,\mathit{y0}\,x^{3} + 2\,y^{3}\,a^{4}\,x\,\mathit{y0} - 2\,y\,a^{6}\,x\,
\mathit{y0} - 14\,y\,a^{4}\,b^{2}\,\mathit{y0}\,x \\
\mytaba - 4\,y\,a^{4}\,x\,\mathit{y0}\,r^{2} + 5\,y^{2}\,a^{4}\,x\,r^{2} - 4\,\mathit{y0}\,b^{2}\,x^{3}\,a^{2}\,y + 22\,\mathit{y0
}\,b^{4}\,x\,y\,a^{2} + 6\,a^{4}\,x^{3}\,r^{2} \\
\mytaba + 3\,a^{6}\,x\,r^{2} - 2\,b^{4}\,x\,\mathit{y0}\,y^{3} - 3\,a^{4}\,x\,r^{4} - b^{2}\,x^{3}\,a^{2}\,y^{2} - b^{2}\,y^{4}\,x
\,a^{2} - 6\,b^{4}\,y^{2}\,a^{2}\,\mathit{x0} \\
\mytaba - 12\,b^{2}\,a^{4}\,\mathit{x0}\,r^{2} + 6\,b^{4}\,x^{2}\,a^{2}\,\mathit{x0} + 6\,\mathit{x0}\,b^{4}\,r^{2}\,a^{2} + 6\,
\mathit{x0}\,a^{6}\,r^{2} + 6\,\mathit{x0}\,a^{4}\,b^{2}\,y^{2}\\
\mytaba + 2\,b^{4}\,x\,y\,\mathit{y0}\,r^{2} + 12\,\mathit{x0}\,a^{4}\,b^{4} + 3\,b^{2}\,x^{3}\,a^{2}\,r^{2} - 4\,y\,a^{2}\,b^{2}
\,x\,\mathit{y0}\,r^{2} - 3\,b^{2}\,x\,r^{4}\,a^{2} \\
\mytaba - 6\,b^{2}\,x^{2}\,r^{2}\,a^{2}\,\mathit{x0} + 4\,b^{2}\,y^{2}\,x\,a^{2}\,r^{2} + 9\,b^{4}\,r^{2}\,a^{2}\,x) }
}
\end{maplelatex}

\end{smallmaplegroup}
\begin{smallmaplegroup}
\begin{mapleinput}
\mapleinline{active}{1d}{g8:=factor(reducedG[8]);}{%
}
\end{mapleinput}

\mapleresult
\begin{maplelatex}
\mapleinline{inert}{2d}{g8 :=
b^2*(4*a^4*y^3-2*a^4*x^2*y-4*a^2*y^3*x^2-2*y^5*a^2-2*a^6*y-2*a^2*b^2*y
0*x^2+4*a^2*b^2*y0*y^2+4*a^2*b^2*y0*r^2+2*a^4*y*b^2-2*b^2*y^3*a^2-2*b^
4*y*a^2-b^4*y0*r^2+b^4*y0*a^2+b^4*y0*y^2-3*a^4*y0*b^2+2*a^4*y*r^2+2*b^
2*y*a^2*r^2+4*b^2*y*a^2*x^2+2*a^6*y0+4*x^2*a^2*y*r^2+3*x^2*a^2*y^2*y0-
2*x^2*a^2*y0*r^2-3*a^2*y^2*y0*r^2-b^2*y^2*y0*x^2+2*b^2*y^2*y0*r^2+b^2*
x^2*y0*r^2+4*y^3*a^2*r^2-2*y*r^4*a^2+r^4*a^2*y0-3*a^4*y0*r^2-b^2*y^4*y
0-b^2*r^4*y0+2*a^2*y^4*y0+a^2*y0*x^4+3*x^2*a^4*y0-4*a^4*y^2*y0-2*x^4*y
*a^2);}{%
\maplemultiline{
\mathit{g8} := b^{2}(4\,a^{4}\,y^{3} - 2\,a^{4}\,x^{2}\,y - 4\,a^{2}\,y^{3}\,x^{2} - 2\,y^{5}\,a^{2} - 2\,a^{6}\,y - 2\,a^{2}\,b
^{2}\,\mathit{y0}\,x^{2} + 4\,a^{2}\,b^{2}\,\mathit{y0}\,y^{2}\\
\mytaba + 4\,a^{2}\,b^{2}\,\mathit{y0}\,r^{2} + 2\,a^{4}\,y\,b^{2} - 2\,b^{2}\,y^{3}\,a^{2} - 2\,b^{4}\,y\,a^{2} - b^{4}\,\mathit{
y0}\,r^{2} + b^{4}\,\mathit{y0}\,a^{2} + b^{4}\,\mathit{y0}\,y^{2} \\
\mytaba - 3\,a^{4}\,\mathit{y0}\,b^{2} + 2\,a^{4}\,y\,r^{2} + 2\,b^{2}\,y\,a^{2}\,r^{2} + 4\,b^{2}\,y\,a^{2}\,x^{2} + 2\,a^{6}\,
\mathit{y0} + 4\,x^{2}\,a^{2}\,y\,r^{2} \\
\mytaba + 3\,x^{2}\,a^{2}\,y^{2}\,\mathit{y0} - 2\,x^{2}\,a^{2}\,\mathit{y0}\,r^{2} - 3\,a^{2}\,y^{2}\,\mathit{y0}\,r^{2} - b^{2}
\,y^{2}\,\mathit{y0}\,x^{2} + 2\,b^{2}\,y^{2}\,\mathit{y0}\,r^{2}\\
\mytaba + b^{2}\,x^{2}\,\mathit{y0}\,r^{2} + 4\,y^{3}\,a^{2}\,r^{2} - 2\,y\,r^{4}\,a^{2} + r^{4}\,a^{2}\,\mathit{y0} - 3\,a^{4}\,
\mathit{y0}\,r^{2} - b^{2}\,y^{4}\,\mathit{y0} - b^{2}\,r^{4}\,\mathit{y0} \\
\mytaba + 2\,a^{2}\,y^{4}\,\mathit{y0} + a^{2}\,\mathit{y0}\,x^{4} + 3\,x^{2}\,a^{4}\,\mathit{y0} - 4\,a^{4}\,y^{2}\,\mathit{y0}
 - 2\,x^{4}\,y\,a^{2}) }
}
\end{maplelatex}

\end{smallmaplegroup}
\begin{smallmaplegroup}
\begin{mapleinput}
\mapleinline{active}{1d}{g9:=factor(reducedG[9]);}{%
}
\end{mapleinput}

\mapleresult
\begin{maplelatex}
\mapleinline{inert}{2d}{g9 :=
b^2*(-4*b^2*x^2*a^2*y^2*y0-2*x^2*a^2*b^2*y0*r^2-2*b^4*y^4*y0+b^4*y^3
*a^2-b^2*y^5*a^2-3*y*a^4*b^4+4*b^4*y^2*y0*a^2+4*b^4*y^2*y0*r^2-b^2*x^2
*a^2*y^3+6*x^2*a^4*y*b^2-2*a^2*y^2*b^2*y0*r^2+2*a^4*y^2*y0*b^2+2*b^2*y
^3*a^2*r^2-y*b^4*r^2*a^2-y*b^2*r^4*a^2+8*y*a^4*b^2*r^2+b^2*x^2*a^2*y*r
^2+2*r^4*a^2*b^2*y0+2*a^2*b^4*y0*r^2-2*a^6*y0*b^2+a^6*y*b^2-2*r^4*b^4*
y0-3*a^4*x^4*y+4*a^6*y^3-2*a^4*y^5-5*x^2*a^4*y^3-2*a^6*r^2*y0+a^6*r^2*
y+5*a^4*r^2*y^3-3*a^4*r^4*y+2*a^4*y^4*y0-4*a^6*y^2*y0+2*a^6*x^2*y0+2*a
^8*y0+6*a^4*x^2*y*r^2+4*a^4*x^2*y^2*y0-4*a^4*r^2*y0*y^2-a^6*x^2*y-2*a^
8*y);}{%
\maplemultiline{
\mathit{g9} := b^{2}( - 4\,b^{2}\,x^{2}\,a^{2}\,y^{2}\,\mathit{y0} - 2\,x^{2}\,a^{2}\,b^{2}\,\mathit{y0}\,r^{2} - 2\,b^{
4}\,y^{4}\,\mathit{y0} + b^{4}\,y^{3}\,a^{2} - b^{2}\,y^{5}\,a^{2} - 3\,y\,a^{4}\,b^{4} \\
\mytaba + 4\,b^{4}\,y^{2}\,\mathit{y0}\,a^{2} + 4\,b^{4}\,y^{2}\,\mathit{y0}\,r^{2} - b^{2}\,x^{2}\,a^{2}\,y^{3} + 6\,x^{2}\,a^{4}
\,y\,b^{2} - 2\,a^{2}\,y^{2}\,b^{2}\,\mathit{y0}\,r^{2} \\
\mytaba + 2\,a^{4}\,y^{2}\,\mathit{y0}\,b^{2} + 2\,b^{2}\,y^{3}\,a^{2}\,r^{2} - y\,b^{4}\,r^{2}\,a^{2} - y\,b^{2}\,r^{4}\,a^{2} + 
8\,y\,a^{4}\,b^{2}\,r^{2} + b^{2}\,x^{2}\,a^{2}\,y\,r^{2} \\
\mytaba + 2\,r^{4}\,a^{2}\,b^{2}\,\mathit{y0} + 2\,a^{2}\,b^{4}\,\mathit{y0}\,r^{2} - 2\,a^{6}\,\mathit{y0}\,b^{2} + a^{6}\,y\,b^{
2} - 2\,r^{4}\,b^{4}\,\mathit{y0} - 3\,a^{4}\,x^{4}\,y \\
\mytaba + 4\,a^{6}\,y^{3} - 2\,a^{4}\,y^{5} - 5\,x^{2}\,a^{4}\,y^{3} - 2\,a^{6}\,r^{2}\,\mathit{y0} + a^{6}\,r^{2}\,y + 5\,a^{4}
\,r^{2}\,y^{3} - 3\,a^{4}\,r^{4}\,y \\
\mytaba + 2\,a^{4}\,y^{4}\,\mathit{y0} - 4\,a^{6}\,y^{2}\,\mathit{y0} + 2\,a^{6}\,x^{2}\,\mathit{y0} + 2\,a^{8}\,\mathit{y0
} + 6\,a^{4}\,x^{2}\,y\,r^{2} + 4\,a^{4}\,x^{2}\,y^{2}\,\mathit{y0} \\
\mytaba - 4\,a^{4}\,r^{2}\,\mathit{y0}\,y^{2} - a^{6}\,x^{2}\,y - 2\,a^{8}\,y) }
}
\end{maplelatex}

\end{smallmaplegroup}
\begin{smallmaplegroup}
\begin{mapleinput}
\mapleinline{active}{1d}{g10:=factor(reducedG[10]);}{%
}
\end{mapleinput}

\mapleresult
\begin{maplelatex}
\mapleinline{inert}{2d}{g10 :=
b^2*(2*a^6*b^4-2*a^4*b^6-16*a^4*y^4*b^2+2*a^4*x^6-4*b^4*r^4*a^2+2*b
^6*r^2*a^2-4*b^6*y0*y^3+3*a^4*x^4*y^2-2*a^8*y^2+4*a^6*y^4-4*a^6*r^2*b^
2+4*a^4*b^2*x^2*y^2+9*a^6*y^2*b^2+2*b^4*a^4*r^2+6*b^4*x^2*a^4-6*y^2*a^
4*x^2*r^2+2*r^6*b^2*a^2-4*x^2*a^6*r^2+6*x^2*a^4*r^4+6*b^2*x^4*y^2*a^2-
6*a^4*x^4*r^2+3*y^2*a^2*r^4*b^2+13*y^4*b^2*a^2*x^2-12*y^4*b^2*a^2*r^2-
2*r^6*a^4+2*r^4*a^6+a^4*y^4*r^2+3*a^4*y^2*r^4+5*a^6*y^2*r^2+2*a^6*x^4+
2*b^2*r^2*a^2*x^4-4*x^2*b^2*r^4*a^2-9*y^2*b^2*r^2*a^2*x^2+7*y^6*a^2*b^
2-6*x^4*a^4*b^2-4*b^4*x^2*r^2*a^2+2*y^5*a^4*y0-4*y^3*a^6*y0+2*y*a^8*y0
+6*b^4*y^5*y0+6*b^2*r^2*a^2*x^2*y*y0-4*b^2*x^2*a^2*y^3*y0+14*b^2*y^3*a
^2*y0*r^2-6*b^2*a^2*y*r^4*y0-4*b^4*x^2*y*y0*r^2+4*b^4*x^2*y^3*y0-8*b^2
*y^5*a^2*y0+6*a^6*x^2*y*y0-6*a^6*y*r^2*y0-12*b^4*y^3*y0*r^2+6*b^4*y*r^
4*y0+2*a^4*b^2*r^4-2*y^6*a^4+18*a^4*y^3*b^2*y0-10*a^6*y*y0*b^2+8*a^4*y
0*b^4*y+b^4*y^2*a^2*r^2+4*b^6*y*y0*r^2-8*a^4*y*b^2*y0*x^2+28*a^4*y*b^2
*y0*r^2-26*a^2*y*b^4*y0*r^2+4*a^4*b^2*x^2*r^2-18*a^4*b^2*y^2*r^2-12*a^
2*b^4*y0*y^3-x^2*y^4*a^4-4*x^2*a^6*b^2+3*y^4*b^4*a^2-7*a^4*b^4*y^2+6*y
^2*b^6*a^2-12*b^4*x^2*y^2*a^2-5*a^6*y^2*x^2);}{%
\maplemultiline{
\mathit{g10} :=  b^{2}(2\,a^{6}\,b^{4} - 2\,a^{4}\,b^{6} - 16\,a^{4}\,y^{4}\,b^{2} + 2\,a^{4}\,x^{6} - 4\,b^{4}\,r^{4}\,a^{2
} + 2\,b^{6}\,r^{2}\,a^{2} - 4\,b^{6}\,\mathit{y0}\,y^{3} \\
\mytaba + 3\,a^{4}\,x^{4}\,y^{2} - 2\,a^{8}\,y^{2} + 4\,a^{6}\,y^{4} - 4\,a^{6}\,r^{2}\,b^{2} + 4\,a^{4}\,b^{2}\,x^{2}\,y^{2} + 9
\,a^{6}\,y^{2}\,b^{2} + 2\,b^{4}\,a^{4}\,r^{2} \\
\mytaba + 6\,b^{4}\,x^{2}\,a^{4} - 6\,y^{2}\,a^{4}\,x^{2}\,r^{2} + 2\,r^{6}\,b^{2}\,a^{2} - 4\,x^{2}\,a^{6}\,r^{2} + 6\,x^{2}\,a
^{4}\,r^{4} + 6\,b^{2}\,x^{4}\,y^{2}\,a^{2} \\
\mytaba - 6\,a^{4}\,x^{4}\,r^{2} + 3\,y^{2}\,a^{2}\,r^{4}\,b^{2} + 13\,y^{4}\,b^{2}\,a^{2}\,x^{2} - 12\,y^{4}\,b^{2}\,a^{2}\,r^{2
} - 2\,r^{6}\,a^{4} + 2\,r^{4}\,a^{6} \\
\mytaba + a^{4}\,y^{4}\,r^{2} + 3\,a^{4}\,y^{2}\,r^{4} + 5\,a^{6}\,y^{2}\,r^{2} + 2\,a^{6}\,x^{4} + 2\,b^{2}\,r^{2}\,a^{2}\,x^{4}
 - 4\,x^{2}\,b^{2}\,r^{4}\,a^{2} \\
\mytaba - 9\,y^{2}\,b^{2}\,r^{2}\,a^{2}\,x^{2} + 7\,y^{6}\,a^{2}\,b^{2} - 6\,x^{4}\,a^{4}\,b^{2} - 4\,b^{4}\,x^{2}\,r^{2}\,a^{2}
 + 2\,y^{5}\,a^{4}\,\mathit{y0} - 4\,y^{3}\,a^{6}\,\mathit{y0}\\
\mytaba + 2\,y\,a^{8}\,\mathit{y0} + 6\,b^{4}\,y^{5}\,\mathit{y0} + 6\,b^{2}\,r^{2}\,a^{2}\,x^{2}\,y\,\mathit{y0} - 4\,b^{2}\,x^{2
}\,a^{2}\,y^{3}\,\mathit{y0} + 14\,b^{2}\,y^{3}\,a^{2}\,\mathit{y0}\,r^{2} \\
\mytaba - 6\,b^{2}\,a^{2}\,y\,r^{4}\,\mathit{y0} - 4\,b^{4}\,x^{2}\,y\,\mathit{y0}\,r^{2} + 4\,b^{4}\,x^{2}\,y^{3}\,\mathit{y0} - 
8\,b^{2}\,y^{5}\,a^{2}\,\mathit{y0} + 6\,a^{6}\,x^{2}\,y\,\mathit{y0} \\
\mytaba - 6\,a^{6}\,y\,r^{2}\,\mathit{y0} - 12\,b^{4}\,y^{3}\,\mathit{y0}\,r^{2} + 6\,b^{4}\,y\,r^{4}\,\mathit{y0} + 2\,a^{4}\,
b^{2}\,r^{4} - 2\,y^{6}\,a^{4} + 18\,a^{4}\,y^{3}\,b^{2}\,\mathit{y0} \\
\mytaba - 10\,a^{6}\,y\,\mathit{y0}\,b^{2} + 8\,a^{4}\,\mathit{y0}\,b^{4}\,y + b^{4}\,y^{2}\,a^{2}\,r^{2} + 4\,b^{6}\,y\,\mathit{
y0}\,r^{2} - 8\,a^{4}\,y\,b^{2}\,\mathit{y0}\,x^{2} \\
\mytaba + 28\,a^{4}\,y\,b^{2}\,\mathit{y0}\,r^{2} - 26\,a^{2}\,y\,b^{4}\,\mathit{y0}\,r^{2} + 4\,a^{4}\,b^{2}\,x^{2}\,r^{2} - 18
\,a^{4}\,b^{2}\,y^{2}\,r^{2} - 12\,a^{2}\,b^{4}\,\mathit{y0}\,y^{3} \\
\mytaba - x^{2}\,y^{4}\,a^{4} - 4\,x^{2}\,a^{6}\,b^{2} + 3\,y^{4}\,b^{4}\,a^{2} - 7\,a^{4}\,b^{4}\,y^{2} + 6\,y^{2}\,b^{6}\,a^{2}
 - 12\,b^{4}\,x^{2}\,y^{2}\,a^{2} \\
\mytaba - 5\,a^{6}\,y^{2}\,x^{2}) }
}
\end{maplelatex}

\end{smallmaplegroup}
\begin{smallmaplegroup}
\begin{mapleinput}
\mapleinline{active}{1d}{g11:=factor(reducedG[11]);}{%
}
\end{mapleinput}

\mapleresult
\begin{maplelatex}
\mapleinline{inert}{2d}{g11 :=
b^2*(-12*b^4*x^4*a^4+18*y^6*b^4*a^2+14*b^6*x^2*a^4-6*b^4*x^2*a^6-5*a
^4*b^8+15*y^2*b^8*a^2-27*y^2*b^6*x^2*a^2+6*a^6*y^4*b^2-8*a^6*y^2*x^2*b
^2-6*b^4*x^2*a^4*y^2-11*a^4*y^2*b^6+22*a^6*y^2*b^2*r^2+15*a^4*b^2*x^4*
y^2+24*b^6*y^5*y0+2*a^4*x^6*b^2+16*a^6*x^2*b^2*y*y0-24*b^4*x^2*a^4*y*y
0-8*a^6*x^2*b^2*r^2+6*b^4*x^2*a^4*r^2-9*a^4*x^4*b^2*r^2-33*b^4*r^2*a^2
*x^2*y^2+18*b^4*r^2*a^2*x^2*y*y0-29*a^4*x^2*b^2*r^2*y^2-4*a^4*x^2*b^2*
r^2*y*y0+12*a^4*x^2*b^2*r^4+3*b^4*r^2*a^2*x^4-10*b^4*r^4*a^2*x^2-6*b^4
*x^2*a^2*y^3*y0+8*x^2*b^6*y^3*y0-8*x^2*b^6*y*y0*r^2-42*b^6*y^3*y0*a^2-
48*b^6*y^3*y0*r^2-7*y^4*a^4*b^2*x^2+50*b^4*y^3*a^2*y0*r^2+9*b^4*x^4*a^
2*y^2-9*b^6*x^2*r^2*a^2+43*b^4*x^2*y^4*a^2-22*b^4*y^5*a^2*y0+50*b^4*y^
3*a^4*y0-29*b^4*y^4*a^2*r^2+10*b^6*y^2*r^2*a^2+4*b^4*y^2*r^4*a^2-73*b^
4*y^2*a^4*r^2+7*r^6*a^2*b^4+a^8*b^4+86*y*a^4*b^4*y0*r^2-78*y*a^2*b^6*y
0*r^2-10*y^3*b^8*y0-28*y*a^6*b^4*y0+24*y*a^4*b^6*y0+10*y*b^8*y0*r^2-28
*a^2*y*r^4*b^4*y0+4*a^4*y^3*r^2*b^2*y0-18*a^6*y*r^2*y0*b^2+4*a^4*y*r^4
*b^2*y0-9*a^4*y^4*r^2*b^2+14*a^4*y^2*r^4*b^2+24*y*b^6*r^4*y0-2*a^8*b^2
*x^2-6*y^3*a^6*y0*b^2+4*y*a^8*y0*b^2+x^4*a^4*y^4+y^4*b^6*a^2-4*a^8*y^2
*b^2+22*a^6*y^2*b^4-2*y^3*b^2*r^2*a^2*y0*x^2-2*y^4*a^4*x^2*r^2-2*x^4*b
^4*y*y0*r^2-b^2*r^2*a^2*x^4*y^2-r^8*a^2*b^2+r^8*a^4+a^8*x^4-2*y^8*b^2*
a^2+a^4*x^8-4*a^4*x^6*r^2+6*a^4*x^4*r^4+b^2*r^2*a^2*x^6-3*b^2*r^4*a^2*
x^4-2*a^6*r^6+2*a^6*x^6+a^8*r^4-y^6*b^2*a^2*x^2+2*x^2*y*b^4*r^4*y0-2*r
^6*a^4*y^2+5*r^6*b^2*y^2*a^2-2*a^8*x^2*r^2+y^4*a^4*r^4-9*y^4*b^2*r^4*a
^2+2*a^2*y^7*b^2*y0-2*a^6*y^2*r^4-4*y^5*b^2*a^2*y0*r^2+2*y^3*a^2*r^4*b
^2*y0+4*x^2*a^6*y^2*r^2+6*x^2*a^4*y^2*r^4+5*y^4*b^2*x^2*a^2*r^2-7*x^2*
b^2*y^2*r^4*a^2+2*y^5*b^2*a^2*y0*x^2-4*x^2*b^4*y^3*y0*r^2+7*y^6*b^2*a^
2*r^2+3*r^6*b^2*x^2*a^2+2*x^2*b^4*y^5*y0+4*b^2*x^4*y^4*a^2+2*a^4*x^6*y
^2-6*a^6*x^4*r^2-4*r^6*a^4*x^2+2*x^4*b^4*y^3*y0+3*b^2*x^6*a^2*y^2+6*r^
4*a^6*x^2-6*a^4*x^4*r^2*y^2-2*a^8*b^2*r^2-41*a^4*y^4*b^4-2*x^4*a^6*y^2
+4*a^6*b^6-11*b^6*r^4*a^2+5*b^8*r^2*a^2-10*b^4*a^6*r^2+2*b^4*a^4*r^4+7
*r^2*a^4*b^6+8*r^4*a^6*b^2-5*r^6*a^4*b^2);}{%
\maplemultiline{
\mathit{g11} := b^{2}( - 12\,b^{4}\,x^{4}\,a^{4} + 18\,y^{6}\,b^{4}\,a^{2} + 14\,b^{6}\,x^{2}\,a^{4} - 6\,b^{4}\,x^{2}\,a^{6}
 - 5\,a^{4}\,b^{8} + 15\,y^{2}\,b^{8}\,a^{2} \\
\mytaba - 27\,y^{2}\,b^{6}\,x^{2}\,a^{2} + 6\,a^{6}\,y^{4}\,b^{2} - 8\,a^{6}\,y^{2}\,x^{2}\,b^{2} - 6\,b^{4}\,x^{2}\,a^{4}\,y^{2}
 - 11\,a^{4}\,y^{2}\,b^{6} \\
\mytaba + 22\,a^{6}\,y^{2}\,b^{2}\,r^{2} + 15\,a^{4}\,b^{2}\,x^{4}\,y^{2} + 24\,b^{6}\,y^{5}\,\mathit{y0} + 2\,a^{4}\,x^{6}\,b^{2}
 + 16\,a^{6}\,x^{2}\,b^{2}\,y\,\mathit{y0} \\
\mytaba - 24\,b^{4}\,x^{2}\,a^{4}\,y\,\mathit{y0} - 8\,a^{6}\,x^{2}\,b^{2}\,r^{2} + 6\,b^{4}\,x^{2}\,a^{4}\,r^{2} - 9\,a^{4}\,x^{4
}\,b^{2}\,r^{2} - 33\,b^{4}\,r^{2}\,a^{2}\,x^{2}\,y^{2} \\
\mytaba + 18\,b^{4}\,r^{2}\,a^{2}\,x^{2}\,y\,\mathit{y0} - 29\,a^{4}\,x^{2}\,b^{2}\,r^{2}\,y^{2} - 4\,a^{4}\,x^{2}\,b^{2}\,r^{2}
\,y\,\mathit{y0} + 12\,a^{4}\,x^{2}\,b^{2}\,r^{4} \\
\mytaba + 3\,b^{4}\,r^{2}\,a^{2}\,x^{4} - 10\,b^{4}\,r^{4}\,a^{2}\,x^{2} - 6\,b^{4}\,x^{2}\,a^{2}\,y^{3}\,\mathit{y0} + 8\,x^{2}\,
b^{6}\,y^{3}\,\mathit{y0} - 8\,x^{2}\,b^{6}\,y\,\mathit{y0}\,r^{2} \\
\mytaba - 42\,b^{6}\,y^{3}\,\mathit{y0}\,a^{2} - 48\,b^{6}\,y^{3}\,\mathit{y0}\,r^{2} - 7\,y^{4}\,a^{4}\,b^{2}\,x^{2} + 50\,b^{4}
\,y^{3}\,a^{2}\,\mathit{y0}\,r^{2} + 9\,b^{4}\,x^{4}\,a^{2}\,y^{2} \\
\mytaba - 9\,b^{6}\,x^{2}\,r^{2}\,a^{2} + 43\,b^{4}\,x^{2}\,y^{4}\,a^{2} - 22\,b^{4}\,y^{5}\,a^{2}\,\mathit{y0} + 50\,b^{4}\,y^{3}
\,a^{4}\,\mathit{y0} - 29\,b^{4}\,y^{4}\,a^{2}\,r^{2} \\
\mytaba + 10\,b^{6}\,y^{2}\,r^{2}\,a^{2} + 4\,b^{4}\,y^{2}\,r^{4}\,a^{2} - 73\,b^{4}\,y^{2}\,a^{4}\,r^{2} + 7\,r^{6}\,a^{2}\,b^{4}
 + a^{8}\,b^{4} + 86\,y\,a^{4}\,b^{4}\,\mathit{y0}\,r^{2} \\
\mytaba - 78\,y\,a^{2}\,b^{6}\,\mathit{y0}\,r^{2} - 10\,y^{3}\,b^{8}\,\mathit{y0} - 28\,y\,a^{6}\,b^{4}\,\mathit{y0} + 24\,y\,a^{
4}\,b^{6}\,\mathit{y0} + 10\,y\,b^{8}\,\mathit{y0}\,r^{2} \\
\mytaba - 28\,a^{2}\,y\,r^{4}\,b^{4}\,\mathit{y0} + 4\,a^{4}\,y^{3}\,r^{2}\,b^{2}\,\mathit{y0} - 18\,a^{6}\,y\,r^{2}\,\mathit{y0}
\,b^{2} + 4\,a^{4}\,y\,r^{4}\,b^{2}\,\mathit{y0} \\
\mytaba - 9\,a^{4}\,y^{4}\,r^{2}\,b^{2} + 14\,a^{4}\,y^{2}\,r^{4}\,b^{2} + 24\,y\,b^{6}\,r^{4}\,\mathit{y0} - 2\,a^{8}\,b^{2}\,x^{
2} - 6\,y^{3}\,a^{6}\,\mathit{y0}\,b^{2} \\
\mytaba + 4\,y\,a^{8}\,\mathit{y0}\,b^{2} + x^{4}\,a^{4}\,y^{4} + y^{4}\,b^{6}\,a^{2} - 4\,a^{8}\,y^{2}\,b^{2} + 22\,a^{6}\,y^{2
}\,b^{4} - 2\,y^{3}\,b^{2}\,r^{2}\,a^{2}\,\mathit{y0}\,x^{2} \\
\mytaba - 2\,y^{4}\,a^{4}\,x^{2}\,r^{2} - 2\,x^{4}\,b^{4}\,y\,\mathit{y0}\,r^{2} - b^{2}\,r^{2}\,a^{2}\,x^{4}\,y^{2} - r^{8}\,a
^{2}\,b^{2} + r^{8}\,a^{4} + a^{8}\,x^{4} - 2\,y^{8}\,b^{2}\,a^{2} \\
\mytaba + a^{4}\,x^{8} - 4\,a^{4}\,x^{6}\,r^{2} + 6\,a^{4}\,x^{4}\,r^{4} + b^{2}\,r^{2}\,a^{2}\,x^{6} - 3\,b^{2}\,r^{4}\,a^{2}\,x
^{4} - 2\,a^{6}\,r^{6} + 2\,a^{6}\,x^{6} + a^{8}\,r^{4} \\
\mytaba - y^{6}\,b^{2}\,a^{2}\,x^{2} + 2\,x^{2}\,y\,b^{4}\,r^{4}\,\mathit{y0} - 2\,r^{6}\,a^{4}\,y^{2} + 5\,r^{6}\,b^{2}\,y^{2}\,
a^{2} - 2\,a^{8}\,x^{2}\,r^{2} + y^{4}\,a^{4}\,r^{4} \\
\mytaba - 9\,y^{4}\,b^{2}\,r^{4}\,a^{2} + 2\,a^{2}\,y^{7}\,b^{2}\,\mathit{y0} - 2\,a^{6}\,y^{2}\,r^{4} - 4\,y^{5}\,b^{2}\,a^{2}\,
\mathit{y0}\,r^{2} + 2\,y^{3}\,a^{2}\,r^{4}\,b^{2}\,\mathit{y0} \\
\mytaba + 4\,x^{2}\,a^{6}\,y^{2}\,r^{2} + 6\,x^{2}\,a^{4}\,y^{2}\,r^{4} + 5\,y^{4}\,b^{2}\,x^{2}\,a^{2}\,r^{2} - 7\,x^{2}\,b^{2}
\,y^{2}\,r^{4}\,a^{2} + 2\,y^{5}\,b^{2}\,a^{2}\,\mathit{y0}\,x^{2} \\
\mytaba - 4\,x^{2}\,b^{4}\,y^{3}\,\mathit{y0}\,r^{2} + 7\,y^{6}\,b^{2}\,a^{2}\,r^{2} + 3\,r^{6}\,b^{2}\,x^{2}\,a^{2} + 2\,x^{2}\,b
^{4}\,y^{5}\,\mathit{y0} + 4\,b^{2}\,x^{4}\,y^{4}\,a^{2} \\
\mytaba + 2\,a^{4}\,x^{6}\,y^{2} - 6\,a^{6}\,x^{4}\,r^{2} - 4\,r^{6}\,a^{4}\,x^{2} + 2\,x^{4}\,b^{4}\,y^{3}\,\mathit{y0} + 3\,b^{
2}\,x^{6}\,a^{2}\,y^{2} + 6\,r^{4}\,a^{6}\,x^{2} \\
\mytaba - 6\,a^{4}\,x^{4}\,r^{2}\,y^{2} - 2\,a^{8}\,b^{2}\,r^{2} - 41\,a^{4}\,y^{4}\,b^{4} - 2\,x^{4}\,a^{6}\,y^{2} + 4\,a^{6}\,b
^{6} - 11\,b^{6}\,r^{4}\,a^{2} \\
\mytaba + 5\,b^{8}\,r^{2}\,a^{2} - 10\,b^{4}\,a^{6}\,r^{2} + 2\,b^{4}\,a^{4}\,r^{4} + 7\,r^{2}\,a^{4}\,b^{6} + 8\,r^{4}\,a^{6}\,b
^{2} - 5\,r^{6}\,a^{4}\,b^{2}) }
}
\end{maplelatex}

\end{smallmaplegroup}
\begin{smallmaplegroup}
\begin{mapleinput}
\mapleinline{active}{1d}{g12:=factor(reducedG[12]);}{%
}
\end{mapleinput}

\mapleresult
\begin{maplelatex}
\mapleinline{inert}{2d}{g12 :=
b^2*(-3*y^6*b^4*a^2*x^2-9*a^4*r^6*b^4+2*a^6*b^8+2*a^8*r^4*b^2-19*a^6
*r^2*b^6+32*a^6*b^4*r^4-13*a^8*b^4*r^2+20*a^4*b^8*r^2-10*a^4*b^6*r^4+8
*a^8*b^6+24*y^2*b^10*a^2-15*a^6*r^6*b^2+10*a^6*x^6*b^2-17*a^6*x^4*y^2*
b^2-39*y^2*b^8*x^2*a^2+23*a^6*y^2*x^2*b^4-16*a^8*x^2*y^2*b^2-37*a^4*y^
4*x^2*b^4+24*x^2*y^4*a^6*b^2-16*b^10*y^3*y0+2*b^4*x^6*a^4+94*y^4*b^6*x
^2*a^2+9*y^2*b^6*x^4*a^2-52*b^6*x^2*a^4*y^2+48*b^4*x^4*a^4*y^2-8*b^10*
a^4-18*a^6*b^4*x^4+8*b^10*r^2*a^2-22*b^8*r^4*a^2+21*b^6*r^6*a^2+16*b^1
0*y*y0*r^2-8*r^8*a^2*b^4+8*r^8*a^4*b^2-16*b^6*x^4*a^4+21*b^8*x^2*a^4+2
*y^3*b^4*r^2*a^2*y0*x^2-20*y^4*a^4*b^2*x^2*r^2-2*a^10*b^4-2*a^10*x^4+4
*a^10*b^2*x^2+a^4*x^10-3*a^8*x^6-19*a^8*b^4*x^2+14*a^8*x^4*b^2-16*a^4*
b^2*y^6*x^2-2*a^6*x^6*y^2+5*y^4*a^4*x^6+4*a^4*x^8*y^2-6*x^4*a^6*y^4+6*
x^4*a^8*y^2+2*y^6*a^4*x^4+6*b^6*x^2*a^6-14*y^4*b^8*a^2+24*b^6*y^6*a^2+
20*a^6*y^2*b^6-4*x^4*b^6*y*y0*r^2-8*b^4*r^2*a^2*x^4*y^2-8*a^4*x^6*b^2*
r^2+24*a^4*x^4*b^2*r^4+b^4*r^2*a^2*x^6-10*b^4*r^4*a^2*x^4+4*x^4*b^6*y^
3*y0+3*b^4*x^6*a^2*y^2+4*b^4*a^4*x^2*r^4+40*r^4*a^6*x^2*b^2-10*r^2*a^4
*x^2*b^6+17*r^6*b^4*x^2*a^2-24*r^6*a^4*x^2*b^2-20*a^4*x^4*b^2*r^2*y^2-
30*x^2*b^4*y^2*r^4*a^2-2*x^2*a^2*y*r^4*b^4*y0+40*x^2*a^4*y^2*r^4*b^2-1
2*b^6*x^2*r^4*a^2+8*x^2*b^6*y^5*y0+6*b^4*x^4*y^4*a^2+3*b^6*x^4*a^2*r^2
-35*a^6*x^4*b^2*r^2+16*y^4*b^4*x^2*a^2*r^2+116*y^3*a^2*b^6*y0*r^2-4*y^
5*b^4*a^2*y0*r^2+2*y^3*a^2*r^4*b^4*y0+34*x^2*a^6*y^2*b^2*r^2-20*x^2*b^
6*y^3*y0*r^2+10*b^8*x^2*y0*y^3-13*b^8*x^2*r^2*a^2-50*y^5*b^6*y0*a^2-10
*y^5*b^6*y0*r^2+39*y^6*b^4*a^2*r^2-8*y^4*a^4*b^2*x^4-67*y^4*a^4*b^4*r^
2-31*y^4*b^6*a^2*r^2-8*y^3*a^6*b^4*y0+8*y^3*b^6*r^4*y0+24*y^4*a^6*r^2*
b^2-16*y^6*a^4*r^2*b^2+28*y^4*a^4*r^4*b^2-56*y^4*b^4*r^4*a^2+3*b^4*a^4
*x^4*r^2+6*b^4*a^6*x^2*r^2+4*a^8*y^2*b^4-10*y^8*b^4*a^2+100*a^4*y^3*b^
6*y0-17*a^6*y^2*b^2*r^4-16*a^8*y^2*b^2*r^2-178*a^4*b^6*y^2*r^2+48*a^4*
b^4*r^4*y^2-94*b^8*y^3*y0*a^2-114*b^8*y^3*y0*r^2+32*b^8*y^2*a^2*r^2-14
*b^6*r^4*y^2*a^2-20*r^6*a^4*b^2*y^2+35*r^6*b^4*y^2*a^2+6*y*a^8*b^4*y0-
54*y*a^6*b^6*y0+56*y*b^8*r^4*y0-12*b^4*y*a^4*x^2*y0*r^2+30*b^4*y*a^6*x
^2*y0-2*b^6*y^3*a^2*y0*x^2-48*b^6*y*a^4*y0*x^2+32*b^6*y*a^2*x^2*y0*r^2
-82*b^4*y^2*a^4*x^2*r^2-10*b^8*x^2*y*y0*r^2-70*b^6*y^2*x^2*a^2*r^2+8*a
^4*y^3*b^4*y0*r^2-34*y*a^6*b^4*y0*r^2+170*y*a^4*b^6*y0*r^2+10*y*a^4*r^
4*b^4*y0-66*y*a^2*b^6*r^4*y0-152*y*b^8*a^2*y0*r^2+77*a^6*y^2*b^4*r^2+5
8*b^8*y^5*y0+48*y*a^4*b^8*y0-52*a^4*y^4*b^6+2*b^4*y^9*y0+4*b^6*y^7*y0+
6*b^4*y^7*y0*x^2-12*a^6*y^4*b^4+12*x^2*a^6*y^4*r^2+10*x^2*y^3*b^4*r^4*
y0-12*x^2*a^8*y^2*r^2+4*r^8*a^4*y^2-7*r^8*b^2*y^2*a^2+2*y^6*a^4*r^4-2*
y^3*b^4*r^6*y0-4*y^6*a^4*x^2*r^2-6*b^4*y^7*y0*r^2+11*b^2*y^6*x^4*a^2+4
*b^2*y^8*x^2*a^2+4*b^2*y^8*a^2*r^2-13*b^2*y^6*r^4*a^2+6*b^4*y^5*r^4*y0
+5*r^8*a^4*x^2-4*r^8*b^2*x^2*a^2-16*a^8*b^2*x^2*r^2+9*a^8*x^4*r^2-9*a^
8*r^4*x^2-5*a^4*x^8*r^2+10*a^4*x^6*r^4+b^2*r^2*a^2*x^8-4*b^2*r^4*a^2*x
^6+2*x^6*b^4*y^3*y0+3*b^2*x^8*a^2*y^2+6*r^6*b^2*x^4*a^2-10*r^6*a^4*x^4
+6*x^4*b^4*y^5*y0+10*b^2*x^6*y^4*a^2-12*x^4*b^2*y^2*r^4*a^2+24*x^4*a^4
*y^2*r^4-5*y^4*b^2*x^4*a^2*r^2-10*x^4*b^4*y^3*y0*r^2-6*a^6*x^2*y^2*r^4
-15*y^4*a^4*x^4*r^2-2*x^6*b^4*y*y0*r^2-2*b^2*r^2*a^2*x^6*y^2+4*r^4*x^4
*b^4*y*y0+15*x^2*y^4*a^4*r^4-20*x^2*y^4*b^2*r^4*a^2-16*x^2*r^6*a^4*y^2
+18*x^2*r^6*b^2*y^2*a^2+12*x^2*y*b^6*r^4*y0+6*a^6*x^4*r^2*y^2-16*a^4*x
^6*r^2*y^2-14*b^4*y^5*x^2*y0*r^2+2*b^2*y^6*x^2*a^2*r^2-6*y^4*a^6*r^4-5
*y^4*r^6*a^4+15*y^4*r^6*b^2*a^2+6*a^8*y^2*r^4+2*a^6*r^6*y^2+2*y*b^4*r^
6*a^2*y0-2*y*b^6*r^6*y0+20*y^6*a^4*b^4+4*a^10*x^2*r^2-2*r^6*x^2*b^4*y*
y0-2*a^10*r^4+4*a^10*b^2*r^2+3*a^8*r^6+r^10*a^2*b^2-r^10*a^4);}{%
\maplemultiline{
\mathit{g12} := b^{2}( - 3\,y^{6}\,b^{4}\,a^{2}\,x^{2} - 9\,a^{4}\,r^{6}\,b^{4} + 2\,a^{6}\,b^{8} + 2\,a^{8}\,r^{4}\,b^{2} - 
19\,a^{6}\,r^{2}\,b^{6} + 32\,a^{6}\,b^{4}\,r^{4} \\
\mytaba - 13\,a^{8}\,b^{4}\,r^{2} + 20\,a^{4}\,b^{8}\,r^{2} - 10\,a^{4}\,b^{6}\,r^{4} + 8\,a^{8}\,b^{6} + 24\,y^{2}\,b^{10}\,a^{2
} - 15\,a^{6}\,r^{6}\,b^{2} \\
\mytaba + 10\,a^{6}\,x^{6}\,b^{2} - 17\,a^{6}\,x^{4}\,y^{2}\,b^{2} - 39\,y^{2}\,b^{8}\,x^{2}\,a^{2} + 23\,a^{6}\,y^{2}\,x^{2}\,b^{
4} - 16\,a^{8}\,x^{2}\,y^{2}\,b^{2} \\
\mytaba - 37\,a^{4}\,y^{4}\,x^{2}\,b^{4} + 24\,x^{2}\,y^{4}\,a^{6}\,b^{2} - 16\,b^{10}\,y^{3}\,\mathit{y0} + 2\,b^{4}\,x^{6}\,a^{4
} + 94\,y^{4}\,b^{6}\,x^{2}\,a^{2} \\
\mytaba + 9\,y^{2}\,b^{6}\,x^{4}\,a^{2} - 52\,b^{6}\,x^{2}\,a^{4}\,y^{2} + 48\,b^{4}\,x^{4}\,a^{4}\,y^{2} - 8\,b^{10}\,a^{4} - 18
\,a^{6}\,b^{4}\,x^{4} + 8\,b^{10}\,r^{2}\,a^{2} \\
\mytaba - 22\,b^{8}\,r^{4}\,a^{2} + 21\,b^{6}\,r^{6}\,a^{2} + 16\,b^{10}\,y\,\mathit{y0}\,r^{2} - 8\,r^{8}\,a^{2}\,b^{4} + 8\,r^{
8}\,a^{4}\,b^{2} - 16\,b^{6}\,x^{4}\,a^{4} \\
\mytaba + 21\,b^{8}\,x^{2}\,a^{4} + 2\,y^{3}\,b^{4}\,r^{2}\,a^{2}\,\mathit{y0}\,x^{2} - 20\,y^{4}\,a^{4}\,b^{2}\,x^{2}\,r^{2} - 2
\,a^{10}\,b^{4} - 2\,a^{10}\,x^{4} \\
\mytaba + 4\,a^{10}\,b^{2}\,x^{2} + a^{4}\,x^{10} - 3\,a^{8}\,x^{6} - 19\,a^{8}\,b^{4}\,x^{2} + 14\,a^{8}\,x^{4}\,b^{2} - 16\,a^{4
}\,b^{2}\,y^{6}\,x^{2} \\
\mytaba - 2\,a^{6}\,x^{6}\,y^{2} + 5\,y^{4}\,a^{4}\,x^{6} + 4\,a^{4}\,x^{8}\,y^{2} - 6\,x^{4}\,a^{6}\,y^{4} + 6\,x^{4}\,a^{8}\,y
^{2} + 2\,y^{6}\,a^{4}\,x^{4} + 6\,b^{6}\,x^{2}\,a^{6} \\
\mytaba - 14\,y^{4}\,b^{8}\,a^{2} + 24\,b^{6}\,y^{6}\,a^{2} + 20\,a^{6}\,y^{2}\,b^{6} - 4\,x^{4}\,b^{6}\,y\,\mathit{y0}\,r^{2} - 
8\,b^{4}\,r^{2}\,a^{2}\,x^{4}\,y^{2} \\
\mytaba - 8\,a^{4}\,x^{6}\,b^{2}\,r^{2} + 24\,a^{4}\,x^{4}\,b^{2}\,r^{4} + b^{4}\,r^{2}\,a^{2}\,x^{6} - 10\,b^{4}\,r^{4}\,a^{2}\,x
^{4} + 4\,x^{4}\,b^{6}\,y^{3}\,\mathit{y0} \\
\mytaba + 3\,b^{4}\,x^{6}\,a^{2}\,y^{2} + 4\,b^{4}\,a^{4}\,x^{2}\,r^{4} + 40\,r^{4}\,a^{6}\,x^{2}\,b^{2} - 10\,r^{2}\,a^{4}\,x^{2
}\,b^{6} + 17\,r^{6}\,b^{4}\,x^{2}\,a^{2} \\
\mytaba - 24\,r^{6}\,a^{4}\,x^{2}\,b^{2} - 20\,a^{4}\,x^{4}\,b^{2}\,r^{2}\,y^{2} - 30\,x^{2}\,b^{4}\,y^{2}\,r^{4}\,a^{2} - 2\,x^{2
}\,a^{2}\,y\,r^{4}\,b^{4}\,\mathit{y0} \\
\mytaba + 40\,x^{2}\,a^{4}\,y^{2}\,r^{4}\,b^{2} - 12\,b^{6}\,x^{2}\,r^{4}\,a^{2} + 8\,x^{2}\,b^{6}\,y^{5}\,\mathit{y0} + 6\,b^{4}
\,x^{4}\,y^{4}\,a^{2} + 3\,b^{6}\,x^{4}\,a^{2}\,r^{2} \\
\mytaba - 35\,a^{6}\,x^{4}\,b^{2}\,r^{2} + 16\,y^{4}\,b^{4}\,x^{2}\,a^{2}\,r^{2} + 116\,y^{3}\,a^{2}\,b^{6}\,\mathit{y0}\,r^{2} - 
4\,y^{5}\,b^{4}\,a^{2}\,\mathit{y0}\,r^{2} \\
\mytaba + 2\,y^{3}\,a^{2}\,r^{4}\,b^{4}\,\mathit{y0} + 34\,x^{2}\,a^{6}\,y^{2}\,b^{2}\,r^{2} - 20\,x^{2}\,b^{6}\,y^{3}\,\mathit{
y0}\,r^{2} + 10\,b^{8}\,x^{2}\,\mathit{y0}\,y^{3} \\
\mytaba - 13\,b^{8}\,x^{2}\,r^{2}\,a^{2} - 50\,y^{5}\,b^{6}\,\mathit{y0}\,a^{2} - 10\,y^{5}\,b^{6}\,\mathit{y0}\,r^{2} + 39\,y
^{6}\,b^{4}\,a^{2}\,r^{2} - 8\,y^{4}\,a^{4}\,b^{2}\,x^{4} \\
\mytaba - 67\,y^{4}\,a^{4}\,b^{4}\,r^{2} - 31\,y^{4}\,b^{6}\,a^{2}\,r^{2} - 8\,y^{3}\,a^{6}\,b^{4}\,\mathit{y0} + 8\,y^{3}\,b^{6}
\,r^{4}\,\mathit{y0} + 24\,y^{4}\,a^{6}\,r^{2}\,b^{2} \\
\mytaba - 16\,y^{6}\,a^{4}\,r^{2}\,b^{2} + 28\,y^{4}\,a^{4}\,r^{4}\,b^{2} - 56\,y^{4}\,b^{4}\,r^{4}\,a^{2} + 3\,b^{4}\,a^{4}\,x^{4
}\,r^{2} + 6\,b^{4}\,a^{6}\,x^{2}\,r^{2} \\
\mytaba + 4\,a^{8}\,y^{2}\,b^{4} - 10\,y^{8}\,b^{4}\,a^{2} + 100\,a^{4}\,y^{3}\,b^{6}\,\mathit{y0} - 17\,a^{6}\,y^{2}\,b^{2}\,r^{
4} - 16\,a^{8}\,y^{2}\,b^{2}\,r^{2} \\
\mytaba - 178\,a^{4}\,b^{6}\,y^{2}\,r^{2} + 48\,a^{4}\,b^{4}\,r^{4}\,y^{2} - 94\,b^{8}\,y^{3}\,\mathit{y0}\,a^{2} - 114\,b^{8}\,y
^{3}\,\mathit{y0}\,r^{2} + 32\,b^{8}\,y^{2}\,a^{2}\,r^{2} \\
\mytaba - 14\,b^{6}\,r^{4}\,y^{2}\,a^{2} - 20\,r^{6}\,a^{4}\,b^{2}\,y^{2} + 35\,r^{6}\,b^{4}\,y^{2}\,a^{2} + 6\,y\,a^{8}\,b^{4}\,
\mathit{y0} - 54\,y\,a^{6}\,b^{6}\,\mathit{y0} \\
\mytaba + 56\,y\,b^{8}\,r^{4}\,\mathit{y0} - 12\,b^{4}\,y\,a^{4}\,x^{2}\,\mathit{y0}\,r^{2} + 30\,b^{4}\,y\,a^{6}\,x^{2}\,
\mathit{y0} - 2\,b^{6}\,y^{3}\,a^{2}\,\mathit{y0}\,x^{2} \\
\mytaba - 48\,b^{6}\,y\,a^{4}\,\mathit{y0}\,x^{2} + 32\,b^{6}\,y\,a^{2}\,x^{2}\,\mathit{y0}\,r^{2} - 82\,b^{4}\,y^{2}\,a^{4}\,x^{
2}\,r^{2} - 10\,b^{8}\,x^{2}\,y\,\mathit{y0}\,r^{2} \\
\mytaba - 70\,b^{6}\,y^{2}\,x^{2}\,a^{2}\,r^{2} + 8\,a^{4}\,y^{3}\,b^{4}\,\mathit{y0}\,r^{2} - 34\,y\,a^{6}\,b^{4}\,\mathit{y0}\,r
^{2} + 170\,y\,a^{4}\,b^{6}\,\mathit{y0}\,r^{2} \\
\mytaba + 10\,y\,a^{4}\,r^{4}\,b^{4}\,\mathit{y0} - 66\,y\,a^{2}\,b^{6}\,r^{4}\,\mathit{y0} - 152\,y\,b^{8}\,a^{2}\,\mathit{y0}\,
r^{2} + 77\,a^{6}\,y^{2}\,b^{4}\,r^{2} \\
\mytaba + 58\,b^{8}\,y^{5}\,\mathit{y0} + 48\,y\,a^{4}\,b^{8}\,\mathit{y0} - 52\,a^{4}\,y^{4}\,b^{6} + 2\,b^{4}\,y^{9}\,\mathit{
y0} + 4\,b^{6}\,y^{7}\,\mathit{y0} + 6\,b^{4}\,y^{7}\,\mathit{y0}\,x^{2} \\
\mytaba - 12\,a^{6}\,y^{4}\,b^{4} + 12\,x^{2}\,a^{6}\,y^{4}\,r^{2} + 10\,x^{2}\,y^{3}\,b^{4}\,r^{4}\,\mathit{y0} - 12\,x^{2}\,a^{8
}\,y^{2}\,r^{2} + 4\,r^{8}\,a^{4}\,y^{2} \\
\mytaba - 7\,r^{8}\,b^{2}\,y^{2}\,a^{2} + 2\,y^{6}\,a^{4}\,r^{4} - 2\,y^{3}\,b^{4}\,r^{6}\,\mathit{y0} - 4\,y^{6}\,a^{4}\,x^{2}\,
r^{2} - 6\,b^{4}\,y^{7}\,\mathit{y0}\,r^{2} \\
\mytaba + 11\,b^{2}\,y^{6}\,x^{4}\,a^{2} + 4\,b^{2}\,y^{8}\,x^{2}\,a^{2} + 4\,b^{2}\,y^{8}\,a^{2}\,r^{2} - 13\,b^{2}\,y^{6}\,r^{4}
\,a^{2} + 6\,b^{4}\,y^{5}\,r^{4}\,\mathit{y0} \\
\mytaba + 5\,r^{8}\,a^{4}\,x^{2} - 4\,r^{8}\,b^{2}\,x^{2}\,a^{2} - 16\,a^{8}\,b^{2}\,x^{2}\,r^{2} + 9\,a^{8}\,x^{4}\,r^{2} - 9\,a
^{8}\,r^{4}\,x^{2} - 5\,a^{4}\,x^{8}\,r^{2} \\
\mytaba + 10\,a^{4}\,x^{6}\,r^{4} + b^{2}\,r^{2}\,a^{2}\,x^{8} - 4\,b^{2}\,r^{4}\,a^{2}\,x^{6} + 2\,x^{6}\,b^{4}\,y^{3}\,\mathit{
y0} + 3\,b^{2}\,x^{8}\,a^{2}\,y^{2} \\
\mytaba + 6\,r^{6}\,b^{2}\,x^{4}\,a^{2} - 10\,r^{6}\,a^{4}\,x^{4} + 6\,x^{4}\,b^{4}\,y^{5}\,\mathit{y0} + 10\,b^{2}\,x^{6}\,y^{4}
\,a^{2} - 12\,x^{4}\,b^{2}\,y^{2}\,r^{4}\,a^{2} \\
\mytaba + 24\,x^{4}\,a^{4}\,y^{2}\,r^{4} - 5\,y^{4}\,b^{2}\,x^{4}\,a^{2}\,r^{2} - 10\,x^{4}\,b^{4}\,y^{3}\,\mathit{y0}\,r^{2} - 6
\,a^{6}\,x^{2}\,y^{2}\,r^{4} - 15\,y^{4}\,a^{4}\,x^{4}\,r^{2} \\
\mytaba - 2\,x^{6}\,b^{4}\,y\,\mathit{y0}\,r^{2} - 2\,b^{2}\,r^{2}\,a^{2}\,x^{6}\,y^{2} + 4\,r^{4}\,x^{4}\,b^{4}\,y\,\mathit{y0}
 + 15\,x^{2}\,y^{4}\,a^{4}\,r^{4} \\
\mytaba - 20\,x^{2}\,y^{4}\,b^{2}\,r^{4}\,a^{2} - 16\,x^{2}\,r^{6}\,a^{4}\,y^{2} + 18\,x^{2}\,r^{6}\,b^{2}\,y^{2}\,a^{2} + 12\,x^{
2}\,y\,b^{6}\,r^{4}\,\mathit{y0} \\
\mytaba + 6\,a^{6}\,x^{4}\,r^{2}\,y^{2} - 16\,a^{4}\,x^{6}\,r^{2}\,y^{2} - 14\,b^{4}\,y^{5}\,x^{2}\,\mathit{y0}\,r^{2} + 2\,b^{2}
\,y^{6}\,x^{2}\,a^{2}\,r^{2} - 6\,y^{4}\,a^{6}\,r^{4} \\
\mytaba - 5\,y^{4}\,r^{6}\,a^{4} + 15\,y^{4}\,r^{6}\,b^{2}\,a^{2} + 6\,a^{8}\,y^{2}\,r^{4} + 2\,a^{6}\,r^{6}\,y^{2} + 2\,y\,b^{4}
\,r^{6}\,a^{2}\,\mathit{y0} - 2\,y\,b^{6}\,r^{6}\,\mathit{y0} \\
\mytaba + 20\,y^{6}\,a^{4}\,b^{4} + 4\,a^{10}\,x^{2}\,r^{2} - 2\,r^{6}\,x^{2}\,b^{4}\,y\,\mathit{y0} - 2\,a^{10}\,r^{4} + 4\,a^{10
}\,b^{2}\,r^{2} + 3\,a^{8}\,r^{6} \\
\mytaba + r^{10}\,a^{2}\,b^{2} - r^{10}\,a^{4}) }
}
\end{maplelatex}

\end{smallmaplegroup}
\begin{smallmaplegroup}
\begin{mapleinput}
\mapleinline{active}{1d}{g13:=factor(reducedG[13]);}{%
}
\end{mapleinput}

\mapleresult
\begin{maplelatex}
\mapleinline{inert}{2d}{g13 := y^2-2*y*y0+y0^2+x^2-2*x*x0+x0^2-r^2;}{%
\[
\mathit{g13} := y^{2} - 2\,y\,\mathit{y0} + \mathit{y0}^{2} + x^{2} - 2\,x\,\mathit{x0} + \mathit{x0}^{2} - r^{2}
\]
}
\end{maplelatex}

\end{smallmaplegroup}
\begin{smallmaplegroup}
\begin{mapleinput}
\mapleinline{active}{1d}{g14:=factor(reducedG[14]);}{%
}
\end{mapleinput}

\mapleresult
\begin{maplelatex}
\mapleinline{inert}{2d}{g14 := b^2*y0*x-b^2*y0*x0-a^2*x0*y+a^2*x0*y0;}{%
\[
\mathit{g14} := b^{2}\,\mathit{y0}\,x - b^{2}\,\mathit{y0}\,\mathit{x0} - a^{2}\,\mathit{x0}\,y + a^{2}\,\mathit{x0}\,\mathit{y0}
\]
}
\end{maplelatex}

\end{smallmaplegroup}
\begin{smallmaplegroup}
\begin{mapleinput}
\mapleinline{active}{1d}{g15:=factor(reducedG[15]);}{%
}
\end{mapleinput}

\mapleresult
\begin{maplelatex}
\mapleinline{inert}{2d}{g15 := b^2*(x0*y0*x+y*y0^2-y0*x^2-y0*y^2-y0*a^2+y0*r^2+y*a^2);}{%
\[
\mathit{g15} := b^{2}\,(\mathit{x0}\,\mathit{y0}\,x + y\,\mathit{y0}^{2} - \mathit{y0}\,x^{2} - \mathit{y0}\,y^{2} - \mathit{y0}\,
a^{2} + \mathit{y0}\,r^{2} + y\,a^{2})
\]
}
\end{maplelatex}

\end{smallmaplegroup}

\section{Appendix: Reduced \Grobner basis for the ideal\\ $I_S = \langle h_1, h_2, h_3, h_4, h_5, g \rangle \subset\BR[x, y, a, b, r]$}
\label{AppendC}
In order to solve system (\ref{eq:Fellipse}), we have found a reduced \Grobner basis for $I_S$ using the total term order {\tt tdeg} for 
$y > x > a > b > r$ with ties broken by inverse lexicographical order. We used package {\tt Groebner} from {\tt Maple 8} supplemented by a small custom package {\tt RJGrobner} written by the Authors. \cite{RJGrobner} The reduced basis contains fourteen homogeneous polynomials with degrees ranging from $10$ through $15$ and the number of monomials terms ranging from $51$ through $112.$ We will not display these polynomials here. However, here is the factorization the fourth polynomial $b_4$ in the basis:
\beq
\begin{split}
b_4 &= 
b^{2}\,x(4\,y^{2}\,r^{4}\,a^{2} - 5\,y^{4}\,a^{2}\,r^{2} + a^{4}\,r^{2}\,y^{2} - 2\,b^{2}\,x^{2}\,r^{4} + b^{2}\,x^{4}\,r^{2} - 9
\,x^{2}\,b^{2}\,y^{2}\,a^{2} - 11\,b^{2}\,y^{2}\,a^{2}\,r^{2} 
\\
& \mytabb - 2\,x^{2}\,b^{4}\,r^{2} + 4\,a^{6}\,r^{2} - 2\,a^{8} - b^{6}\,a^{2} - a^{4}\,b^{4} + 4\,a^{6}\,b^{2} + 6\,y^{2}\,a^{6} - 
a^{4}\,r^{4} - r^{6}\,a^{2} 
\\
& \mytabb - 4\,y^{4}\,b^{4} - 6\,y^{4}\,a^{4} - 4\,a^{6}\,x^{2} + 2\,y^{6}\,a^{2} - a^{4}\,x^{4} - b^{6}\,y^{2} + b^{6}\,r^{2} + r^{6}\,b^{2} - 4\,b^{4}\,r^{4} 
\\
& \mytabb + 2\,a^{4}\,b^{2}\,x^{2} - 9\,a^{4}\,b^{2}\,y^{2} + 2\,x^{2}\,b^{4}\,y^{2} + 5\,b^{4}\,y^{2}\,a^{2} - a^{4}\,x^{2}\,y^{2}
 - 3\,b^{2}\,a^{2}\,x^{4} + 6\,b^{2}\,y^{4}\,a^{2} 
\\
& \mytabb + 3\,x^{2}\,b^{4}\,a^{2} - y^{6}\,b^{2} + a^{2}\,x^{6} - 2\,b^{2}\,x^{2}\,a^{2}\,r^{2} + 2\,a^{4}\,x^{2}\,r^{2} - 3\,a^{2}
\,x^{4}\,r^{2} + 3\,a^{2}\,x^{2}\,r^{4} 
\\
& \mytabb - 8\,x^{2}\,y^{2}\,a^{2}\,r^{2} + 4\,b^{2}\,x^{2}\,y^{2}\,r^{2} + 8\,b^{4}\,y^{2}\,r^{2} + 5\,b^{2}\,r^{4}\,a^{2} + 5\,b
^{4}\,r^{2}\,a^{2} - 3\,b^{2}\,r^{4}\,y^{2} 
\\
& \mytabb + 3\,b^{2}\,y^{4}\,r^{2} - 10\,b^{2}\,a^{4}\,r^{2} - b^{2}\,x^{4}\,y^{2} - 2\,y^{4}\,b^{2}\,x^{2} + 4\,a^{2}\,x^{4}\,y^{2}
 + 5\,a^{2}\,y^{4}\,x^{2})
\end{split}
\eeq
It appears, that since $b>0,$ one possible solution is $x=0.$ Upon substituting $x=0$ into the list $F$ shown in (\ref{eq:Fellipse}), one gets
\beq
\begin{split}
&[0,  - y\,A (2\,b^{2}\,y^{4} + b^{4}\,y^{2} - 4\,b^{2}\,r^{2}\,y^{2} - 4\,b^{2}\,y^{2}\,a^{2} + y^{2}\,a^{2}\,r^{2} - b^{4}\,a^{2} - b^{4}\,r^{2} 
      + 2\,b^{2}\,r^{4} + 2\,a^{4}\,b^{2} - a^{4}\,r^{2} - r^{4}\,a^{2}),
\\  
& \mytabb - A( - 2\,b^{2}\,y^{4} + r^{2}\,y^{4} - b^{4}\,y^{2} + 2\,b^{2}\,r^{2}\,y^{2} + 4\,b^{2}\,y^{2}\,a^{2} - 3\,y^{2}\,a^{2}\,r^{2} 
                 - 2\,r^{4}\,y^{2} + b^{4}\,a^{2} - b^{4}\,r^{2} - 2\,a^{4}\,b^{2} 
\\
& \mytabb \mytabb + 2\,b^{2}\,a^{2}\,r^{2} + r^{6} + 2\,a^{4}\,r^{2} - 3\,r^{4}\,a^{2}),  
\\
& \mytabb (y + r - a)\,(y + r + a)\,(a - y + r)\,( - a - y + r)\,( - y^{2} + r^{2} + a^{2} - 2\,b^{2})\,A,
\\ 
& \mytabb A( - y^{4}\,a^{2} + 2\,b^{2}\,y^{4} + b^{4}\,y^{2} - 4\,b^{2}\,r^{2}\,y^{2} - 2\,b^{2}\,y^{2}\,a^{2} + 3\,y^{2}\,a^{2}\,r^{2} 
               + 2\,a^{4}\,y^{2} + b^{4}\,a^{2} - b^{4}\,r^{2} - 2\,b^{2}\,a^{2}\,r^{2} 
\\
& \mytabb \mytabb + 2\,b^{2}\,r^{4} + 3\,a^{4}\,r^{2} - 2\,r^{4}\,a^{2} - a^{6}),  
\\
& \mytabb (y + r - a)\,(y + r + a)\,(a - y + r)\,( - a - y + r)\,A^{2}]
\label{eq:sympF}
\end{split}
\eeq
where
\beq
A=  - y^{2}\,b^{2} - b^{4} - r^{2}\,a^{2} + b^{2}\,a^{2} + b^{2}\,r^{2}.
\eeq
It can be easily verified that setting $y = -r - a$ or $y = r+a$ in (\ref{eq:sympF}) does not yield any critical values for $r,$ as the system becomes, respectively,
\begin{align}
F_1 &= [0, \, - 2\,(r + a)\,(b^{2} + r\,a)^{4}\,r\,a, \,2\,(b^{2} + r\,a)^{4}\,r\,(r + a), \,0, \, - 2\,(b^{2} + r\,a)^{4}\,a\,(r + a), \,0],
\\[0.5ex]
F_4 &= [0, \,2\,(r + a)\,(b^{2} + r\,a)^{4}\,r\,a, \,2\,(b^{2} + r\,a)^{4}\,r\,(r + a), \,0, \, - 2\,(b^{2} + r\,a)^{4}\,a\,(r + a), \,0]
\end{align} 
as $a>b>r>0$ and is insolvable. However, setting $y=-r+a$ or $y=r-a$ does yield $r_{crit} = \frac{b^2}{a}$ and makes the system solvable in each case:
\begin{align}
F_2 &= [0, \,2\,(r - a)\,( - b^{2} + r\,a)^{4}\,r\,a, \,2\,( - b^{2} + r\,a)^{4}\,r\,(r - a), \,0, \,2\,( - b^{2} + r\,a)^{4}\,a\,(r - a), \,0],
\\[0.5ex]
F_3 &= [0, \, - 2\,(r - a)\,( - b^{2} + r\,a)^{4}\,r\,a, \,2\,( - b^{2} + r\,a)^{4}\,r\,(r - a), \,0, \,2\,( - b^{2} + r\,a)^{4}\,a\,(r - a), \,0]
\end{align}
Another possibility for solving the system is to set
\beq
y^2 = \frac{-b^4-r^2 a^2+b^2 a^2+b^2 r^2}{b^2}
\label{eq:SP12}
\eeq
Thus, 
\beq
\left(x=0,\,y = \pm \sqrt{\frac{-b^4-r^2 a^2+b^2 a^2+b^2 r^2}{b^2}}\right)
\eeq 
provide two solutions to the system without any condition on $r$ except that $r>0.$ These are the two virtual singular points found earlier that are always present for $a>b>r>0.$ Observe also that the critical value of $r$ satisfies 
$$
r_{crit} = \frac{b^2}{a} = \left(\frac{b}{a}\right) b < b \quad \mbox{since} \quad a> b.
$$
Another possibility to solve the system is to set
\beq
y^2 = r^{2} + a^{2} - 2\,b^{2}
\eeq
When we substitute the above expressions for $y^2$ into $F,$ we get
\beq
\begin{split}
& \left[0, \, - 6\,y\,( - b^{2} + r\,a)^{2}\,(b^{2} + r\,a)^{2}\,b^{2},\,6\,( - b^{2} + r\,a)^{2}\,(b^{2} + r\,a)^{2}\,( - b + r)\,(b + r), \,0,
\right.  
\\
& \mytaba 6\,( - b^{2} + r\,a)^{2}\,(b^{2} + r\,a)^{2}\,(b - a)\,(b + a),  
\\
& \mytaba \left. (a + r - y)\,(a + r + y)\,(r - a + y)\,( - a + r - y)\,( - b^{2} + r\,a)^{2}\,(b^{2} + r\,a)^{2} \right] 
\end{split}
\eeq
The above shows that in this case in order to satisfy the equations of the variety of singular points $\bV(I_S)$ we must also have 
$r = r_{crit} = \frac{b^2}{a}.$ Hence we get
\beq
\left(x=0,\, y = \pm \left(\frac{a^2-b^2}{a}\right) = \frac{c^2}{a}\right) \quad \mbox{when} \quad r = \frac{b^2}{a},
\eeq
which is the special case of the above singular points discussed above.

Finally, we have special cases that we are not interested in but we mention them for completeness. They result when we set $y=0$ in the list $F$ (in addition to $x=0.$ The first case is when $a=b,$ that is, when the ellipse is a circle. The second case is when $r=b$ or $r=a.$

Thus, to summarize, we have these two cases of singular points when $x=0:$
\begin{enumerate}
\item $\left( x=0,\,y^2 = -\frac{-b^2 r^2 +b^4+ a^2 r^2 - a^2 b^2}{b^2} \right)$ and $r$ is arbitrary as long as $b > r > 0.$ This yields singular points that are always present.
\item $\left( x=0,\, y = \pm (-r +a) = \pm \left(\frac{c^2}{a}\right) \right)$ when $r = r_{crit} = \frac{b^2}{a}$ and $c^2 = a^2 - b^2.$ This yields singular points in the special case $r = r_{crit}.$
\end{enumerate}

In order to complete the analysis of all solutions of the system $F,$ we need to consider the case when $x,y \neq 0.$ The fourteen-polynomial reduced \Grobner basis found for $I_S$ when solving system (\ref{eq:Fellipse}), could be simplified by diving each polynomial by $x, y, a^2,$ or $b^2,$ as appropriate. Let $I \subset \BR[x,y,a,b,r]$ denote the ideal generated by these fourteen simplified polynomials. In an effort to eliminate variable 
$y$ from some polynomials, we computed\footnote{We used {\tt gbasis} procedure from the {\tt Groebner} package with the total order {\tt lexdeg([y],[r,a,b,x])} indicating that the variable $y$ needed to be eliminated.} the \Grobner basis for the first elimination ideal 
$I_x = I \cap \BR[x,a,b,r]$ for the lexicographic order $y > r > a > b > x.$ The first polynomial in this new thirteen-polynomial homogeneous basis for $I_x$ is of degree $12.$ It is the square of the following polynomial $g_x \subset \BR[x,a,b,r]$ of degree $6$ in $x:$
\beq
\begin{split}
g_x &= 
(b^{6} - 3\,b^{4}\,a^{2} + 3\,a^{4}\,b^{2} - a^{6})\,x^{6} + (6\,b^{6}\,a^{2} + 3\,a^{6}\,r^{2} - 6\,b^{2}\,a^{4}\,r^{2} 
       + 3\,b^{4}\,r^{2}\,a^{2} - 3\,b^{8} - 3\,a^{4}\,b^{4})\,x^{4} 
\\
& \mytaba  + ( - 21\,b^{4}\,a^{4}\,r^{2} + 3\,b^{10} - 3\,b^{8}\,a^{2} + 21\,b^{6}\,r^{2}\,a^{2} + 3\,a^{4}\,b^{2}\,r^{4} - 3\,r^{4}\,a^{6})\,x^{2}
\\
& \mytaba \mytaba - 3\,b^{4}\,a^{4}\,r^{4} - b^{12} + 3\,b^{8}\,r^{2}\,a^{2} + a^{6}\,r^{6}
\end{split}
\eeq
Upon substituting $z=x^2$ this becomes an easily solvable cubic equation for $z$ that has the following solution:
\beq
z  = \frac {3\,b^{2}\,a^{2}\,r^{2} - r^{2}\,a^{2}\,(a\,b\,r)^{\frac23} - 3\,b^{3}\,(a\,b\,r)^{\frac13}\,a\,r 
     + b^{4}\,(a\,b\,r)^{\frac23}}{(a\,b\,r)^{\frac23}\,(-a^{2} + b^{2})} = - \frac{\delta}{(a\,b\,r)^{\frac23}\,(-a^{2} + b^{2})}
\eeq
where $\delta$ is given in (\ref{eq:delta}). Thus, we have derived the $x$-coordinates of the singular points $S_5,S_6,S_7$ and $S_8$ listed in 
(\ref{eq:S5}), (\ref{eq:S6}), (\ref{eq:S7}), and (\ref{eq:S8}). 

To find the $y$ coordinates of these four singular points, we proceed in a similar fashion except that we eliminate the $x$ variable. Thus, compute the first elimination ideal $I_y = I \cap \BR[y,a,b,r]$ for the lexicographic order $x > r > a > b > y.$ 
The first polynomial in this new six-polynomial homogeneous basis for $I_y$ is of degree $12.$ It is the square of the following polynomial 
$g_y \subset \BR[y,a,b,r]$ of degree $6$ in $y:$
\beq
\begin{split}
g_y &= ( - b^{6} + 3\,b^{4}\,a^{2} - 3\,a^{4}\,b^{2} + a^{6})\,y^{6} + (3\,b^{2}\,a^{4}\,r^{2} + 3\,b^{6}\,r^{2} - 3\,a^{4}\,b^{4} 
        + 6\,a^{6}\,b^{2} - 3\,a^{8} - 6\,b^{4}\,r^{2}\,a^{2})\,y^{4} 
\\
& \mytaba + (3\,b^{4}\,a^{2}\,r^{4} - 3\,b^{6}\,r^{4} - 21\,b^{4}\,a^{4}\,r^{2} + 3\,a^{10} + 21\,a^{6}\,b^{2}\,r^{2} - 3\,a^{8}\,b^{2})\,y^{2} 
         + b^{6}\,r^{6} 
\\
& \mytaba \mytaba - 3\,b^{4}\,a^{4}\,r^{4} - a^{12} + 3\,a^{8}\,b^{2}\,r^{2} 
\end{split}
\eeq
Upon substituting $z=y^2$ this becomes an easily solvable cubic equation for $z$ that has the following solution:
\beq
z = \frac { - a^{4}\,(a\,b\,r)^{\frac23} + 3\,(a\,b\,r)^{\frac13}\,b\,r\,a^{3} - 3\,b^{2}\,a^{2}\,r^{2} 
            + b^{2}\,r^{2}\,(a\,b\,r)^{\frac23}}{(a\,b\,r)^{\frac23}\,(a^{2} - b^{2})^{2}} 
  = \frac{(-a^2+b^2)\,\varepsilon}{(a\,b\,r)^{2/3}\,(-a^2+b^2)^2}
\eeq
where $\varepsilon$ is given in (\ref{eq:epsilon}). Thus, we have derived the $y$-coordinates of the singular points $S_5,S_6,S_7$ and $S_8$ listed in 
(\ref{eq:S5}), (\ref{eq:S6}), (\ref{eq:S7}), and (\ref{eq:S8}).

\vskip+5pt
\noindent{\bf Acknowledgment:}

The first author, R.A., gratefully acknowledges financial support from the College of Arts and Sciences and the Department of Mathematics, Tennessee Technological University, to present a preliminary version of this paper at the 5th Annual Hawaii International Conference on Statistics, Mathematics and Related Fields, Honolulu, January 2006.


\vskip 2em
\noindent Submitted: \today; Revised: TBA.
\end{document}